\def\ArxivVersion{1}
\documentclass[twoside,11pt]{article}
\usepackage[utf8]{inputenc}
\usepackage{amsmath,amsthm}
\ifdefined\ArxivVersion
  \usepackage[preprint]{jmlr2e}
\else
  \usepackage{jmlr2e}
\fi
\usepackage{bm,mathrsfs,algorithmic,makecell,tablefootnote,comment,soul}
\usepackage[ruled,linesnumbered]{algorithm2e}
\usepackage{microtype,subcaption,cleveref,xcolor,dsfont,bbm,enumitem,blkarray}
\usepackage{adjustbox}
\hypersetup{hidelinks,pdftitle={Optimal Clustering of Multi-layer Networks under Mixture Block Models},pdfauthor={Zhongyuan Lyu, Ting Li, Dong Xia}}
\makeatletter
\def\@starteditor{}
\def\@endeditor{}
\makeatother
\definecolor{revisionred}{RGB}{190,0,0}
\ifdefined\MarkedVersion
  \DeclareRobustCommand{\ZL}[1]{{\color{revisionred}#1}}
\else
  \DeclareRobustCommand{\ZL}[1]{#1}
\fi
\newcommand{\jasarr}[1]{#1}
\newcommand{\jasaadd}[1]{#1}

\newenvironment{movedenv}{}{}
\colorlet{red}{black}\colorlet{blue}{black}
\colorlet{magenta}{black}\colorlet{cyan}{black}

\allowdisplaybreaks[1]
\SetAlFnt{\small}
\SetAlCapFnt{\small}
\SetAlCapNameFnt{\small}
\SetAlgoNlRelativeSize{-1}
\graphicspath{{./}}
\ShortHeadings{Optimal Clustering of Multi-layer Networks}{Lyu, Li and Xia}
\firstpageno{1}

\def\calA{{\mathcal A}}

\def\calE{{\mathcal E}}

\def\calG{{\mathcal G}}

\def\calI{{\mathcal I}}
\def\calJ{{\mathcal J}}

\def\calN{{\mathcal N}}

\def\calP{{\mathcal P}}

\def\calS{{\mathcal S}}

\def\calV{{\mathcal V}}

\def\bcalA{{\boldsymbol{\mathcal A}}}
\def\bcalB{{\boldsymbol{\mathcal B}}}
\def\bcalC{{\boldsymbol{\mathcal C}}}

\def\bcalE{{\boldsymbol{\mathcal E}}}

\def\bcalX{{\boldsymbol{\mathcal X}}}

\def\EE{{\mathbb E}}

\def\II{{\mathbb I}}

\def\PP{{\mathbb P}}

\def\RR{{\mathbb R}}

\def\ZZ{{\mathbb Z}}

\def\j{{\boldsymbol j}}
\def\k{{\boldsymbol k}}

\DeclareMathOperator*{\argmax}{arg\,max}

\def\calA{{\cal  A}}

\def\calE{{\cal  E}} 
 
\def\calG{{\cal  G}} 
 
\def\calI{{\cal  I}} 
\def\calJ{{\cal  J}}

 \def\scrM{{\mathscr  M}}
\def\calN{{\cal  N}} 
 
\def\calP{{\cal  P}}

\def\calS{{\cal  S}}

\def\calV{{\cal  V}}

\def\bone{\bfm 1}

\newcommand{\bfm}[1]{\ensuremath{\mathbf{#1}}}

\def\ba{\bfm a}   \def\bA{\bfm A}
\def\bb{\bfm b}   \def\bB{\bfm B}  
\def\bc{\bfm c}     
   \def\bD{\bfm D}  
\def\be{\bfm e}   \def\bE{\bfm E}  \def\EE{\mathbb{E}}
  \def\bF{\bfm F}  
\def\bg{\bfm g}     
   \def\bH{\bfm H}  
   \def\bI{\bfm I}  \def\II{\mathbb{I}}

   \def\bL{\bfm L}  
   \def\bM{\bfm M}  
     
   \def\bO{\bfm O}  
\def\bp{\bfm p}   \def\bP{\bfm P}  \def\PP{\mathbb{P}}
   \def\bQ{\bfm Q}  
   \def\bR{\bfm R}  \def\RR{\mathbb{R}}
   \def\bS{\bfm S}  
   \def\bT{\bfm T}  
\def\bu{\bfm u}   \def\bU{\bfm U}  
\def\bv{\bfm v}   \def\bV{\bfm V}  
\def\bw{\bfm w}   \def\bW{\bfm W}  
\def\bx{\bfm x}   \def\bX{\bfm X}  
\def\by{\bfm y}   \def\bY{\bfm Y}  
\def\bz{\bfm z}   \def\bZ{\bfm Z}  \def\ZZ{\mathbb{Z}}

\def\bsig{\boldsymbol{\sigma}}

\def\hat{\widehat}
\def\wt{\widetilde}

\newcommand\brac[1]{\left( #1 \right)}
\newcommand\fro[1]{\| #1 \|_{\rm F}}
\newcommand\op[1]{\left\|#1\right\|}

\def\bcalE{\boldsymbol{\mathcal{E}}}
\newcommand{\E}{\mathbb{E}}

\begin{document}
\title{Optimal Clustering of Multi-layer Networks under Mixture Block Models}
\author{\name Zhongyuan Lyu \email zhongyuan.lyu@sydney.edu.au \\
\addr The University of Sydney Business School\\
The University of Sydney, Sydney NSW 2008, Australia
\AND
\name Ting Li \email liting@sustech.edu.cn \\
\addr Department of Statistics and Data Science\\
Southern University of Science and Technology, Shenzhen, China
\AND
\name Dong Xia \email dxia@cuhk.edu.hk \\
\addr Department of Statistics and Data Science\\
The Chinese University of Hong Kong, Shatin, N.T., Hong Kong SAR, China}
\editor{}
\maketitle
\begingroup
\makeatletter
\long\def\@makefntext#1{\noindent #1}
\makeatother
\footnotetext[0]{Zhongyuan Lyu and Ting Li are co-first authors.}
\endgroup
\begin{abstract}
We study clustering for multi-layer networks, where each observation is an undirected network measured on the same set of vertices. Under the mixture multi-layer stochastic block model (MMSBM), each network is generated from one of two stochastic block models, and the goal is to recover the latent network labels. We establish the minimax optimal network clustering error rate, which takes an exponential form and is characterized by \ZL{the Chernoff information} between the edge probability distributions of the component networks. We propose a novel two-stage estimator that combines tensor-based initialization, node and sample splitting, and a likelihood-based refinement. The refinement step exploits estimated node communities and achieves the minimax rate \ZL{with high probability} under extreme network sparsity. We further extend the likelihood-refinement framework to count-weighted networks through a mixture multi-layer Poisson block model. Numerical studies and two real-data analyses illustrate the finite-sample behavior of the proposed method and its use for both external-label benchmarking and unlabeled exploratory network clustering.

\end{abstract}
\section{Introduction}\label{sec:intro}

\jasarr{Many modern scientific studies produce a collection of networks rather than a single network \citep{kivela2014multilayer,weber2021functional, li2021super}. In such data, each observational unit is itself a network: a subject-specific brain connectivity network \citep{li2021super, relion2019network,paul2020random}, a commodity-specific trade network among countries \citep{lyu2021latent,jing2021community}, or an ecosystem-specific interaction network among species \citep{poisot2016mangal}. These examples share a common statistical feature that the networks are measured on a common vertex set, while the scientific question concerns heterogeneity across the network samples. For instance, one may ask whether brain networks form latent groups related to disease status, whether commodities share similar international trading structures, or whether ecological systems exhibit recurring interaction patterns. Throughout this paper, we consider such node-aligned collections of networks, and we call each observed network sample a layer in the standard terminology of multi-layer network analysis. }

Most of existing literature focused on \emph{community detection} \citep{paul2016consistent, paul2020random, chen2022global,agterberg2022joint, jing2022community}, which aims at clustering ``similar" vertices into non-overlapping groups. When multiple networks are present, finding similarity have substantial scientific values in many fields \citep{jing2021community}. \jasarr{In this paper, we focus on the network-sample clustering. Given $n$ networks observed on the same $d$ vertices, the goal is to cluster the network samples according to their latent connectivity-generating mechanisms. This problem is different from community detection, which clusters vertices within one or several networks. Network-sample clustering treats the whole network as the observational unit. The within-network community structure is not the final target, but it is \ZL{used to estimate the edge probabilities in the refinement step}. Scientific applications often require this network-level view because labels for the networks may be unavailable or unreliable. Clustering a large collection of networks can provide representative network classes, guide downstream comparative analysis, and reveal structural heterogeneity that is difficult to see from individual edges or vertices.}

Broadly speaking, two types of network clustering methods have appeared in the literature. The first type consists of those graphon-based methods \citep{mukherjee2017clustering, sabanayagam2021graphon, shao2022higher}. These methods are non-parametric in nature, offering great practical flexibility in that network sizes can vary and node registration can be absent. The ultimate clustering is performed on a low-dimensional embedding of networks. Nonetheless, it is unclear what the ``optimal" embedding is in practice and the non-parametric framework lacks interpretability of the results of network clustering. The second line of methods are model-based, assuming that the observed networks are generated by some mixture model.  See, e.g., \cite{stanley2016clustering,signorelli2020model,mantziou2021bayesian, chabert2022learning, rebafka2023model}. Model-based approaches facilitate the design of specially crafted clustering algorithms and enjoy better interpretability. However, \jasarr{the earlier model-based works} provided few insights on the theoretical fronts of network clustering. More recently, \cite{jing2021community} proposed a mixture multi-layer stochastic block model (MMSBM), which simultaneously models node communities and network clusters. It is a generalization of SBM \citep{holland1983stochastic} - one of the most popular generative models for a single network. MMSBM allows heterogeneous node community assignments across different networks. A tensor-based clustering algorithm was designed in \cite{jing2021community}, which can consistently recover the node communities and network clusters even when the observed networks are extremely sparse. Some other works \cite{fan2021alma, pensky2021clustering} proposed alternative algorithms attempting to achieve smaller clustering error but with a more stringent network sparsity condition. \cite{chen2022global} studied a generative multi-layer network models and investigated the optimal error rates in node  community detection. \jasarr{Their proposed model does not include latent labels for each network and therefore does not address network-sample clustering.}

This paper aims at studying the fundamental limit of clustering \jasarr{network-valued observations under MMSBM}.  The question is challenging for several reasons. First of all, \ZL{our network clustering method uses} node community detection. Without taking advantage of node community information, multiple networks are barely matrix-valued observations with binary entries. Clustering becomes extremely challenging if only some low-rank structure is present, even for low-rank Gaussian mixture model \citep{lyu2022optimal}. \ZL{The proposed algorithm} thus hinges upon efficient node community detection. Secondly, multi-layer networks can be very sparse, for example, some network may even be a disconnected graph. Intuitively, it is statistically challenging to compare different networks if most of them are disconnected graphs. We shall design novel clustering methods and technical tools to cope with the extreme sparsity in multi-layer networks. Lastly, model-based clustering methods usually require specific design of computational algorithms for both initialization and node community detection. The integrated analysis of both computational and statistical performances is also technically challenging.

Our main theoretical result shows that the minimax optimal network clustering error under MMSBM \ZL{with fixed component parameters and a fixed balance constant $\alpha>1$} has the sharp exponential form \ZL{$\exp(-[1+o(1)]I^{\ast})$} as $I^{\ast}\to\infty$. Here $I^{\ast}$ represents the \ZL{Chernoff information} between the two Bernoulli random matrices characterizing the edge generating mechanism of two SBMs. The formal definition of $I^{\ast}$ under MMSBM can be found in Section~\ref{sec:MMSBM}. Thus $I^{\ast}$ plays the same conceptual role as a signal-to-noise or separation parameter in classical mixture clustering, but it is adapted to network-valued observations and aggregates the distinguishability of the two block-structured network distributions. We first prove that this exponent is achieved by an oracle likelihood rule when the stochastic block model parameters are known, and we then establish a matching minimax lower bound. This identifies the exact statistical benchmark for network-sample clustering under MMSBM. This \ZL{exponential rate improves upon the} polynomial-type error rates attained by prior works \citep{jing2021community,fan2021alma,pensky2021clustering}, thereby filling the void in the optimal error rate of network clustering under MMSBM.

\jasarr{The main algorithmic challenge is to attain this oracle benchmark without knowing the node communities, the block probability matrices, or the network labels.} To that end, we propose a novel two-stage clustering algorithm based on the popular \emph{initialization then refinement} scheme. The algorithm begins with initial node community detection and layer clustering, which are  utilized to estimate the edge probability matrices. The refinement step re-assigns layer label to networks based on the maximum likelihood. The accuracy of estimated edge probability matrices crucially affects the success of refinement procedure. We prove that the two-stage clustering algorithm can achieve the minimax optimal network clustering error rate with high probability and in expectation, as long as the initial node community detection and layer clustering are sufficiently accurate.

It thus remains to design an efficient initial node community detection and layer clustering algorithm under MMSBM. The \ZL{layer} clustering error rate delivered by \cite{jing2021community} involves an additional $\log n$ term, \ZL{which imposes a stronger separation requirement for refinement}. We design a novel initialization algorithm based on tensor decomposition and node and sample switching. The algorithm involves both node and network sample splitting and uses two sub-samples for cross spectral estimate. The sample splitting allows us to derive a sharp concentration inequality for the sum of random sparse matrices. Note that sample splitting introduces loss of efficiency only for initialization, which is recovered subsequently by the refinement procedure. By combining the initialization and two-stage clustering algorithms, we end up with a computationally efficient method which achieves \ZL{the sharp exponential error bound} and allows extreme network sparsity under MMSBM. Besides the theoretical achievements, numerical simulations and real data experiments demonstrate that our new network clustering algorithm considerably outperforms the method proposed in \cite{jing2021community}.

\jasarr{These results connect network clustering with a broader statistical literature on sharp theory in clustering literature. Exponential minimax rates have been established for community detection in stochastic block models \citep{zhang2016minimax,gao2017achieving}, for node community detection in multi-layer networks \citep{chen2022global}, and for Gaussian and low-rank mixture models \citep{lu2016statistical,loffler2019optimality,loffler2020computationally,gao2022iterative,lyu2022optimal}. The present paper provides a network-valued analogue of this theory for clustering network samples rather than vertices or Euclidean observations. In particular, it shows that \ZL{Chernoff information} is the intrinsic information measure governing optimal graph-level clustering under mixture block models.}

We also extend MMSBM to multi-layer weighted networks where edges can carry count-type weights. The model is referred to as the mixture multi-layer Poisson block model (MMPBM) in that the edge-wise weight follows a Poisson distribution. Poisson block model has been studied for weighted networks. See, e.g., \cite{sewell2016latent,chiquet2019variational}.
By analyzing the oracle case and a likelihood-based Lloyd's algorithm, we show that the minimax optimal rate of network clustering under MMPBM also takes the form \ZL{$\exp(-[1+o(1)]I^{\ast})$}, where $I^{\ast}$ is the \ZL{Chernoff information} between the Poisson distributions describing the edge weight distribution. A similar two-stage clustering algorithm is shown to achieve the minimax optimal clustering error rate  if well initialized.

Our contributions are summarized as follows. We study the clustering problem and investigate the fundamental limits for multi-layer binary networks and multi-layer Poisson networks.
By analyzing a likelihood-based Lloyd's algorithm and establishing respective minimax lower bounds, we demonstrate that the minimax optimal clustering error rate in these models takes a universal exponential form where the exponent reflects the \ZL{Chernoff information} between the underlying data distributions. We provide a general two-stage clustering algorithm based on the \emph{initialization-then-refinement} scheme which achieves the minimax optimal clustering error rate \ZL{with high probability}. Initialization under MMSBM is challenging. We design a novel initialization algorithm involving both node and sample switching and cross spectral estimate, which is guaranteed to provide a sufficiently accurate initialization even when the multi-layer networks are extremely sparse.

{The rest of the paper is organized as follows. Section~\ref{sec:MMSBM} introduces the MMSBM and its two-stage estimator. Section~\ref{sec:MMPBM} presents the weighted-network analogue and its two-stage estimator. Section~\ref{sec:num} reports two real-data analyses and points to the deferred numerical studies. The appendices contain the deferred MMPBM theory, the remaining numerical studies, and all proofs.}

The following common notations will frequently appear in later sections. Vectors, matrices, and tensors are denoted by bold lowercase letters, bold uppercase letters, and bold calligraphic letters, respectively. For a vector $\bx$, $x_j$ denotes its $j$-th coordinate; $\bA(i,j)$ and $\bcalA(i,j,k)$ denote matrix and tensor entries. We use $c$ and $C$ to denote generic small and large constants, which may vary from line to line. Their dependence on fixed model bounds is suppressed; any other dependence is stated explicitly.  For nonnegative sequences $x_n$  and $y_n$, we write $x_n\lesssim y_n$ or $y_n\gtrsim x_n$ or $x_n=O(y_n)$ or $y_n=\Omega(x_n)$  if there exists a constant $C>0$ such that $x_n\le Cy_n$, write $x_n\asymp y_n$ if $x_n\lesssim y_n$ and $y_n\lesssim x_n$, and write $x_n\ll y_n$ or $y_n\gg x_n$  or $x_n=o(y_n)$ or $y_n=\omega(x_n)$ if $x_n=O(c_ny_n)$  for some $c_n\rightarrow 0$. Denote $\op{\cdot}$  the $\ell_2$-norm for vectors and operator norms for matrices. Let $\|\cdot\|_{\rm F}$ denote the matrix Frobenius norm.

\section{\texorpdfstring{{Two-stage clustering for multi-layer binary networks}}{Two-stage clustering for mixture multi-layer networks}}\label{sec:MMSBM}

\subsection{Mixture multi-layer stochastic block model}

The mixture multi-layer stochastic block model (MMSBM) was introduced  in \cite{jing2021community} for community detection in multi-layer networks.  Suppose we observe a collection of undirected, binary networks $\{\calG_i\}_{i=1}^n$ of common $d$ nodes, or equivalently, their symmetric adjacency matrices $\{\bX_i\in\{0,1\}^{d\times d},i\in[n]\}$. The essential assumption is that each network sample, or {\it layer}, is sampled independently from {\emph{one}} of {\emph{two}} \footnote{\ZL{Our theory concerns two components; the numerical studies also consider more than two.}} stochastic block model (SBM). Let a vector $\bz^\ast\in[2]^n$ denote the latent label of each network, i.e., $z^{\ast}_i=2$ meaning that the $i$-th network was sampled from the second SBM. The two SBM's are characterized by two {\it community probability matrices} $\bB_{1},\bB_{2}\in(0,1)^{K\times K}$ and {\it community assignment} $\bsig_1, \bsig_2: [d]\mapsto [K]$. Here we assume there are $K$ communities in each SBM for simplicity. Each edge of the network is independently sampled as
\begin{align*}
	\bX_i(j_1,j_2)\sim \text{Bernoulli}\left(\bB_{z^\ast_i}[\bsig_{z^\ast_i}(j_1),\bsig_{z^\ast_i}(j_2)]\right ),\quad i\in[n],\ \ZL{1\le j_1\le j_2\le d}
\end{align*}
Note that we allow self-loops for notational simplicity, whereas there is no essential difficulty in extending our results to the scenario with no self-loops.

After observing all the networks $\calG_i, i\in[n]$, our goal is to cluster them into two groups, i.e., which layers are possibly generated from the same SBM.
The Hamming distance evaluates the goodness of clustering:
\begin{equation}\label{eq:network-hamming-loss}
	h(\bz,\bz^\ast)=\min_{\pi\in  \mathfrak{S}_2}\frac{1}{n}\sum_{i=1}^n\mathbb{I}[ z_i\ne \pi(z^\ast_i)],
\end{equation}
where $\bz$ is a vector in $[2]^n$ and $\mathfrak{S}_k$ stands for the set of all permutations of $[k]$.  It can be viewed as the mis-clustering error of multi-layer networks.  Previously, \cite{jing2021community} proposed a tensor-decomposition based approach for layer clustering and showed that it can consistently recover the layer labels. Our interest is to characterize the minimax optimal clustering error under MMSBM. For ease of exposition,  define $\calI_m(\bz^\ast):=\{i\in[n]:z^\ast_i=m\}$ and $\calS_k(\bsig_m):=\{j\in[d]:\bsig_m(j)=k\}$ as the corresponding index set of $\bz^\ast$ and $\bsig_m$ for $\forall m\in[2] ,k\in[K]$. Denote $n_m:=|\calI_m(\bz^{\ast})|$ so that $n_1+n_2=n$.

Throughout this section, we focus on the  following parameter space of MMSBM:
\begin{align}\label{eq:parameter-space}
	\Theta:=\Theta(n,d,K, \bp,\alpha,\beta,\gamma)=\bigg\{\left(\bz,\{\bB_1,\bB_2\},\{\bsig_1,\bsig_2\}\right):\bz\in[2]^n,|\calI_m(\bz)|\in\left[\frac{n}{2\alpha},\frac{\alpha n}{2}\right], \nonumber
	\\ \bsig_m:[d]\rightarrow[K],|\calS_k(\bsig_m)|\in\left[\frac{d}{\beta K},\frac{\beta d}{K}\right], \bB_m=\bB_m^\top\in(0,1)^{K\times K},\nonumber
	\\  p_m:=\min_{i} \bB_m(i,i)\ge   \gamma^{-1}\max_i\bB_m(i,i),\ \gamma \min_{i\ne j} \bB_m(i,j)\ge q_m:=\max_{i\ne j} \bB_m(i,j), \nonumber
	\\  \max\{p_{m_1}/p_{m_2}, q_{m_1}/q_{m_2}, p_{m_1}/q_{m_2},  q_{m_1}/p_{m_2}\}\leq \gamma,\ \forall m_1, m_2\in[2]
	\bigg\},
\end{align}
where $\alpha,\beta,\gamma\ge 1$ are assumed to be absolute constants and $\bp:=(p_1,q_1, p_2,q_2)$ is a vector containing the boundary probabilities. Denote $\bP_m(j_1,j_2):=\bB_m\left[\bsig_m(j_1),\bsig_m(j_2)\right ]$ the edge probability matrix, i.e., the expectation of adjacency matrix under SBM$\left(\bB_m,\bsig_m\right)$. Note that the network sparsity is characterized by the probabilities $p_m, q_m$ and we will be particularly interested in the extremely sparse case when $p_m, q_m$ are of order $[(n\wedge d)d]^{-1}$.

\begin{movedenv}
\ZL{
Let $\calJ_d:=\{(j_1,j_2): 1\le j_1\le j_2\le d\}$. Define the separation strength by the Chernoff information
\begin{equation}\label{eq:def-I-star-mmsbm}
I^*:
=\max_{0\le s\le1}\bigg\{-\sum_{\omega\in\calJ_d}\log\Big(\bP_1(\omega)^{1-s}\bP_2(\omega)^s+[1-\bP_1(\omega)]^{1-s}[1-\bP_2(\omega)]^s\Big)\bigg\}.
\end{equation}
The corresponding R\'enyi divergence of order $1/2$ \citep{van2014renyi} is
\begin{equation}\label{eq:def-I-half-mmsbm}
I_{1/2}:=-2\sum_{\omega\in\calJ_d}\log\Big(\sqrt{\bP_1(\omega)\bP_2(\omega)}+\sqrt{[1-\bP_1(\omega)][1-\bP_2(\omega)]}\Big).
\end{equation}
Then $I_{1/2}/2\le I^*\le I_{1/2}$.
}
\end{movedenv}

Since network clustering becomes trivial if the two probability matrices $\bP_{m_1}$ and $\bP_{m_2}$ are strikingly different, we focus on the more difficult yet challenging regime where the edge probabilities are homogeneous,  i.e.,  $p_{m_1}/p_{m_2}$ and $q_{m_1}/q_{m_2}$ are bounded for $\forall m_1, m_2\in[2]$.
Homogeneous probability condition as above is typical in existing literature of studying minimax rates in network analysis. See, e.g., \cite{zhang2020theoretical, zhou2018non, gao2017achieving, chen2022global,zhenstructure,fei2019achieving, zhang2023fundamental}.

\subsection{\texorpdfstring{\textcolor{black}{Two-stage network clustering}}{Two-stage network clustering}}\label{subsec:practical}
{Our estimator follows an initialization-then-refinement strategy. Starting from initial layer labels and local community memberships, we estimate the block probability matrices by pooling information across layers and within blocks, and then update each layer label by maximizing the Bernoulli likelihood.}

\begin{algorithm}[!htb]
	\caption{Two-stage network clustering under MMSBM}\label{algorithm}\label{algorithm:simplified}
	\KwData{Network samples $\{\bX_i\}_{i=1}^n$, any initialization algorithm, denoted by $\text{Init}(\cdot )$, for estimating network label $\bz^\ast$ and local community memberships  $\bsig_1,\bsig_2$ .}
	\KwResult{Network label $\hat\bz$}
	$(\widetilde \bz,\widetilde \bsig_1,\widetilde \bsig_2)=\text{Init}(\{\bX_{i}\}_{i=1}^n)$

	\For{$m\in\{1,2\}$}{
		\For{$k,l\in[K]$}{
			\ZL{
\begin{equation}\label{eq:B_kl-est}
\widehat\bB_m(k,l)=
\frac{\sum_{i:\widetilde z_i=m}\sum_{\omega\in\calJ_d}\bX_i(\omega)\mathbf1(\{\widetilde\bsig_m(j_1),\widetilde\bsig_m(j_2)\}=\{k,l\})}
{\sum_{i:\widetilde z_i=m}\sum_{\omega\in\calJ_d}\mathbf1(\{\widetilde\bsig_m(j_1),\widetilde\bsig_m(j_2)\}=\{k,l\})},\quad \omega=(j_1,j_2),
\end{equation}
}
		}
		}
	Let
	\begin{equation}\label{eq:z_i-est}
			\hat z_i=\argmax_{m\in \{1,2\}}\Bigg\{\sum_{j_1\le j_2}\bX_i(j_1,j_2)\log\frac{\hat\bP_{m}(j_1,j_2)}{1-\hat\bP_{m}(j_1,j_2)}+\sum_{j_1\le j_2}\log\left(1-\hat\bP_{m}(j_1,j_2)\right)\Bigg\}
	\end{equation}
		where $\hat\bP_{m}(j_1,j_2):=\hat \bB_m\left[\widetilde \bsig_m(j_1),\widetilde \bsig_m(j_2)\right]$
\end{algorithm}

\subsection{\texorpdfstring{\textcolor{black}{Provable tensor-based initialization}}{Provable tensor-based initialization}}\label{sec:MMSBM-init}

{We next describe a tensor-based initialization for layer clustering and local community detection, adapted from \cite{jing2021community}.}

\begin{movedenv}
To this end, some extra notations  are necessary. Denote $\bcalX\in\{0,1\}^{d\times d\times n}$ the {\it adjacency tensor} such that its $i$-th slice $\bcalX(:,:,i):=\bX_i$,  and let {$\bZ_1, \bZ_2\in\{0,1\}^{d\times K}$} be the local membership matrix of two SBMs, respectively. {More precisely,} $\bZ_m{(j,:)}=\be^\top_{\bsig_m(j)}$ for $m\in [2]$ where $\be_j$ represent the $j$-th canonical basis vector whose dimension varies at different appearances. Under MMSBM, it is readily seen that $\EE(\bX_i|z_i^\ast)=\bZ_{z_i^\ast}\bB_{z_i^\ast}\bZ^\top_{z_i^\ast}$ and the expected adjacency tensor is decomposable in the following format
\begin{align}\label{eq:tensor-representation}
	\E(\bcalX|\bz^\ast)=\bcalB\times_1 \bar\bZ\times_2\bar\bZ\times_3 \bW,
\end{align}
where $\bar\bZ:=(\bZ_1,\bZ_2)\in\{0,1\}^{d\times 2K}$ is called the {\it global membership} matrix, $\bW:=(\be_{z_1^\ast},\be_{z_2^\ast}\cdots,\be_{z_n^\ast})^\top\in\{0,1\}^{n\times 2}$ is the layer label matrix, and $\bcalB\in[0,1]^{2K\times 2K\times 2}$ is a probability tensor with its $1$st and $2$nd slices are defined as
\begin{align*}
	\bcalB(:,:,1):=\left [ \begin{matrix}
\bB_1& 0 \\
0& 0 \\
\end{matrix} \right ] \quad {\rm and}\quad
\bcalB(:,:,2):=\left [ \begin{matrix}
0& 0 \\
0& \bB_2 \\
\end{matrix} \right ].
\end{align*}
The multilinear product in \eqref{eq:tensor-representation} is defined in the way such that for any $i_1,i_2\in[d]$ and $i_3\in[n]$, $\E[\bcalX(i_1,i_2,i_3)|\bz^\ast]=\sum_{j_1=1}^{2K}\sum_{j_2=1}^{2K}\sum_{m\in\{1,2\}}\bcalB(j_1,j_2,m) \bar\bZ(i_1,j_1)\bar\bZ(i_2,j_2) \bW(i_3,m)$. For notational clarity, we use $\text{SVD}_r(\bM)$ to denote the top $r$ left singular vectors of any matrix $\bM$ and $\scrM_k(\bcalA)$ to denote the mode-$k$ matricization of any tensor $\bcalA$. See \cite{kolda2009tensor} for a more comprehensive introduction on tensor algebra.

Due to technical reasons, we focus on the regime $n,d\to\infty$ and $n=O(d)$.
Denote $\bcalB=:\bar p\bcalB_0$ with $\bar a:=(n\wedge d)d(p_1\vee p_2)$ and $\bar p:=\bar a(nd)^{-1}$. Under \eqref{eq:parameter-space}, $\bP_m(u,v)\asymp\bar p$ uniformly over $m\in[2]$ and $u,v\in[d]$. \ZL{Of particular interest is} the sparsity regime $\bar a\ll d$.
Let $\bar{\bZ}:=\bar\bU\bar\bD\bar\bR^\top$ be the compact SVD with $\bar\bU,\bar\bR$ being left/right singular vectors of $\bar\bZ$ and $\bar\bD\in\RR^{r\times r}$ containing all non-zero singular values of $\bar\bZ$, then the population adjacency tensor admits the {\it Tucker decomposition} as
\begin{align}\label{eq:tucker-representation}
	\E(\bcalX|\bz^\ast)=\bar\bcalC\times_1 \bar\bU\times_2\bar\bU\times_3 \bar\bW,
\end{align}

where $\bar\bcalC=\bcalB\times_1 \bar\bD \bar\bR^\top\times_2 \bar\bD \bar\bR^\top\times_3\bD_n^{1/2}\in \RR^{r\times r\times 2}$ is a core tensor, $\bD_{n}=\text{diag}(n_1,n_2)\in\RR^{2\times 2}$, $n_m=|\calI_m(\bz^\ast)|$ is the number of layers in $m$-th clusters and $\bar\bW:=\bW\bD_{n}^{-1/2}\in\RR^{n\times 2}$. The matrices $\bar\bU$ and $\bar\bW$ are often referred to as the singular vectors of $\E(\bcalX|\bz^\ast)$,  which provide the essential information of global community memberships and layer clusters.

The tensor-based spectral initialization algorithm proposed in \cite{jing2021community} works as follows. The mode-$1$ singular vectors $\bar \bU$ are estimated using the SVD of the sum of adjacency matrices $\sum_i \bX_i$. This simple yet useful method is blessed by the non-negativity of the adjacency matrices. The empirical mode-$1$ singular vectors are then truncated by the regularization operator defined as
 $\calP_\delta(\bU):=\text{SVD}_r(\bU_*)$ with $\bU_*(i,:):=\bU(i,:)\cdot \min\{\delta,\op{\bU(i,:)}\}/\op{\bU(i,:)}$ for any $\bU $. \ZL{By construction, $\|\bU_*\|_{2,\infty}\le\delta$.}  We assume that $\|\bar\bU\|_{2,\infty}=O\big((r/d)^{1/2}\big)$,  i.e.,  the majority rows of $\bar\bU$ have comparable norms.
 The adjacency tensor is then multiplied by the regularized estimate of mode-$1$ singular vectors, whose mode-$3$ singular vectors, denoted by $\hat\bW$, are used as the estimate of $\bar\bW$.
The layer clusters can be found by screening the rows of $\hat\bW$ using, say, $K$-means clustering with $K=2$.  Finally, the local community memberships are estimated by spectral clustering on aggregating networks with the same estimated layer labels.  The detailed implementation can be found in Algorithm~\ref{alg:spectral-init}.

\begin{algorithm}[!t]
\caption{Regularized spectral initial clustering (RSpec)} \label{alg:spectral-init}
\KwData{adjacency tensor $\bcalX\in\{0,1\}^{d\times d\times n}$, rank $r$, \ZL{number of communities $K$,} regularization parameters $\delta_1$}
\KwResult{layer clustering: $\tilde \bz$ and local community memberships $(\tilde\bsig_1, \tilde\bsig_2)$}

Initializations of $\bar\bU$ and $\bar\bW$: \;
{$\hat\bU=\text{SVD}_{r}\left(\sum_{i=1}^{n}\bX_i\right )$, $\hat\bW=\text{SVD}_2\left(\scrM_3\left(\bcalX\times_1 \calP_{\delta_1}\left(\hat\bU\right)^\top\times_2\calP_{\delta_1}\left(\hat\bU\right)^\top \right) \right)$}\;

Network clustering: $\tilde \bz \longleftarrow$  K-means clustering on rows of $\hat\bW$\\

Node clustering: $\tilde\bsig_m\longleftarrow$  K-means clustering on rows of ${\rm SVD}_{K}\big(\sum_{i=1}^n\II(\tilde z_i=m)\bX_i\big), m\in[2]$

\end{algorithm}

It has been shown by \cite{jing2021community} that Algorithm~\ref{alg:spectral-init} together with regularized tensor power iterations, referred to as TWIST, can provide a consistent layer clustering.
Unfortunately, the mis-clustering error rate established in \cite{jing2021community} is rather weak such that an immediate application of their result will trivialize our conclusions in Theorem~\ref{thm:main}.  In particular, \cite{jing2021community} showed that TWIST achieves the following layer mis-clustering error rate:
\begin{align*}
	h(\widetilde\bz^{\textsf{\tiny TWIST}},\bz^\ast)=O\Big(\frac{r^2\log(n\vee d)}{\bar a}\Big),
\end{align*}
\ZL{guaranteeing consistency when} $\bar{a}\gg r^2\log(n\vee d)$. The initialization condition by Theorem~\ref{thm:main} requires $\eta_z\cdot \bar a K^2(1+dn^{-1})\ll I^*$ means that
\begin{align}\label{eq:strong-cond}
	I^*\gg r^2K^2\Big(1+\frac{d}{n	}\Big)\log(n\vee d)
\end{align}
\ZL{The separation in \eqref{eq:strong-cond} makes the exponential error bound smaller than $1/n$, leading to exact recovery.}

We now present a simple variant of Algorithm~\ref{alg:spectral-init} that \ZL{removes the logarithmic factor from the error bound for initial layer clustering}.  The basic idea is to introduce independence between $\hat\bU$ and $\bcalX$ in Algorithm~\ref{alg:spectral-init} so that a sharper perturbation bound of $\hat\bW$ can be derived.  In particular, we randomly split the network samples $[n]$ and vertices $[d]$ into two disjoint subsets of approximately equal size, denoted by $\calN^{[0]}$, $\calN^{[1]}$ and $\calV^{[0]}, \calV^{[1]}$, respectively. Here $\calN^{[0]}\cup \calN^{[1]}=[n]$, $\calV^{[0]}\cup \calV^{[1]}=[d]$, and \ZL{$d_0:=|\calV^{[0]}|=\lfloor d/2\rfloor$ and $d_1:=|\calV^{[1]}|=\lceil d/2\rceil$.} We also write $\wt n_0:=|\calN^{[0]}|$ and $\wt n_1:=|\calN^{[1]}|$. We focus on the sub-networks restricted to each subset of vertices. As a result, we end up with four multi-layer sub-networks, whose adjacency tensors are denoted as $\bcalX_+^{[0]}\in \RR_+^{d_0\times d_0\times \wt n_0}$, $\bcalX_+^{[1]}\in\RR_+^{d_1\times d_1\times \wt n_1}$, $\bcalX_-^{[0]}\in\RR_+^{d_1\times d_1\times \wt n_0}$, and $\bcalX_-^{[1]}\in\RR_+^{d_0\times d_0\times \wt n_1}$,  respectively.  The node and layer splitting is displayed as in Figure~\ref{fig:node-layer-splitting}.

\begin{figure}
	\centering
		\includegraphics[width=.80\linewidth]{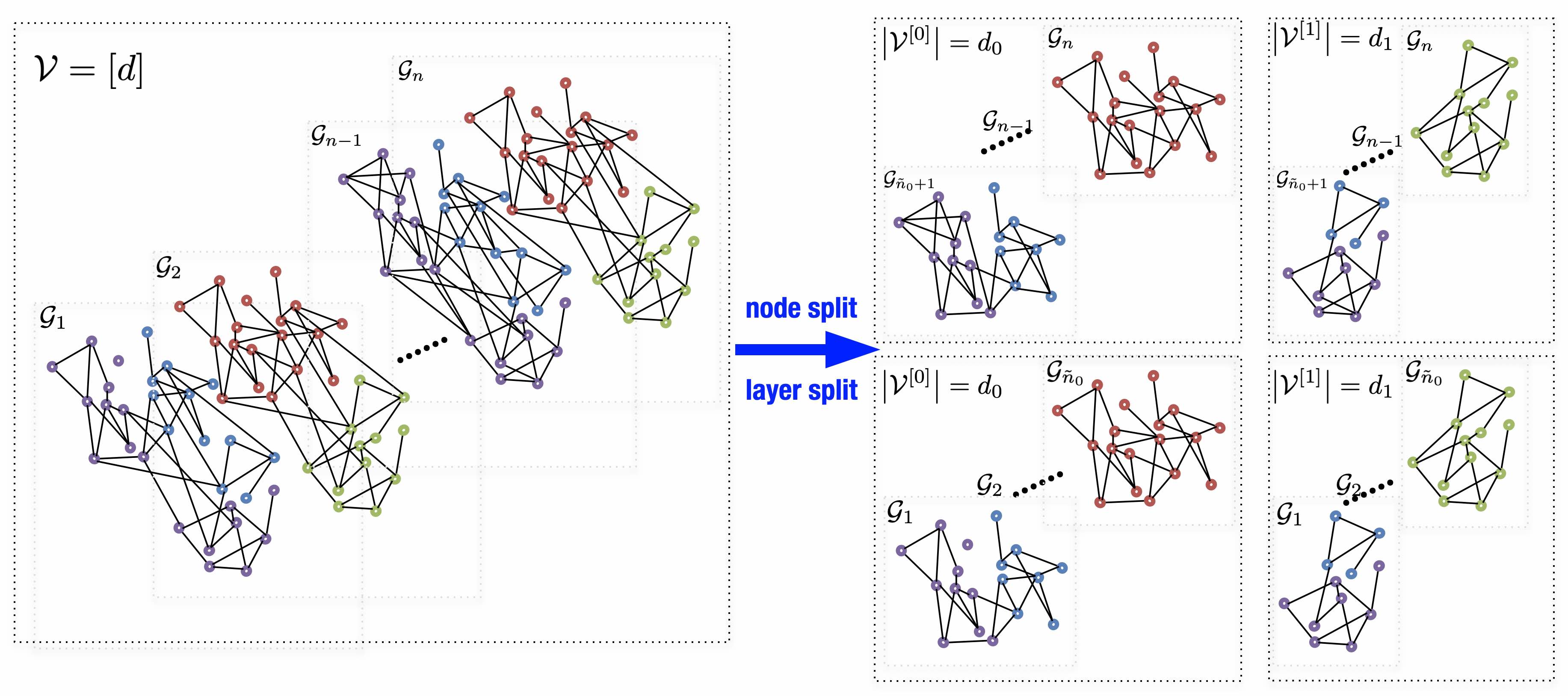}
	\caption{A diagram showing the node and layer sample splitting for initialization under MMSBM.  }
	\label{fig:node-layer-splitting}
\end{figure}

\begin{algorithm}[!htb]
\LinesNumberedHidden
\caption{Tensor-based initialization with node and sample switching}\label{algorithm:init}
\KwData{adjacency tensor $\bcalX\in\{0,1\}^{d\times d\times n}$, rank $r$, number of communities $K$, regularization parameters $\delta_1$.}
\KwResult{network labels $\widetilde\bz$, local memberships $\{\widetilde \bsig_1,\widetilde\bsig_2\}$.}
\textbf{Reference labels}: $\left(\check \bz,\left\{\check \bsig_1,\check \bsig_2\right\}\right )=\text{RSpec}(\bcalX,r,\ZL{K,}\delta_1)$\\
\textbf{Sample splitting}: Let $\bcalX_+^{[0]}\in\{0,1\}^{d_0\times d_0\times \wt n_0}$, $\bcalX_+^{[1]}\in\{0,1\}^{d_1\times d_1\times \wt n_1}$,  $\bcalX_-^{[0]}\in\{0,1\}^{d_1\times d_1\times \wt n_0}$, $\bcalX_-^{[1]}\in\{0,1\}^{d_0\times d_0\times \wt n_1}$ (or $\{\bX_i^{[+]}\in\{0,1\}^{d_k\times d_k}\}_{i\in\calN^{[k]}}$ and $\{\bX_i^{[-]}\in\{0,1\}^{d_{k^\prime}\times d_{k^\prime}}\}_{i\in\calN^{[k]}}$ for $k\ne k^\prime\in\{0,1\}$) be adjacency tensors (or samples) after sample splitting.\\
\textbf{Sample-switching estimation}: \For{$k\ne k^\prime\in\{0,1\}$}{
Initializations of $\bar\bU$ and $\bar\bW$: \;
{$\hat\bU^{[k^\prime]}=\text{SVD}_{r}\left(\sum_{i\in \calN^{[k^\prime]}}\ZL{\bX_i^{[-]}}\right )$, $\hat\bW^{[k]}=\text{SVD}_2\left(\scrM_3\left(\bcalX_+^{[k]}\times_1 \calP_{\delta_1}\left(\hat\bU^{[k^\prime]}\right)^\top\times_2\calP_{\delta_1}\left(\hat\bU^{[k^\prime]}\right)^\top \right) \right)$}.\;

Network clustering: $\widetilde\bz^{[k]}=$  K-means clustering on rows of $\hat\bW^{[k]}$\\
Local memberships recovery: \For{$m\in\{1,2\}$}{
	\quad $\widetilde\bsig_m^{[k^\prime]}$ =  K-means clustering on rows of $\text{SVD}_K\left(\sum_{i\in\calN^{[k]}:\widetilde z_i^{[k]}=m}\bX^{[-]}_i\right)$
}
}
\ZL{
\textbf{Network labels alignment}: \For{$k\in\{0,1\}$}{
$\pi_k\in\argmax_{\pi\in\mathfrak S_2}\sum_{i\in\calN^{[k]}}\II[\pi(\widetilde z_i^{[k]})=\check z_i],\qquad \widetilde z_i\leftarrow\pi_k(\widetilde z_i^{[k]}),\quad i\in\calN^{[k]}$\;
$(\widetilde\bsig_1^{[1-k]},\widetilde\bsig_2^{[1-k]})\leftarrow(\widetilde\bsig_{\pi_k^{-1}(1)}^{[1-k]},\widetilde\bsig_{\pi_k^{-1}(2)}^{[1-k]})$\;
}
\textbf{Local memberships alignment}: \For{$k\in\{0,1\},\,m\in[2]$}{
$\tau_{km}\in\argmax_{\tau\in\mathfrak S_K}\sum_{j\in\calV^{[k]}}\II[\tau(\widetilde\bsig_m^{[k]}(j))=\check\bsig_m(j)]$\;
$\widetilde\bsig_m(j)\leftarrow\tau_{km}(\widetilde\bsig_m^{[k]}(j)),\quad j\in\calV^{[k]}$\;
}
}
\end{algorithm}

For the full model and each restriction
$(\calN,\calV)=(\calN^{[k]},\calV^{[\ell]})$, $k,\ell\in\{0,1\}$,
use the substitutions
\begin{equation}\label{eq:init-restrictions}
\begin{split}
 (n,d,n_m,\bP_m,\bar\bZ)\ \mapsto\
 \Big(|\calN|,|\calV|,\sum_{i\in\calN}\mathbb I(z_i^*=m),\bP_m[\calV,\calV],\bar\bZ[\calV,:]\Big).
\end{split}
\end{equation}
The population factors are recomputed after splitting, with $\bar p$ and
$\bcalB_0$ unchanged. \ZL{Assume that the balance conditions in
\eqref{eq:parameter-space}, $\operatorname{rank}(\bar\bZ)=r$, and
$\|\bar\bU\|_{2,\infty}\le C\sqrt{r/d}$, for a fixed $C>0$, hold for all these models. }
\ZL{In the following theorem and corollary, probabilities and expectations
are conditional on the node and layer split, which is independent of
the observations and satisfies the stated assumptions.}
\ZL{
Let $r:=\operatorname{rank}(\bar\bZ)$ and define
\begin{align*}
	\overline{\bcalB}_0:=\bcalB_0\times_1\bar\bR^\top\times_2\bar\bR^\top,\qquad \sigma_{\min}(\overline{\bcalB}_0)
 :=\min\big\{\sigma_r(\scrM_1(\overline{\bcalB}_0)),
              \sigma_2(\scrM_3(\overline{\bcalB}_0))\big\}.
\end{align*}
}
\begin{theorem}\label{thm:init}\label{thm:initialization}
Suppose the following conditions (A1)--(A3) hold before and after the node and layer splitting \ZL{in \eqref{eq:init-restrictions}}:
\begin{enumerate}
\item[(A1)] \ZL{$\overline{\bcalB}_0$} has Tucker rank $(r,r,2)$.
\item[(A2)] For some $\kappa\ge1$, \ZL{$\operatorname{rank}(\bP_m)=K$ for $m\in[2]$ and}
\[
 \ZL{\max\bigg\{}\frac{\sigma_1(\bar\bD)}{\sigma_r(\bar\bD)}\ZL{,
       \max_{m\in[2]}\frac{\sigma_1(\bP_m)}{\sigma_K(\bP_m)}\bigg\}}\le\kappa.
\]
\item[(A3)] \ZL{For sufficiently large $C$ and sufficiently small $c>0$,}
\[
 \ZL{\sigma_r(n_1\bP_1+n_2\bP_2)\ge C\sqrt{rnd\bar p}},
 \qquad
 \frac{\sigma_1(n_1\bP_1+n_2\bP_2)}{\sigma_r(n_1\bP_1+n_2\bP_2)}
 \ZL{\le c\frac{\sqrt n}{\log n}}.
\]
\end{enumerate}
\ZL{For a sufficiently small constant $c>0$, assume also that}
\begin{equation}\label{eq:init-variable-gap}
\ZL{\frac{\kappa\sqrt K\,d\bar p}{\sigma_2^2(\scrM_3(\bar\bcalC))}
 \left[1+
 \frac{\kappa^2r^3(nd\bar p)\log^2(e\kappa r)\log^4(e(n\vee d))}
      {\sigma_r^2(n_1\bP_1+n_2\bP_2)}\right]\le c}
\end{equation}
\ZL{for the full model and under every substitution in
\eqref{eq:init-restrictions}.}
Furthermore, suppose
\[
 d\ge C_0\log^3(n\vee d),\qquad
 n/\log n\ge C_0\kappa^4r^6,
\]
\[
 \bar a\ge C_0\kappa^8r^4\log^2(\kappa r)\log^6(n\vee d).
\]
For sufficiently large $C_0$ and $C_1$, there exists a constant $C>0$ such that, with probability at least $1-C(n\vee d)^{-3}$, \ZL{there exists $\pi\in\mathfrak S_2$ such that} the output of Algorithm~\ref{algorithm:init} with $\delta_1=C_1\sqrt{r/d}$ satisfies
\begin{equation}\label{eq:init-joint-conclusion}
\begin{aligned}
 \frac1n\sum_{i=1}^n\mathbb I[\widetilde z_i\ne\pi(z_i^*)]
 &
 \ZL{\le C\max_{k\in\{0,1\}}
 \frac{d_k\bar p}
 {\sigma_2^2(\scrM_3(\mathbb E[\bcalX_+^{[k]}\mid\bz^*]))} =:\eta_z},\\
 \max_{m\in[2]}h(\widetilde\bsig_{\ZL{\pi(m)}},\bsig_m)
 &\ZL{=0}.
\end{aligned}
\end{equation}
\end{theorem}

The conditions \ZL{on the Tucker rank in (A1) and the bound on $\sigma_1(\bar\bD)/\sigma_r(\bar\bD)$ in (A2)} were originally proposed by \cite{jing2021community} in the context of MMSBM. \ZL{Condition (A2) also bounds the condition number of each component matrix.} Condition (A3) \ZL{controls the spectral gap of the aggregated adjacency matrix and} allows a diverging condition number of $n_1\bP_1+n_2\bP_2$. \ZL{Proposition~\ref{ir:split-preservation} shows that random splitting preserves the assumptions for the full model with high probability, up to changes in the constants.}

\ZL{Combining the initialization and refinement steps gives the following corollary.}

\begin{corollary}\label{cor:main}
Under the assumptions of Theorem~\ref{thm:init}, suppose
\ZL{$\max_{m\in[2],\omega\in\calJ_d}\bP_m(\omega)\le1-c_0$ for a fixed $c_0\in(0,1)$ and}
\begin{equation}\label{cond:Jmin-main}
 \ZL{\eta_z\to0,
 \qquad \frac{I^*}{K^2}\to\infty.}
\end{equation}
Then, for any $\epsilon\in(0,1)$, with probability at least
$1-\exp(-(I^*)^{1-\epsilon})\ZL{-C(n\vee d)^{-3}}$,
the output of Algorithm~\ref{algorithm}, initialized by
Algorithm~\ref{algorithm:init}, satisfies
\[
 h(\widehat\bz,\bz^*)\le\exp(-[1-o(1)]\ZL{I^*}).
\]
Moreover,
\[
 \mathbb E h(\widehat\bz,\bz^*)
 \le\exp(-[1-o(1)]\ZL{I^*})+\ZL{C(n\vee d)^{-3}}.
\]
\end{corollary}

\end{movedenv}

\subsection{\texorpdfstring{\textcolor{black}{Statistical guarantee for two-stage clustering}}{Statistical guarantee for two-stage clustering}}\label{subsec:theory-two-stage}

\ZL{We first show that accurate estimates of the edge probabilities preserve the sharp exponential error bound.} Denote the sup-norm $\|\bP\|_{\ell_\infty}:=\max_{\omega\in\calJ_d} |\bP(\omega)|$ for a symmetric matrix $\bP$.

\begin{lemma}\label{lem:optimal-rate-dev}
Assume \eqref{eq:parameter-space}.
\ZL{Assume $\max_{m\in[2],\omega\in\calJ_d}\bP_m(\omega)\le1-c_0$ for a fixed $c_0\in(0,1)$.}
Suppose that a network $\bX$ is sampled from the edge probability matrix $\bP_1$, i.e., $\bX(\omega)\stackrel{{\rm ind.}}{\sim} {\rm Bern}(\bP_1(\omega))$ for $\forall \omega\in\calJ_d$.  Let $\hat\bP_1$ and $\hat\bP_2$ be estimated edge probability matrices, possibly computed from the full sample.  Suppose that $\bP_1=\bZ_1\bB_1\bZ_1^{\top},  \hat \bP_1=\hat \bZ_1 \hat \bB_1\hat\bZ_1^{\top}$, where $\bZ_1,  \hat\bZ_1,  \bB_1,  \hat \bB_1$ satisfy the conditions defined in (\ref{eq:parameter-space}).
For $m\in[2]$, write $\PP_m$ for the joint data law when the target layer has component $m$. Let $\hat z$ be the MLE of $\bX$'s label based on $\hat\bP_1$ and $\hat \bP_2$:
\begin{equation}\label{eq:MLE-emp}
\hat z:=\underset{z\in\{1,2\}}{\arg\max}\ \sum_{\omega\in\calJ_d} \bX(\omega)\log \hat\bP_z(\omega)+\sum_{\omega\in\calJ_d}\big[1-\bX(\omega)\big]\log\big(1-\hat\bP_z(\omega)\big)
\end{equation}
If there exists an event $\calE$ such that, on $\calE$,
\[
{\max_{m\in\{1,2\}}\|\hat\bP_m-\bP_m\|_{\ell_\infty}\le \tilde\rho\cdot \frac{I^{\ast}}{d(d+1)/2}}
\]
{for some $\tilde\rho=o(1)$,} {where $I^{\ast}$ is defined in \eqref{eq:def-I-star-mmsbm}} and $I^{\ast}\to\infty$, then
$$
\max_{m\in[2]}\PP_m(\hat z\ne m,\calE)\leq \ZL{\exp\big(-\big[1-o(1)\big]I^*\big).}
$$

\end{lemma}

{The same bound holds if the likelihood comparison in \eqref{eq:MLE-emp} is perturbed by an additive $o(I^{\ast})$ threshold.}

Suppose that the initialization algorithm ensures
\begin{equation}\label{eq:init-cond}
\begin{split}
 \calA_0:=\Big\{\ZL{\exists\pi\in\mathfrak S_2:}\;&\ZL{\frac1n\sum_{i=1}^n\mathbb I[\widetilde z_i\ne\pi(z_i^*)]}\le\eta_z,\max_{m\in[2]}h(\widetilde\bsig_{\ZL{\pi(m)}},\bsig_m)\le\eta_\sigma\Big\},
\end{split}
\end{equation}
with $\PP(\calA_0)\ge 1-C n^{-2}$ for some absolute constant $C>0$. Here $\eta_z,\eta_\sigma\to0$ are deterministic as $n,d\to\infty$. 

The performance of both the initialization algorithm and Algorithm~\ref{algorithm} crucially relies on the network sparsity. We primarily focus on the extreme sparse case and set
\begin{equation}\label{eq:pmqm}
p_1=\frac{a_1}{(n\wedge d)d},\quad p_2=\frac{a_2}{(n\wedge d)d},\quad q_1=\frac{b_1}{(n\wedge d)d},\quad q_2=\frac{b_2}{(n\wedge d)d}.
\end{equation}
Denote $\bar a=a_1\vee a_2$ and $\underline a=a_1\wedge a_2$.

\begin{theorem}\label{thm:main}
Under \eqref{eq:parameter-space} and \eqref{eq:pmqm}, suppose that the initialization algorithm satisfies \eqref{eq:init-cond} with
\begin{equation}\label{eq:main-init-cond}
 \eta_\sigma\log\frac{1}{\eta_\sigma}
 =o\left(\frac{I^{*2}}{\bar a K^2(1+d/n)^2}\right),
 \qquad
 \eta_\sigma\vee\eta_z
 =o\left(\frac{I^*}{\bar a K^2(1+d/n)}\right),
\end{equation}
\ZL{In addition, suppose
\begin{equation}\label{eq:node-entropy}
 d\eta_\sigma\log(1/\eta_\sigma)=o(I^*),
\end{equation}}
If $\underline a\ge c>0$, $K^2\log n=o(n)$, \ZL{$\max_{m,\omega}\bP_m(\omega)\le1-c_0$ for a fixed $c_0\in(0,1)$}, and $I^*\to\infty$ as $n,d\to\infty$, then, for any $\epsilon\in(0,1)$, there exists an event $\calA_1$ with $\PP(\calA_1)\ge1-\exp(-(\ZL{I^*})^{1-\epsilon})$ such that, on $\calA_0\cap\calA_1$, the output of Algorithm~\ref{algorithm} satisfies
\[
 h(\widehat\bz,\bz^*)\le\exp(-[1-o(1)]\ZL{I^*}).
\]
Moreover,
\[
 \mathbb E h(\widehat\bz,\bz^*)\le \exp(-[1-o(1)]\ZL{I^*})
 +\ZL{O(n^{-2})}.
\]
\end{theorem}
Let us illustrate the implications of Theorem~\ref{thm:main} and the initialization requirement \eqref{eq:main-init-cond} for special cases. The case of particular interest is $p_m, q_m=o(1)$ so that
\begin{align*}
I^{\ast}\asymp \sum_{\omega\in\calJ_d} \frac{\big[\bP_1(\omega)-\bP_2(\omega)\big]^2}{\bP_1(\omega)\vee \bP_2(\omega)}\asymp& \frac{(p_1-p_2)^2}{p_1\vee p_2}\cdot \frac{d^2}{K} + \frac{(q_1-q_2)^2}{q_1\vee q_2}\cdot d^2\\
\asymp & \frac{(a_1-a_2)^2}{K(a_1\vee a_2)}\Big(1+\frac{d}{n}\big)+\frac{(b_1-b_2)^2}{b_1\vee b_2}\Big(1+\frac{d}{n}\Big),
\end{align*}
\ZL{when $\bsig_1=\bsig_2$ and $\bB_m=q_m\mathbf1\mathbf1^\top+(p_m-q_m)\bI_K$.}  Assume $K=O(1)$ for simplicity,  the condition $I^{\ast}\to\infty$ holds as long as
\begin{align}\label{eq:I_cond_simplify1}
\bigg[\frac{(a_1-a_2)^2}{a_1\vee a_2}+\frac{(b_1-b_2)^2}{b_1\vee b_2}\bigg]\bigg(1+\frac{d}{n}\bigg)\to \infty\quad {\rm as}\quad n, d\to\infty.
\end{align}
\ZL{On the high-probability event in Theorem~\ref{thm:main},} the mis-clustering error rate becomes
\begin{align}\label{eq:opt-rate-simplify1}
h(\widehat\bz,\bz^*)\leq \exp\bigg(  -\Omega\Big(\frac{(a_1-a_2)^2}{a_1\vee a_2}+\frac{(b_1-b_2)^2}{b_1\vee b_2}\Big)\cdot \Big(1+\frac{d}{n}\Big) \bigg)
\end{align}
The bound~\eqref{eq:opt-rate-simplify1} implies that networks can be more accurately clustered \ZL{as the network size $d$ increases relative to
the number of layers $n$}. Moreover, under \ZL{\eqref{eq:init-cond}, \eqref{eq:main-init-cond} and \eqref{eq:node-entropy}}, exact network clustering is possible if
\[
 \ZL{I^*>(1+\delta)\log n\qquad\text{for a fixed }\delta>0.}
\]

The initialization requirement in the aforementioned case is also inspiring. Note that $-x\log x= o(x^{1/(1+\delta)})$ as $x\to0$ for any constant $\delta>0$. For the sake of clarity, assume that $(b_1-b_2)^2/(b_1\vee b_2)\asymp (a_1-a_2)^2/(a_1\vee a_2)$. Then, \ZL{a sufficient condition for \eqref{eq:main-init-cond} is}
$$
\eta_{\sigma}=o\bigg(\frac{1}{\bar{a}^{1+\delta}}\cdot \Big[\frac{(a_1-a_2)^2}{\bar{a}}\Big]^{2(1+\delta)}\bigg),\qquad \eta_{\sigma},  \eta_z=o\bigg(\frac{1}{\bar{a}}\cdot \Big[\frac{(a_1-a_2)^2}{\bar{a}}\Big]\bigg),
$$
for any $\delta>0$.
They suggest that a stronger initialization condition \ZL{is sufficient} if $a_1-a_2$ is smaller, i.e., the two edge probability matrices are closer. Furthermore, if $a_1-a_2\asymp \bar{a}$, the network sparsity condition (\ref{eq:I_cond_simplify1}) \ZL{is satisfied if} $\bar{a}\to\infty$ and the initialization requirements reduce to $\eta_{\sigma},  \eta_z=o(1)$, namely, a consistent initial clustering \ZL{satisfies \eqref{eq:main-init-cond}. Condition~\eqref{eq:node-entropy} holds when $\eta_\sigma=0$.}

\subsection{\texorpdfstring{\textcolor{black}{Oracle property for likelihood-based refinement}}{Oracle property for likelihood-based refinement}}\label{subsec:oracle}
\begin{movedenv}
Lloyd's algorithm \citep{lloyd1982least} is a simple yet popular clustering algorithm. The algorithm alternates between re-estimating the cluster centers and re-assigning cluster labels. It has been demonstrated that Lloyd's algorithm can achieve statistically optimal clustering error under various mixture models such as GMM \citep{lu2016statistical} and the low-rank mixture model (LrMM, \cite{lyu2022optimal, lyu2023optimal}).

In this section, we explore the limit of Lloyd's algorithm under MMSBM by considering the oracle case. Under Gaussian mixture model, the oracle case means the situation where the population cluster centers are known beforehand and the label of a new observation can be decided by comparing its distances from the known population cluster centers. Similarly, for MMSBM, oracle case refers to the situation that the edge probability matrices $\bP_1$ and $\bP_2$ are already given. If an observed network $\bX\in\{0,1\}^{d\times d}$ is known sampled from either $\bP_1$ or $\bP_2$, the log-likelihood can be written as
\begin{align*}
	\ell(\bX|z):=\sum_{\omega\in\calJ_d}\bX(\omega)\log\bP_z(\omega)+\sum_{\omega\in\calJ_d}[1-\bX(\omega)]\log\left(1-\bP_z(\omega)\right)
\end{align*}
{where $z$ takes values in $\{1,2\}$.  The likelihood-based Lloyd's algorithm then assigns the label to $\bX$ which maximizes the log-likelihood $\ell(\bX|z)$,  i.e.,  $\hat z:=\arg\max_{z\in\{1,2\}}\ell (\bX| z)$.}

Without loss of generality,  assume the true label of $\bX$ is $1$.  \ZL{This network can be mis-clustered by Lloyd's algorithm only if}
\begin{equation}\label{eq-mis-event}
		\sum_{\omega\in\calJ_d}\bX(\omega)\log\frac{\bP_{2}(\omega)[1-\bP_{1}(\omega)]}{\bP_{1}(\omega)[1-\bP_{2}(\omega)]}\ZL{\ge}\sum_{\omega\in\calJ_d}\log\frac{1-\bP_{1}(\omega)}{1-\bP_{2}(\omega)}.
\end{equation}
By studying the probability of event defined in \eqref{eq-mis-event},  we get the oracle mis-clustering error rate of Lloyd's algorithm,  which is formally stated in the following lemma.  Note that the result holds for any given edge probability matrices,  which are not necessarily from stochastic block models.
\begin{lemma}\label{lem:optimal-rate-chernoff}
Suppose that $\bP_1$ and $\bP_2$ are known. For $m\in[2]$, let $\PP_m$ be the law under which
$\bX(\omega)\stackrel{\rm ind.}{\sim}{\rm Bern}(\bP_m(\omega))$, $\omega\in\calJ_d$.
For $\widehat z\in\argmax_{m\in[2]}\ell(\bX\mid m)$,
\[
 \max_{m\in[2]}\PP_m(\widehat z\ne m)\le\ZL{\exp\brac{-I^*}},
\]
where $I^*$ is defined in \eqref{eq:def-I-star-mmsbm}.
\end{lemma}

\ZL{The quantity $I_{1/2}$ in \eqref{eq:def-I-half-mmsbm} is the R\'enyi divergence of order $1/2$ between the two joint Bernoulli distributions.}  The R\'enyi divergence has commonly played the role of SNR in literature of network community detection \cite{zhang2016minimax,gao2017achieving}.  For instance,  the optimal mis-clustering rate for community detection in a two-community single SBM whose within-clustering and between-cluster edge probabilities are $p$ and $q$,  respectively,  is
$$
\exp\big(-[1-o(1)]dI_{0}/2\big),
$$
where $I_{0}:=-2\log \left(\sqrt{pq}+\sqrt{1-p}\sqrt{1-q}\right)$.  It implies that a consistent clustering requires SNR condition $dI_{0}\rightarrow \infty $.

\paragraph{Comparison with Gaussian mixture model}
For simplicity, consider the two-component Gaussian mixture model
$\frac{1}{2}N\big({\rm vec}(\bP_1),\sigma^{\ZL{2}}\bI\big)
+\frac{1}{2}N\big({\rm vec}(\bP_2),\sigma^{\ZL{2}}\bI\big)$,
where the population cluster centers are $\bP_1$ and $\bP_2$, respectively.
The difficulty of clustering is essentially characterized by the \emph{separation strength}:
\[
\Delta:=\fro{\bP_1-\bP_2}.
\]
It has been shown \citep{zhang2022leave,lu2016statistical,loffler2019optimality,lyu2022optimal}
that the minimax optimal mis-clustering error rate is
$\exp\big(-(1-o(1))\Delta^2/(8\sigma^2)\big)$.
In particular, this minimax rate holds for GMM \citep{loffler2019optimality},
sparse GMM (i.e., $\bP_1$ and $\bP_2$ are sparse, \cite{loffler2020computationally}),
and low-rank GMM (i.e., $\bP_1$ and $\bP_2$ are low-rank, \cite{lyu2022optimal}).
Furthermore, \cite{lyu2022optimal} \ZL{established an exponential error bound for a Lloyd's algorithm
under sub-Gaussian noise, with leading exponent $\Delta^2/(8\sigma_{\textsf{sg}}^2)$,
where $\sigma_{\textsf{sg}}^2$ is a variance proxy.}
Under MMSBM, the noise is Bernoulli and thus sub-Gaussian.
\ZL{If the edge probabilities are bounded away from zero and one,}
then $\sigma_{\textsf{sg}}\asymp1$ and the exponent
$\Delta^2/\sigma_{\textsf{sg}}^2\asymp\|\bP_1-\bP_2\|_{\rm F}^2\asymp I^*$.
\ZL{Thus the Gaussian and Bernoulli exponents have the same order in this regime.}
However, when $p_m,q_m=o(1)$ (the most interesting scenario for network analysis),
\ZL{the smallest common variance proxy satisfies
$\sigma_{\textsf{sg}}^2\asymp[\log(1/\bar p)]^{-1}$}
\citep{arbel2020strict}, so that
$\Delta^2/(8\sigma_{\textsf{sg}}^2)\asymp\ZL{\Delta^2\log(1/\bar p)}$.
The exponent in Lemma~\ref{lem:optimal-rate-chernoff} becomes
\[
 I^*\asymp\sum_{\omega\in\calJ_d}
 \frac{[\bP_1(\omega)-\bP_2(\omega)]^2}{\bP_1(\omega)\vee\bP_2(\omega)}
 \ZL{\asymp\frac{\Delta^2}{\bar p}\gg\Delta^2\log(1/\bar p)}.
\]
\ZL{For sparse networks, the sub-Gaussian bound has a smaller exponent. The low-rank Lloyd's algorithm uses Euclidean distances, whereas our refinement uses the Bernoulli likelihood.}

\end{movedenv}

\subsection{\texorpdfstring{\textcolor{black}{Minimax lower bound}}{Minimax lower bound}}\label{subsec:lower-bound}
\begin{movedenv}

We now characterize the minimax lower bound of clustering networks under MMSBM.
Recall the working parameter space $\Theta(n,d,K, \bp,\alpha,\beta,\gamma)$ and Hamming distance $h(\bz_1, \bz_2)$ defined in \eqref{eq:parameter-space} and \eqref{eq:network-hamming-loss},  respectively.  \ZL{We consider fixed} $(\bB_1,  \bB_2,  \bsig_1,  \bsig_2)$ that satisfies the constraint in \eqref{eq:parameter-space}.  Denote $\bB_{1,2}:=\{\bB_1,  \bB_2\}$,  $\bsig_{1,2}:=\{\bsig_1, \bsig_2\}$, and
$$
\Theta_{\bz}^{(\bB_{1,2},\bsig_{1,2})}
 :=\left\{\bz\in[2]^n:\frac{n}{2\alpha}\le\sum_{i=1}^n\mathbb I(z_i=m)
 \le\frac{\alpha n}{2},\ m\in[2]\right\}.
$$
The minimax lower bound in the following theorem presents an exponential rate, which is concordant with the oracle mis-clustering rate in Lemma~\ref{lem:optimal-rate-chernoff}.

\begin{theorem}\label{thm:minimax-lower-bound}
For fixed $\bB_{1,2}$ and $\bsig_{1,2}$ satisfying the constraints in \eqref{eq:parameter-space}, let $\bP_1=\bB_1\circ\bsig_1$ and $\bP_2=\bB_2\circ\bsig_2$, and let $I^*$ be defined in \eqref{eq:def-I-star-mmsbm}.
If $n\to\infty$, $I^*\to\infty$, and \ZL{$\alpha>1$ is fixed}, then
\[
 \inf_{\widehat\bz}\sup_{\bz^*\in\Theta_{\bz}^{(\bB_{1,2},\bsig_{1,2})}}
 \mathbb E h(\widehat\bz,\bz^*)
 \ge\exp\big(-[1+o(1)]\ZL{I^*}\big),
\]
where the infimum is taken over all possible clustering algorithms working on the multi-layer networks sampled from MMSBM with parameters $(\bz^*,\bB_{1,2},\bsig_{1,2})$.
\ZL{The lower bound holds even when the two component matrices are known.}
\end{theorem}
\ZL{The case of equal cluster sizes ($\alpha=1$) is treated separately in Proposition~\ref{prop:repair-exact-balance} in the appendix.}

\end{movedenv}

\section{\texorpdfstring{\textcolor{black}{Two-stage clustering for multi-layer count networks}}{Two-stage clustering for mixture multi-layer Poisson networks}}\label{sec:MMPBM}

\subsection{Mixture of Poisson Block Models}
Count-type data arise in many scientific and econometric applications such as,  4D-STEM imaging \citep{han2022optimal},  cell counts in spatial-transcriptome study \citep{moses2022museum},  and trade flow networks \citep{lyu2021latent, cai2022generalized}.  A Poisson block model (PBM) can be viewed as a generalization of SBM to count-type data.  We view PBM as a model for weighted networks.  A weighted network $\calG$ is observed with $d$ vertices,  each of which belongs to one of $K$ communities.  Let $\bsig:[d]\mapsto [K]$ denote the community assignment defined similarly as in SBM.   There is a symmetric matrix $\bB\in \RR_+^{K\times K}$ carrying the Poisson weight.  The network $\calG$ carries random weight,  or equivalently,  the entries of its adjacency matrix $\bX$,   generated from
$$
\bX(\omega)\stackrel{{\rm ind.}}{\sim}\ {\rm Poisson}\big(\bB[\bsig(j_1),\bsig(j_2)]\big),\quad \omega=(j_1,j_2)\in\calJ_d
$$
which is referred to as the PBM$(\bsig, \bB)$.

Suppose we observe a collection of weighted networks $\{\calG_i\}_{i=1}^n$ of common $d$ nodes, i.e., a multi-layer weighted network. Each network is either sampled from PBM$(\bsig_1, \bB_1)$ or PBM$(\bsig_2, \bB_2)$. Denote $\bz^{\ast}\in[2]^n$ the layer label vector meaning that $\calG_i\sim {\rm PBM}(\bsig_{z_i^{\ast}}, \bB_{z_i^{\ast}})$. Let $\bX_i$ denote the adjacency matrix of $\calG_i$ whose entries follow Poisson distributions as
\begin{equation}\label{eq:MMPBM}
\bX_i(\omega)\stackrel{{\rm ind.}}{\sim}\ {\rm Poisson}\big(\bB_{z^{\ast}_i}[\bsig_{z_i^{\ast}}(j_1),\bsig_{z_i^{\ast}}(j_2)]\big),\quad i\in[n],\quad\omega=(j_1,j_2)\in\calJ_d.
\end{equation}
We call it the mixture multi-layer Poisson block model (MMPBM), which models both layer clusters and nodal communities in multi-layer weighted networks.
After observing the networks $\calG_i, i\in[n]$, our goal is to find their latent labels.   The notations $\calI_m(\bz^{\ast}), \calS_k(\bsig_m)$, and $n_1,n_2$ are defined similarly as in Section~\ref{sec:MMSBM}.

Throughout this section, we focus on the  following parameter space of MMPBM:
\begin{align}\label{eq:parameter-space-MMPBM}
	\Theta:=\Theta(n,d,K, \theta_0,\alpha,\beta,\gamma)=\bigg\{\left(\bz,\{\bB_1,\bB_2\},\{\bsig_1,\bsig_2\}\right):\bz\in[2]^n,|\calI_m(\bz)|\in\left[\frac{n}{2\alpha},\frac{\alpha n}{2}\right], \nonumber
	\\ \bsig_m:[d]\rightarrow[K],|\calS_k(\bsig_m)|\in\left[\frac{d}{\beta K},\frac{\beta d}{K}\right], \bB_m=\bB_m^\top\in\RR_+^{K\times K},\nonumber
	\\  \min_{m\in[2]}\min_{i,j}\bB_m(i,j)\geq \theta_0 \quad {\rm and}\quad \max_{m\in[2]}\max_{i,j}\bB_m(i,j)\leq \gamma \theta_0
	\bigg\},
\end{align}
where $\alpha,\beta,\gamma\ge 1$ are assumed to be absolute constants and $\theta_0>0$ is a constant characterizing the overall edge intensity. The condition that $\theta_0$ is bounded away from zero is mild. See, \cite{han2022optimal} and \cite{cai2022generalized}.
Denote $\bP_m(j_1,j_2):=\bB_m\left[\bsig_m(j_1),\bsig_m(j_2)\right ]$ the edge intensity matrix, i.e., the expectation of adjacency matrix under PBM$\left(\bB_m,\bsig_m\right)$.

\begin{movedenv}
\ZL{
Define the separation strength under MMPBM by the Chernoff information
\begin{equation}\label{eq:def-I-star-mmpbm}
 I^*:=I^*(\bB_{1,2},\bsig_{1,2})=\max_{0\le s\le1}\sum_{\omega\in\calJ_d}
 \left[(1-s)\bP_1(\omega)+s\bP_2(\omega)-\bP_1(\omega)^{1-s}\bP_2(\omega)^s\right].
\end{equation}
The maximization uses a single $s$ for all edges. The R\'enyi divergence of order $1/2$ is
\begin{equation}\label{eq:def-I-half-mmpbm}
 I_{1/2}:=\sum_{\omega\in\calJ_d}\big[\sqrt{\bP_1(\omega)}-\sqrt{\bP_2(\omega)}\big]^2,
\end{equation}
and again $I_{1/2}/2\le I^*\le I_{1/2}$.
}
\end{movedenv}

\subsection{\texorpdfstring{\textcolor{black}{Two-stage clustering under MMPBM}}{Two-stage clustering under MMPBM}}
{For MMPBM, replace the Bernoulli likelihood in Algorithm~\ref{algorithm} with the Poisson likelihood.}

We first study the accuracy of estimated edge intensity matrices that ensures the optimality of the likelihood-based refinement.

\begin{lemma}\label{lem:MMPBM-optimal-rate-dev}
Assume \eqref{eq:parameter-space-MMPBM}.
Suppose that a network $\bX$ is sampled from the edge intensity matrix $\bP_1$, i.e., $\bX(\omega)\stackrel{{\rm ind.}}{\sim} {\rm Poisson}(\bP_1(\omega))$ for $\forall \omega\in\calJ_d$.  Let $\hat\bP_1$ and $\hat\bP_2$ be estimated edge intensity matrices, possibly computed from the full sample. For $m\in[2]$, write $\PP_m$ for the joint data law when the target layer has component $m$. Let $\hat z$ be the MLE of $\bX$'s label based on $\hat\bP_1$ and $\hat \bP_2$:
\begin{equation}\label{eq:MMPBM-MLE-emp}
\hat z:=\underset{z\in\{1,2\}}{\arg\max}\ \sum_{\omega\in\calJ_d} \bigg[-\hat \bP_z(\omega)+\bX(\omega)\log \hat \bP_z(\omega)\bigg]
\end{equation}
If there exists an event $\calE$ such that, on $\calE$,
\[
{\max_{m\in\{1,2\}}\|\hat\bP_m-\bP_m\|_{\ell_\infty}\le \tilde\rho\cdot \frac{I^{\ast}}{d(d+1)/2}}
\]
{for some $\tilde\rho=o(1)$,} {where $I^{\ast}$ is defined in \eqref{eq:def-I-star-mmpbm}} and $I^{\ast}\to\infty$, then
$$
\max_{m\in[2]}\PP_m(\hat z\ne m,\calE)\leq \ZL{\exp\big(-\big[1-o(1)\big]I^*\big).}
$$

\end{lemma}

\ZL{Define $\calA_0^{\rm P}$ as in \eqref{eq:init-cond} under MMPBM, and assume
\begin{equation}\label{eq:MMPBM-init-cond}
 {\color{black}\PP(\calA_0^{\rm P})\ge 1-C(n\vee d)^{-2}},
 \qquad \eta_z,\eta_\sigma\to0.
\end{equation}}
Here $\eta_z$ and $\eta_\sigma$ are deterministic.

\begin{algorithm}[htbp]
	\caption{Two-stage network clustering under MMPBM}\label{alg:MMPBM-two-stage}
	\KwData{network samples $\{\bX_i\}_{i=1}^n$, $K$-the number of communities in each network, $\text{Init}(\cdot )$-an initialization algorithm for estimating network labels and local community memberships. }
	\KwResult{network label $\hat\bz\in[2]^n$}
	Run Algorithm~\ref{algorithm} but replace \eqref{eq:z_i-est} with
		\begin{equation}\label{eq:MMPBM-z_i-est}
			\hat z_i=\argmax_{m\in [2]}\Bigg\{\sum_{\omega\in\calJ_d} \bigg[-\hat \bP_m(\omega)+\bX_i(\omega)\log \hat \bP_m(\omega)\bigg]\Bigg\}
		\end{equation}
\end{algorithm}

\begin{theorem}\label{thm:MMPBM-main}
Under \eqref{eq:parameter-space-MMPBM}, suppose that the initialization algorithm satisfies \eqref{eq:MMPBM-init-cond} with
\begin{equation}\label{eq:MMPBM-main-init-cond}
 \eta_\sigma\log\frac{1}{\eta_\sigma}
 =o\left(\frac{I^{*2}}{K^2d^2\theta_0(d/n)}\right),
 \qquad \eta_\sigma\vee\eta_z
 =o\left(\frac{I^*}{K^2d^2\theta_0}\right).
\end{equation}
If $\theta_0\ge c>0$, $K^2\log n=o(n)$, and $I^*\to\infty$ as $n,d\to\infty$, then, for any $\epsilon\in(0,1)$, there exists an event $\calA_1^{\rm P}$ with $\PP(\calA_1^{\rm P})\ge1-\exp(-(\ZL{I^*})^{1-\epsilon})$ such that, on $\calA_0^{\rm P}\cap\calA_1^{\rm P}$, the output of Algorithm~\ref{alg:MMPBM-two-stage} satisfies
\[
 h(\widehat\bz,\bz^*)\le\exp(-[1-o(1)]\ZL{I^*}).
\]
Moreover,
\[
 \mathbb E h(\widehat\bz,\bz^*)\le \exp(-[1-o(1)]\ZL{I^*})
 +\ZL{O((n\vee d)^{-2})}.
\]
\end{theorem}

\noindent\textit{Deferred theoretical results.} The oracle analysis and minimax lower bound for the MMPBM are deferred to the appendices for completeness.

\section{Numerical experiments}\label{sec:num}

\subsection{Simulation studies}

\noindent\textit{Deferred simulation results.} The simulation studies are deferred to the appendices for completeness.

\subsection{Real-data analysis}
\jasarr{We present two real-data analyses with different roles. The COBRE brain-connectivity data provide an externally labeled benchmark for unsupervised graph-level clustering, where diagnosis labels are withheld from the algorithms and used only for post-hoc evaluation. The food-trading data provide an unlabeled exploratory example in which commodity-specific networks are clustered and interpreted through their estimated network structures.}

We test our proposed method using the COBRE dataset \citep{aine2017multimodal, paul2020random}, which consists of resting-state fMRI experiments conducted on both diagnosed schizophrenia patients and healthy control subjects. Here the brain network is defined by the connectivity among regions of interest in a subject's brain. We focus on functional connectivity, a statistical measure of pairwise dependence between brain regions derived from fMRI data. With class labels available, \cite{paul2020random} further examined the differences in community structure between schizophrenia patients and healthy controls.

The processed connectomic dataset is downloaded from \cite{relion2019network}. It contains connectivity networks of 54 schizophrenia patients and 70 healthy controls, each represented as a $263\times 263$ matrix. Weak connections are converted to 0 by thresholding at the median of the mean network across all samples, and the remaining strong connections are converted to 1.

\jasarr{We compare the proposed refinement method with three unsupervised baselines. TWIST denotes the tensor-based multi-layer network clustering method of \cite{jing2021community}; Tucker denotes HOSVD/Tucker decomposition followed by clustering of the layer-mode embedding; and M3-SC denotes spectral clustering applied to the mode-3 matricization of the adjacency tensor.} Table~\ref{table:brain} reports the clustering errors of all compared unsupervised methods against the true labels. \jasarr{The error is the permutation-invariant mis-clustering error in \eqref{eq:network-hamming-loss}, with diagnosis labels used only after clustering for validation.} The proposed refinement method achieves the smallest error, $0.395$, outperforming TWIST ($0.444$), Tucker ($0.435$), and M3-SC ($0.492$). \jasarr{The majority-class baseline, which assigns every subject to the healthy-control group, has error $54/(54+70)=0.435$.} Although a higher accuracy was reported in \cite{relion2019network}, that analysis used supervised learning, whereas the present task is unsupervised clustering.  \jasarr{The improvement over the majority baseline is modest in absolute terms, and we interpret this example as a real-data benchmark for unsupervised graph-level clustering rather than as evidence of clinical diagnostic utility.}

\begin{table}[htbp]
\begin{center}
\begin{tabular}{|l|l|l|l|l|}
\hline
Method & Refinement & TWIST & Tucker & M3-SC \\ \hline
Error & \textbf{0.395} & 0.444 & 0.435 & 0.492 \\ \hline
\end{tabular}
\caption{Clustering error of the COBRE brain-connectivity data.}
\label{table:brain}
\end{center}
\end{table}

\subsubsection{Worldwide food trading networks}

As an unlabeled exploratory example, we next analyze worldwide food and agricultural commodity trading networks. The data were collected by \cite{de2015structural} and have been adopted in multi-layer network analysis \citep{jing2021community,fan2021alma}. The data contains an economic network in which layers represent different products, nodes are countries and edges at each layer represent trading relationships of a specific food product among countries.

We focus on the trading data in 2010 only. Following the preprocessing procedure in \cite{jing2021community}, we obtained a 30-layers network with 99 nodes at each layer. Each layer represents trading relationships between 99 countries/regions worldwide with respect to one of the 30 different food products. Together they form a $3$rd-order tensor of dimension $99\times99\times30$.

\begin{figure}[!htbp]
\centering
\begin{subfigure}[b]{.48\linewidth}
	\includegraphics[width=\linewidth]{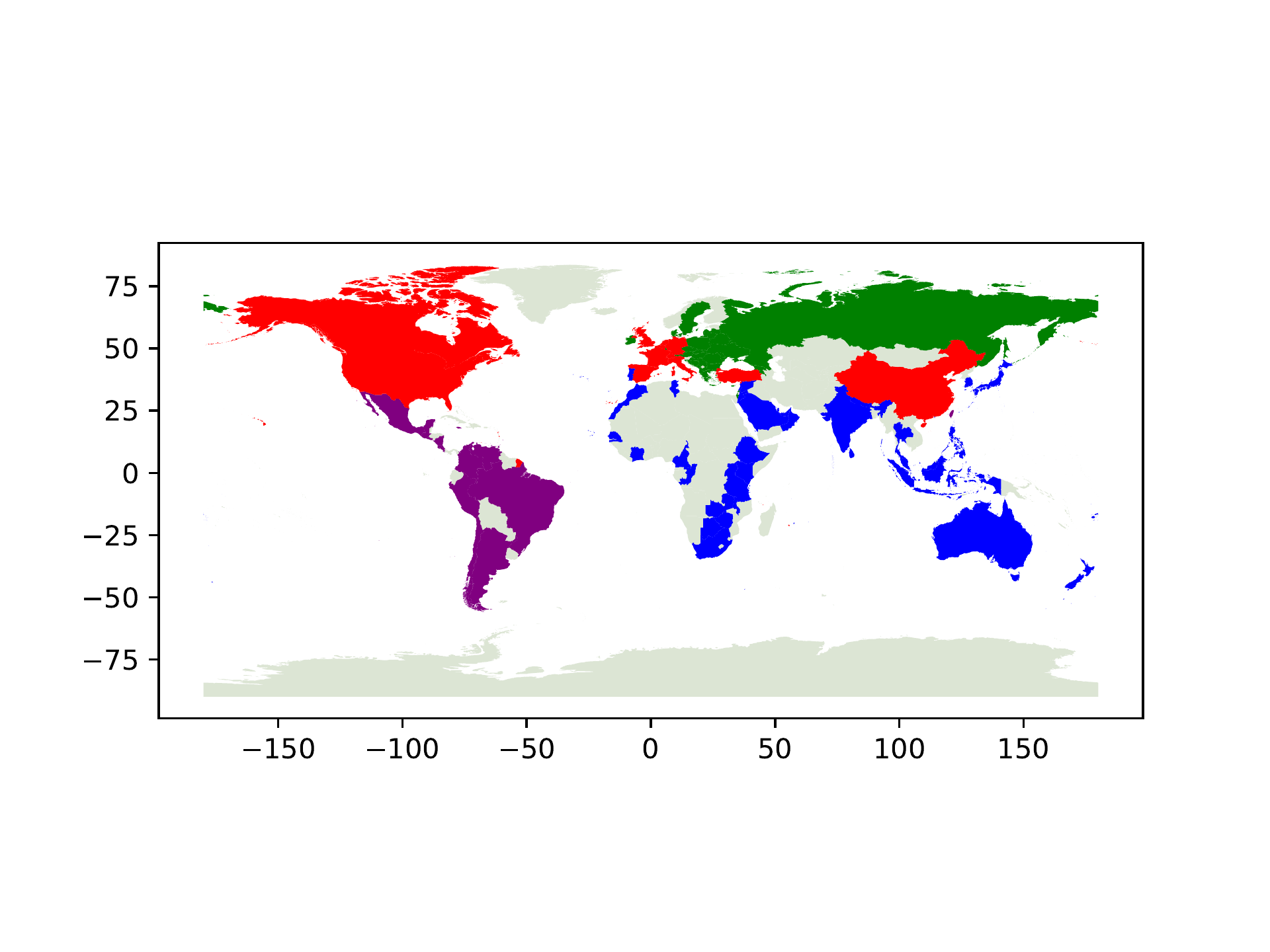}
	\caption{\small{Membership of countries for networks in cluster 1.}}
\end{subfigure}
\begin{subfigure}[b]{.48\linewidth}
	\includegraphics[width=\linewidth]{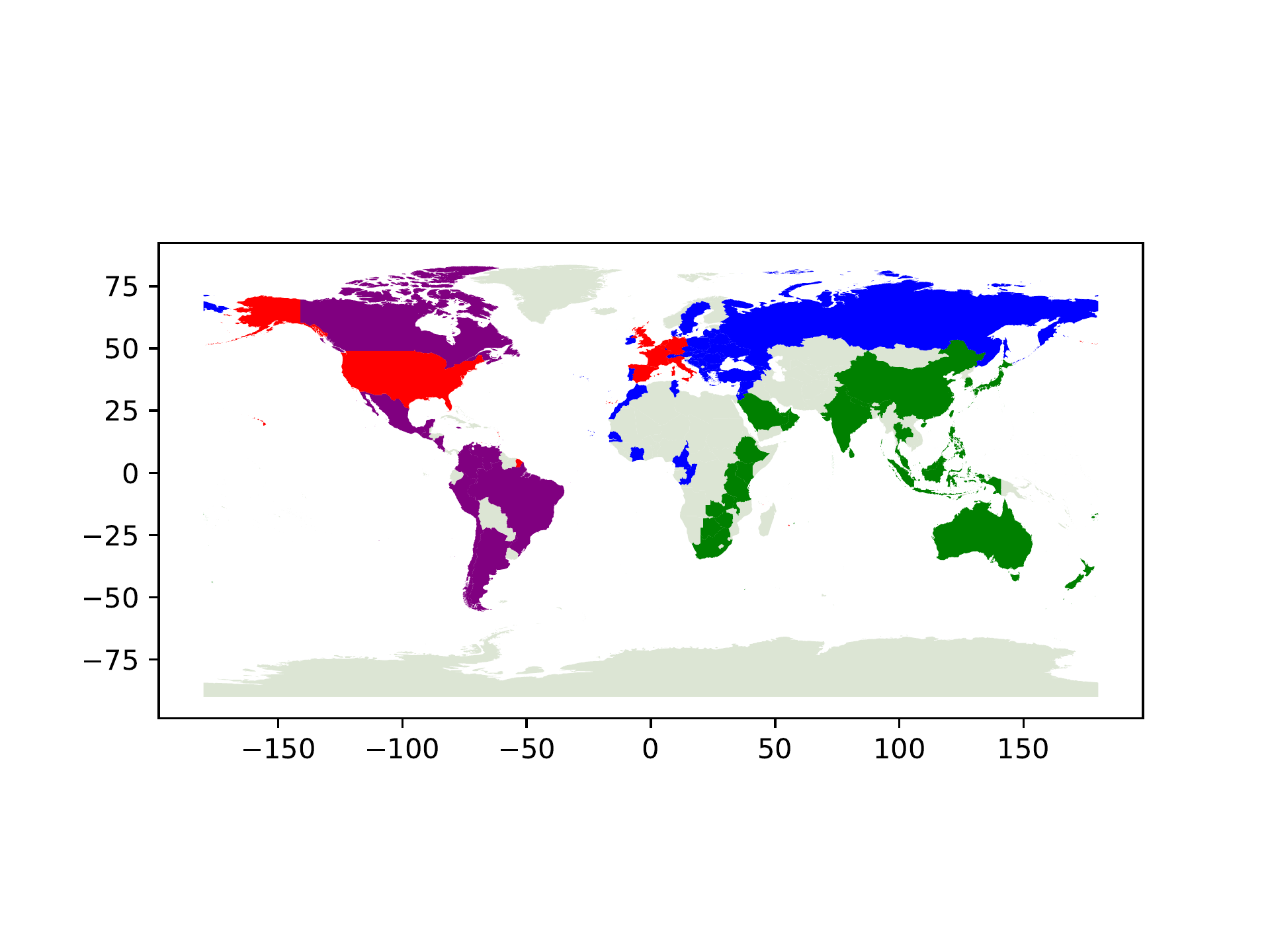}
	\caption{\small{Membership of countries for networks  in cluster 2.}}
\end{subfigure}
\caption{Memberships of countries on {two different types of food trading networks.}}
\label{fig:country_embedding}
\end{figure}

\begin{figure}[!htbp]
\centering
\begin{subfigure}[b]{.48\linewidth}
	\includegraphics[width=\linewidth]{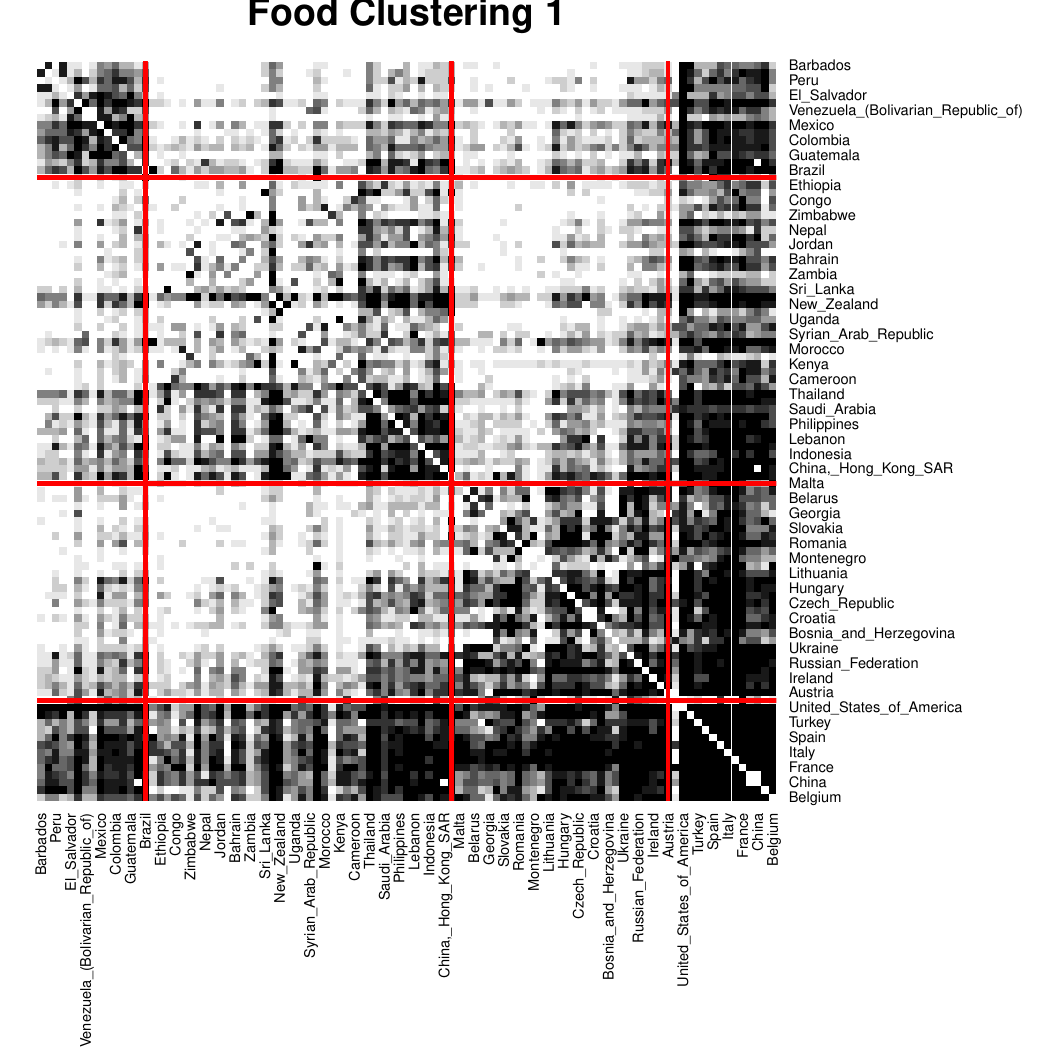}
	\caption{\small{Heat map of networks in cluster 1.}}
\end{subfigure}
\begin{subfigure}[b]{.48\linewidth}
	\includegraphics[width=\linewidth]{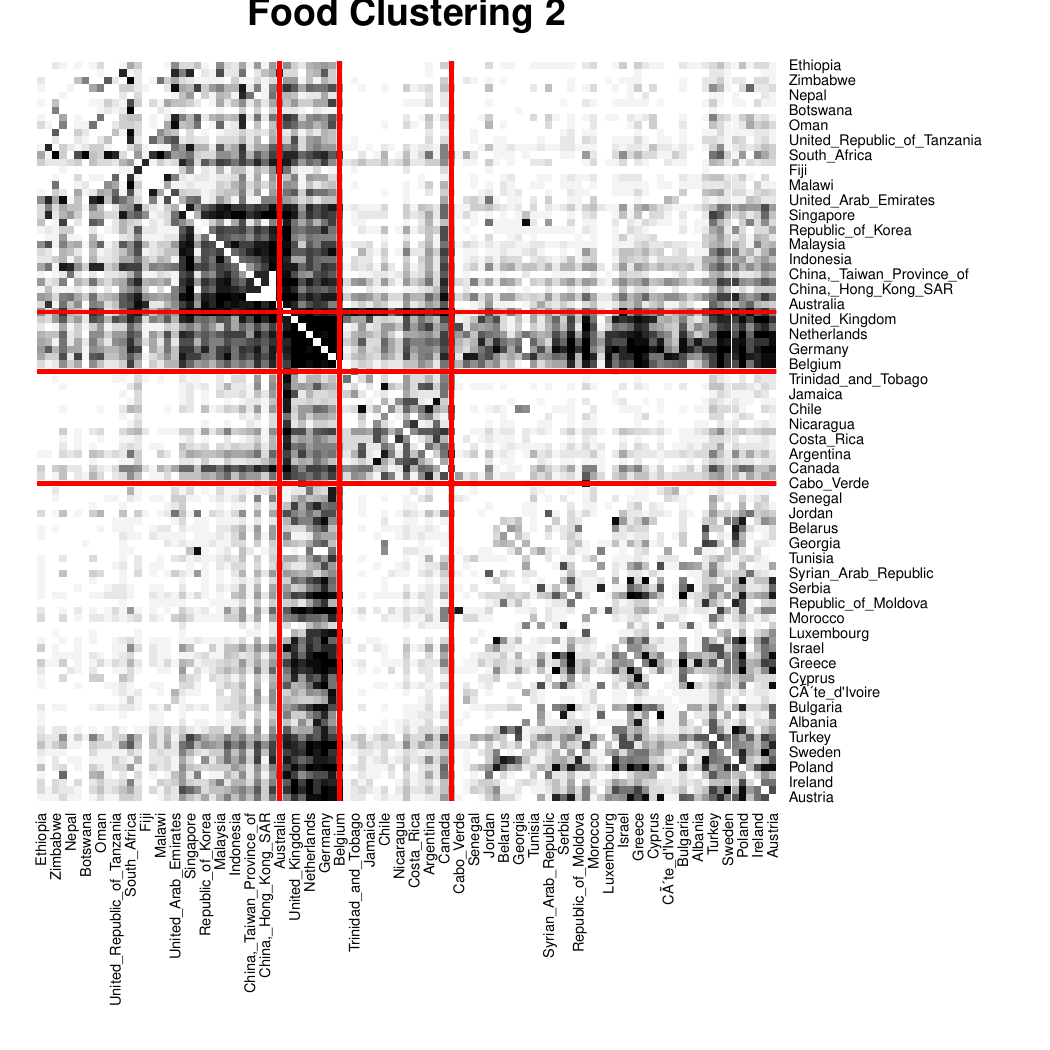}
	\caption{\small{Heat map of networks in cluster 2.}}
\end{subfigure}
\caption{Heat maps of two types of networks.}
\label{fig:food_networks}
\end{figure}

We apply the refinement algorithm to the data tensor. The resulting two clusters of layers are listed in Table \ref{food_cluster}. We then apply the spectral method to the sum of each cluster of networks separately (here we have two clusters) to find the community structures for each cluster, in order to obtain the clustering result of countries. The membership of 99 countries with clustering results from $K$-means with $K=4$ are shown in Figure \ref{fig:country_embedding}.
For the two types of networks, we plot in Figure \ref{fig:food_networks} the sum of adjacency matrices with nodes arranged according to the community labels to visualize different community structures of two network types.

\begin{table}
	\begin{center}
		\begin{tabular}{ | l  p{10cm} |}
			\hline

			Food cluster 1: &
			\tiny{Beverages, non alcoholic;
				Food prep nes;
				Chocolate products nes;
				Crude materials;
				Fruit, prepared nes;
				Beverages, distilled alcoholic;
				Juice, fruit nes;
				Pastry;
				Sugar confectionery;
				Wine.
			}    \\ \hline
			Food cluster 2: &
			\tiny{Cheese, whole cow milk;
				Cigarettes;
				Flour, wheat;
				Beer of barley;
				Cereals, breakfast;
				Coffee, green;
				Milk, skimmed dried;
				Maize;
				Macaroni;
				Oil, palm;
				Milk, whole dried;
				Oil, essential nes;
				Rice, milled;
				Sugar refined;
				Tea;
				Spices, nes;
				Vegetables, preserved nes;
				Waters,ice etc;
				Vegetables, fresh nes;
				Tobacco, unmanufactured;
			}  \\
			\hline
		\end{tabular}
		\caption{Estimated clusters of commodity-specific trading networks in the 2010 worldwide food and agricultural commodity trade data.}
	\label{food_cluster}
\end{center}
\end{table}

We make the following remarks from Table \ref{food_cluster}, Figures \ref{fig:country_embedding} and  \ref{fig:food_networks}.
\begin{enumerate}
\item
Food types are mainly divided into two categories. The first category comprises food products that have a longer shelf life, facilitating a more globalized trading pattern. Conversely, the second category includes foods that typically have a shorter shelf life or whose transportation costs are high relative to their overall value. Thus, the trade of these foods is predominantly conducted within nearby regions.

\item From the first category, a handful of countries, such as China, Canada, the United Kingdom, the United States, France, and Germany, are notably active in trading both internationally and amongst themselves. However, outside of these hub nations, other countries are grouped mainly by geographical location, such as America, East Europe, Africa, and Australia. This pattern demonstrates the dominance of these hub countries in the global food trade, while other nations engage more in regional trading to minimize delivery costs.

\item In the second category, countries are primarily clustered based on their geographical location. Apart from the USA, which is grouped with European countries, the four main clusters consist of America, Western Europe, Eastern Europe and West Africa, and nations surrounding the Indian Ocean - Asia, East Africa, and Australia. Regional trading is crucial in this category as it helps to keep food costs low due to reduced transportation costs and maintains food freshness due to faster delivery times.
\end{enumerate}

\section{Discussions}\label{sec:discuss}

Several directions merit further study. The likelihood-refinement argument in Section~\ref{sec:MMSBM} and \ref{sec:MMPBM} suggests that the \ZL{Chernoff} exponent may persist more broadly for discrete or weighted network mixtures, provided that one can estimate the layer-specific edge distributions accurately enough. Extending the present analysis beyond Bernoulli and Poisson edges requires distribution-specific concentration and plug-in likelihood control. It is also of interest to develop initialization methods that remain effective when assortative and disassortative layers coexist, and to clarify the gap, if any, between the statistical and computational limits in more general multi-layer models. \jasarr{Extending the theory to directed networks is an important direction, but it requires asymmetric edge-probability matrices, possibly separate sender and receiver memberships, and different concentration arguments.}

\noindent\textit{Supplementary material.} The appendices contain the deferred MMPBM theory, the additional numerical studies, and all proofs.

\noindent\textit{Data availability statement.} The data supporting the findings of this study were derived from public domain resources. The brain-connectivity data are from the COBRE dataset (\url{https://fcon_1000.projects.nitrc.org/indi/retro/cobre.html}), described in \cite{aine2017multimodal} and processed as in \cite{relion2019network}. The worldwide food-trading network data are derived from public FAO trade data (\url{https://www.fao.org/faostat/}), following \cite{de2015structural,jing2021community}.
\acks{Dong Xia's research was partially supported by Hong Kong RGC Grant GRF 16302020 and GRF 16301622.}
\clearpage
\appendix
\setcounter{theorem}{0}
\renewcommand{\thetheorem}{A.\arabic{theorem}}
\renewcommand{\theHtheorem}{appendix.\arabic{theorem}}
\setcounter{equation}{0}
\renewcommand{\theequation}{A.\arabic{equation}}
\renewcommand{\theHequation}{appendix.\arabic{equation}}
\setcounter{figure}{0}
\renewcommand{\thefigure}{A.\arabic{figure}}
\renewcommand{\theHfigure}{appendix.\arabic{figure}}
\setcounter{table}{0}
\renewcommand{\thetable}{A.\arabic{table}}
\renewcommand{\theHtable}{appendix.\arabic{table}}
\setcounter{algocf}{0}
\renewcommand{\thealgocf}{A.\arabic{algocf}}
\section{Further Results for the MMPBM}

\subsection{\texorpdfstring{\textcolor{blue}{Oracle property for likelihood-based refinement}}{Oracle property for likelihood-based refinement}}
\begin{movedenv}

Suppose that the edge intensity matrices $\bP_1$ and $\bP_2$ are known in advance,  if an observed network $\bX\in \ZZ_+^{d\times d}$ is known sampled from either $\bP_1$ or $\bP_2$, the log-likelihood can be written as
$$
\ell(\bX|z):=\sum_{\omega\in\calJ_d} \bigg\{-\bP_z(\omega)+\bX(\omega)\log \bP_z(\omega)-\log\big(\bX(\omega)!\big)\bigg\},
$$
where $z$ takes values in $[2]$.  The likelihood-based Lloyd's algorithm assigns the label to $\bX$ which maximizes the log-likelihood $\ell(\bX|z)$.  Without loss of generality,  assume the true label of $\bX$ is $1$.  This network is mis-clustered by Lloyd's algorithm \ZL{only if}
$$
\sum_{\omega\in\calJ_d}\big[\bP_1(\omega)-\bP_2(\omega)\big]\ZL{\ge}\sum_{\omega\in\calJ_d} \bX(\omega)\log \frac{\bP_1(\omega)}{\bP_2(\omega)}
$$
The following lemma characterizes the oracle mis-clustering error rate of Lloyd's algorithm.

\begin{lemma}\label{lem:oracle-MMPBM}
Suppose that $\bP_1$ and $\bP_2$ are known. For $m\in[2]$, let $\PP_m$ be the law under which
$\bX(\omega)\stackrel{\rm ind.}{\sim}{\rm Poisson}(\bP_m(\omega))$, $\omega\in\calJ_d$.
For $\widehat z\in\argmax_{m\in[2]}\ell(\bX\mid m)$,
\[
 \max_{m\in[2]}\PP_m(\widehat z\ne m)\le\ZL{\exp\brac{-I^*}},
\]
\jasaadd{where $I^*$ is defined in \eqref{eq:def-I-star-mmpbm}.}
\end{lemma}

\ZL{The oracle bound implies consistency when $I^{\ast}\to\infty$. For fixed $\alpha>1$, necessity for uniform consistency follows from the anchored-label testing reduction below: its equal-prior testing error is at least $e^{-2I^*}/4$.} We remark that Poisson is sub-exponential, while the distributions under MMSBM (Section~\ref{sec:MMSBM}) and LrMM \citep{lyu2022optimal} are sub-Gaussian.

\end{movedenv}

\begingroup
\subsection{Minimax lower bound}
For fixed component memberships and intensity matrices, retain the
parameter set $\Theta_{\bz}^{(\bB_{1,2},\bsig_{1,2})}$ of
\eqref{eq:parameter-space-MMPBM}, with \ZL{a fixed layer-balance constant
$\alpha>1$}.
\begin{theorem}\label{thm:lower-MMPBM}
For each fixed pair satisfying \eqref{eq:parameter-space-MMPBM},
suppose \ZL{$n\to\infty$} and $\ZL{I^*}\to\infty$. \ZL{If $\alpha>1$ is fixed}, then
\[
 \inf_{\widehat\bz}\sup_{\bz^*\in
 \Theta_{\bz}^{(\bB_{1,2},\bsig_{1,2})}}
 \EE h(\widehat\bz,\bz^*)
 \ge\ZL{\exp(-[1+o(1)]I^*)}.
\]
The infimum is over all estimators $\widehat\bz:\mathbb N_0^{n\times d(d+1)/2}\to[2]^n$, with the fixed component pair available as input.
\end{theorem}
Theorem~\ref{thm:MMPBM-main} gives the upper bound when the component
parameters are unknown.
\endgroup
\section{Simulation Studies}

\subsection{Simulation studies}
\ZL{These experiments use the practical full-sample refinement implementation.}
We conduct several simulations to test the performance of the refinement algorithm on the MMSBM with different choices of network sparsity, ``out-in" ratio, number of layers and the size of each layer. We use K-means as the clustering algorithm.
We compare the mis-clustering rate across 150 replications of each experiment to identify the average performance.

We generate the data according to MMSBM in the following fashion. \ZL{For an experiment with $m$ network types, the underlying class $s_l$ of layer $l$ is sampled uniformly from $[m]$. The experiments below include $m>2$ as empirical extensions; the theoretical results in this paper concern two components.}  The membership $z_i^j$ for node $i$ in layer type $j$ is generated from the multinomial distribution with $\PP(z_i^j=s)=1/K,\ s=1,\cdots,K.$ We choose the probability matrix as
$\bB=p\bI_K+q(\bone_K\bone_K^{\top}-\bI_K),$ where $\bone_K$ is a $K$-dimensional all-one vector and $\bI_K$ is the $K\times K$ identity matrix . Let $\alpha=q/p$ be the out-in ratio.

We compare the performance of the refinement algorithm with the TWIST \citep{jing2021community}, Tucker decomposition initialized by HOSVD (HOSVD-Tucker) and spectral clustering applied to the mode-3 matricization of $\bcalX$ (M3-SC). The function ``\textit{tucker}" from the R package "\textit{rTensor}" \cite{li2018rtensor} is used to apply Tucker decomposition for HOSVD-Tucker.

In Simulation 1, the networks are generated with the number of nodes $d=40$, the number of layers $n=40$, number of types of networks $m=2,$ number of communities of each network $K=2$ and out-in ratio of each layer $\alpha=0.75.$ The $p$ of each layer varies from $0.1$ to $0.8$, which reflects the average degree with fixed out-in ratio.

In Simulation 2, the networks are generated as in Simulation 1, except that \ZL{$K=4$, $m=4$, and} the \ZL{within-community edge probability is} $p=0.4,$ the number of layers $n=40$ and the out-in ratio $\alpha$ of each layer varies from $0.1$ to $0.9.$

In Simulation 3, the networks are the same as in Simulation 2, except that \ZL{$m=2$}, the out-in ratio $\alpha=0.75$ and the number of layers $n$ varies from $20$ to $80$.

In Simulation 4, the networks are the same as in Simulation 3, except that \ZL{$K=3$, $m=3$, the number of layers is $n=40$, and the number of nodes in each layer, $d$, varies from $10$ to $100$}.

\begin{figure}
	\centering
	\begin{subfigure}[b]{.48\linewidth}
		\includegraphics[width=\linewidth]{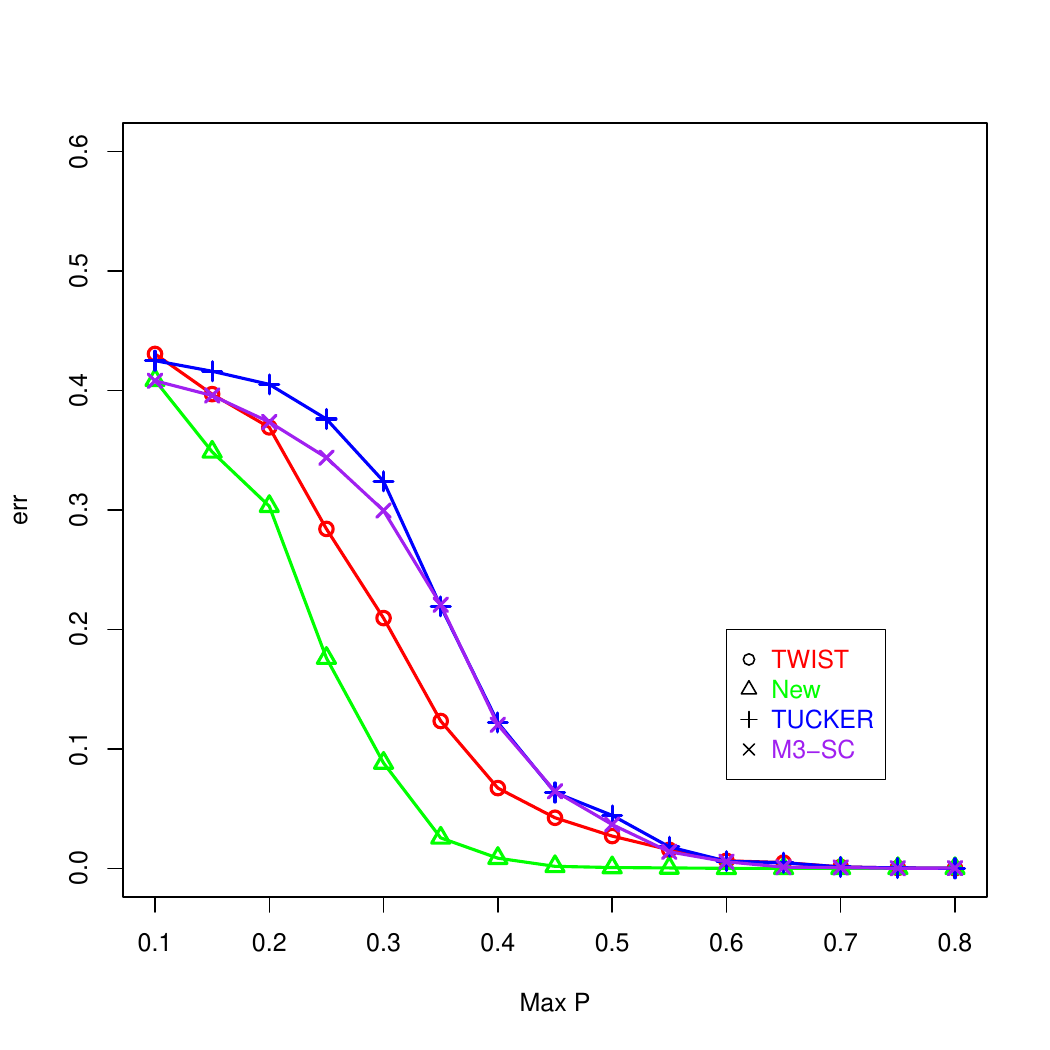}
		\caption{{The result of Simulation 1: $d=40$, $K=2$, $m=2$, $n=40$, $\alpha=0.75$, varying $p.$ }}
	\end{subfigure}
	\begin{subfigure}[b]{.48\linewidth}
		\includegraphics[width=\linewidth]{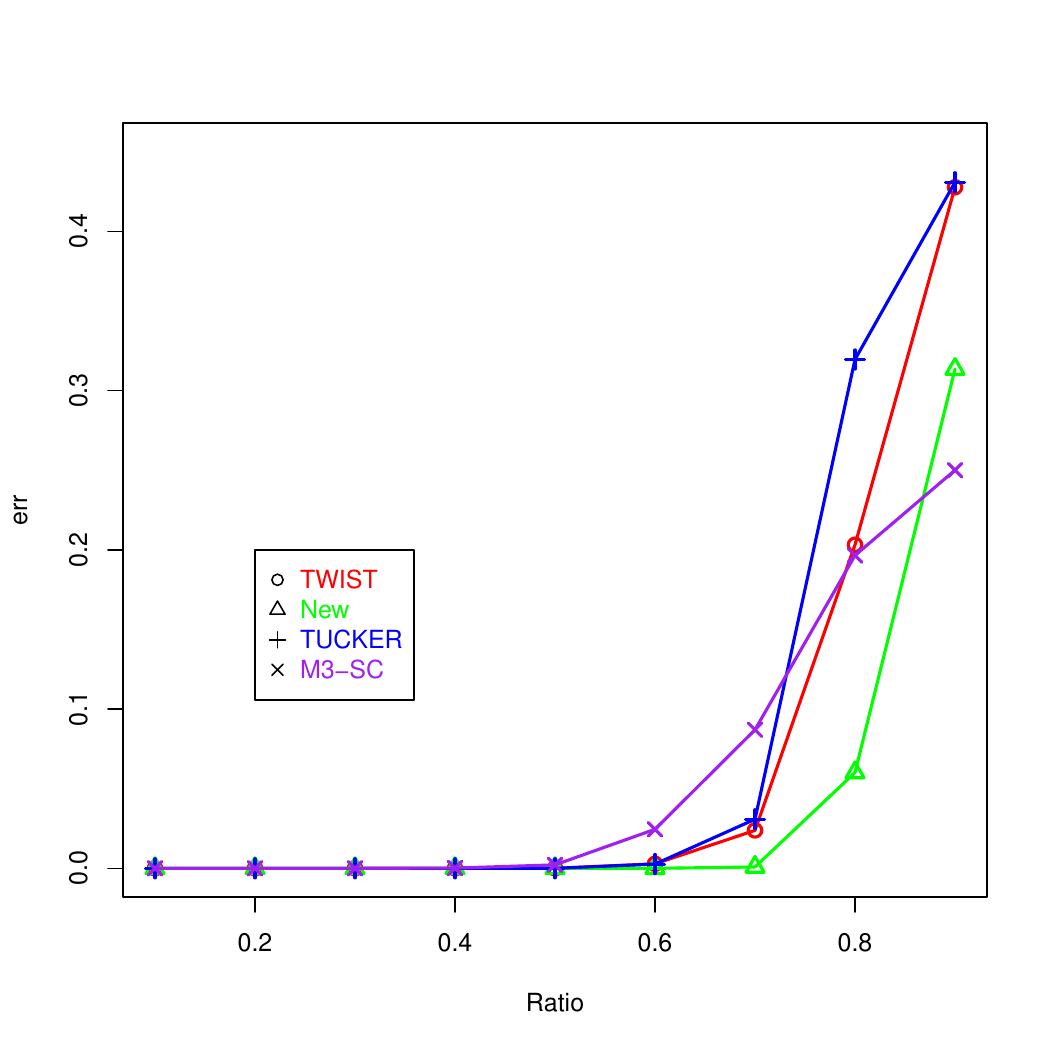}
		\caption{{The result of Simulation 2: $d=40$, $K=4$, $m=4$, $p=0.4$, $n=40$, varying $\alpha.$ }}
	\end{subfigure}
	\begin{subfigure}[b]{.48\linewidth}
		\includegraphics[width=\linewidth]{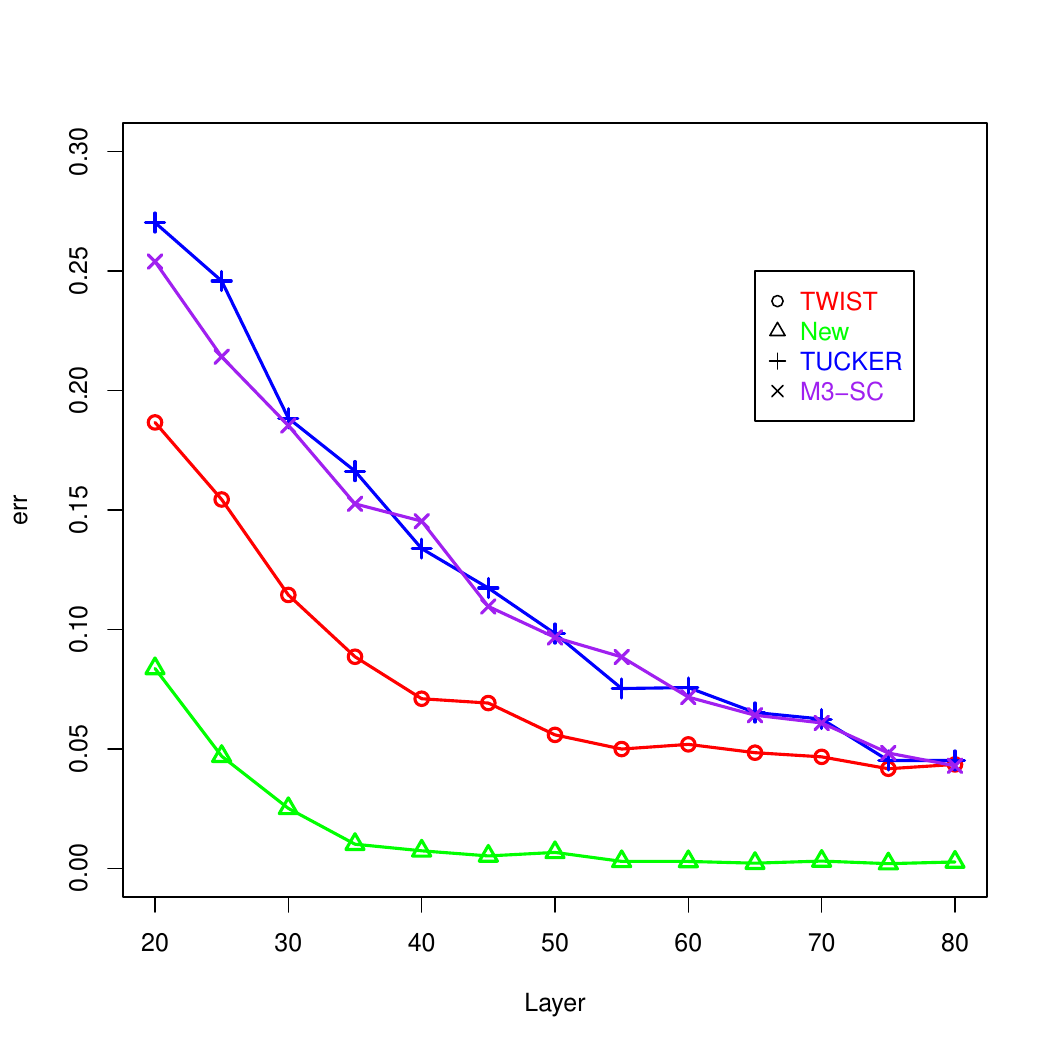}
		\caption{{The result of Simulation 3: $d=40$, $K=4$, $m=2$, $p=0.4$, $\alpha=0.75$, varying $n.$ }}
	\end{subfigure}
	\begin{subfigure}[b]{.48\linewidth}
		\includegraphics[width=\linewidth]{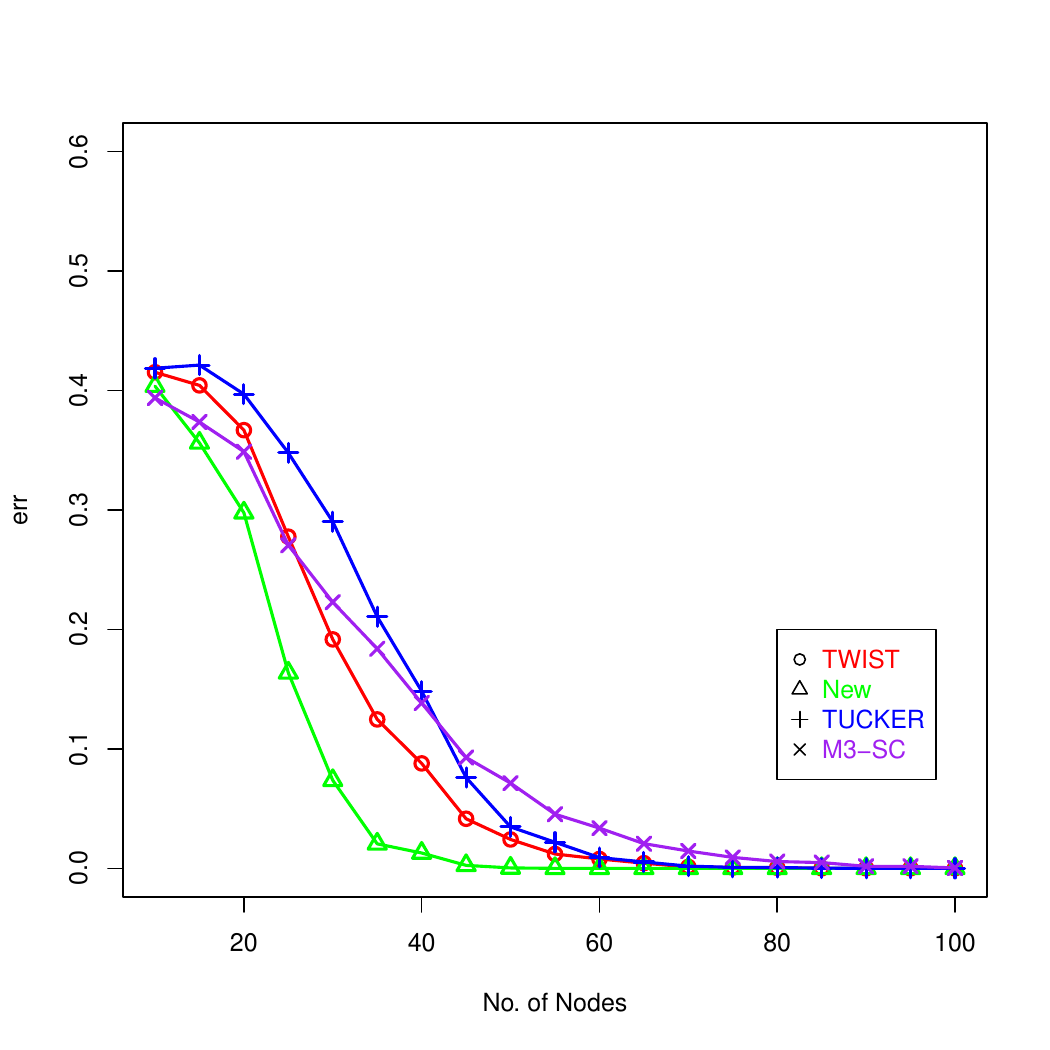}
		\caption{{\ZL{The result of Simulation 4: $K=3$, $m=3$, $p=0.4$, $\alpha=0.75$, $n=40$, varying $d$.} }}
	\end{subfigure}
	\caption{The refinement algorithm is the best overall, particularly when the signal is not strong enough, for instance, $p<0.5$ in (a), $\alpha>0.6$ in (b), \jasaadd{all scenarios in (c), and $d<60$ in (d). From Simulation 3, the other three methods reach at most a 95\% accuracy rate when $n$ increases to 80.} However, the refinement algorithm nearly achieves 100\% accuracy when $n$ is only 40.
	}
	\label{fig:network membership}
\end{figure}

The results of the simulations are shown in Figure \ref{fig:network membership}, from which we can draw several key findings:
\begin{enumerate}
	\item As we anticipated from our theoretical results, mis-clustering rates for all three methods decrease when the average degree of each layer increases, the out-in ratio of each layer decreases, the number of layers increases, and the size of each layer increases.
	\item Our refinement algorithm notably surpasses the other methods in terms of accuracy as the number of layers increases. While the other methods only achieve a 95\% accuracy rate when the number of layers reaches 80, the refinement algorithm nearly hits 100\% accuracy with only 40 layers (as observed in Simulation 3).
	\item The refinement algorithm consistently outperforms the other three methods across all simulations. Its superior performance is especially noticeable when the signal strength is relatively weak, for example when $p<0.5$ in Simulation 1, $\alpha>0.6$ in Simulation 2, in all scenarios in Simulation 3, and when $d<60$ in Simulation 4.
	\item While the refinement algorithm uses TWIST as an initialization step, the improvements obtained from the refinement step are significant, indicating its pivotal role in the performance of the algorithm.
\end{enumerate}

\section{Proof of Main Results}\label{sec:main-proofs}
In the likelihood proofs, identify a symmetric network matrix with
its coordinates in $\calJ_d=\{(j,k):1\le j\le k\le d\}$.
The observation space is $\mathcal X_d=\{0,1\}^{\calJ_d}$ for Bernoulli
observations and $\mathcal X_d=\mathbb N_0^{\calJ_d}$ for Poisson
observations, where $\mathbb N_0=\{0,1,2,\ldots\}$.
For $\bu\in\mathcal X_d$ and $m\in\{1,2\}$, the layer mass function is
\[
 f_m(\bu)=
 \begin{cases}
 \displaystyle\prod_{\omega\in\calJ_d}
 \bP_m(\omega)^{u_\omega}[1-\bP_m(\omega)]^{1-u_\omega},
 &\text{in the Bernoulli model},\\[0.5em]
 \displaystyle\prod_{\omega\in\calJ_d}
 \frac{e^{-\bP_m(\omega)}\bP_m(\omega)^{u_\omega}}{u_\omega!},
 &\text{in the Poisson model}.
 \end{cases}
\]
Write $\bx=(\bX(\omega))_{\omega\in\calJ_d}$ and let
$\mathbb P_m,\mathbb E_m$ denote probability and expectation under $f_m$.
Unless specified otherwise, $\bu\in\mathcal X_d$ and
$\omega\in\calJ_d$ in sums and products.

\ZL{For Bernoulli probabilities, write}
\[
 \ZL{D_{\rm KL}(p\Vert q)
 =p\log\frac pq+(1-p)\log\frac{1-p}{1-q},\qquad p,q\in(0,1).}
\]

\subsection{Proof of Lemma \ref{lem:optimal-rate-chernoff}}
\ZL{
For $m\in[2]$, misclassification implies $f_{3-m}(\bx)\ge f_m(\bx)$,
including when the likelihoods are tied. For every $s\in[0,1]$,
\begin{align*}
 \mathbb P_m(\widehat z\ne m)
 &\le\mathbb E_m\left[\left\{\frac{f_{3-m}(\bx)}{f_m(\bx)}\right\}^{s}\right]\\
 &=\prod_{\omega}\Big\{
 \bP_m(\omega)^{1-s}\bP_{3-m}(\omega)^s
 +[1-\bP_m(\omega)]^{1-s}[1-\bP_{3-m}(\omega)]^s\Big\}.
\end{align*}
Minimizing over $s$ gives $e^{-I^*}$ by
\eqref{eq:def-I-star-mmsbm}; for $m=2$, use the substitution $s\mapsto1-s$.
}

\subsection{Proof of Lemma \ref{lem:oracle-MMPBM}}
\ZL{
For $m\in[2]$, the same likelihood-ratio argument gives, for every $s\in[0,1]$,
\begin{align*}
 \mathbb P_m(\widehat z\ne m)
 &\le\mathbb E_m\left[\left\{\frac{f_{3-m}(\bx)}{f_m(\bx)}\right\}^{s}\right]\\
 &=\exp\left(-\sum_{\omega}\left[
 (1-s)\bP_m(\omega)+s\bP_{3-m}(\omega)
 -\bP_m(\omega)^{1-s}\bP_{3-m}(\omega)^s\right]\right).
\end{align*}
Minimizing over $s$ and using \eqref{eq:def-I-star-mmpbm} proves the claim.
}

\subsection{Proof of Lemma~\ref{lem:optimal-rate-dev}}
\begingroup
Suppose the true label is $1$. Recall that $\hat z$ is defined in
\eqref{eq:MLE-emp}. On the event $\calE$ in the lemma, choose $0<\tilde\rho=o(1)$ such that
\[
 \max_{m\in\{1,2\}}\max_{\omega\in\calJ_d}
 |\hat\bP_m(\omega)-\bP_m(\omega)|
 \le\tilde\rho\frac{I^*}{d(d+1)/2}.
\]
By definition, on the event $\calE$,
\begin{align*}
\II(\hat z=2)\II(\calE)\leq&\II\left\{\sum_{\omega\in\calJ_d}\bar\bX(\omega)\log\frac{\bP_{2}(\omega)[1-\bP_{1}(\omega)]}{\bP_{1}(\omega)[1-\bP_{2}(\omega)]}\ZL{\ge}\sum_{\omega\in\calJ_d}D_{KL}\left(\bP_{1}(\omega)||\bP_{2}(\omega)\right)-\tilde\rho^{1/3}I^*\right\}\\
&+\II(\calE)\II\Bigg\{\sum_{\omega\in\calJ_d}\bX(\omega)\left\{\log\frac{\hat\bP_{2}(\omega)[1-\hat\bP_{1}(\omega)]}{\hat\bP_{1}(\omega)[1-\hat\bP_{2}(\omega)]}-\log\frac{\bP_{2}(\omega)[1-\bP_{1}(\omega)]}{\bP_{1}(\omega)[1-\bP_{2}(\omega)]}\right\}\\
&\qquad\qquad\qquad \ZL{\ge}\sum_{\omega\in\calJ_d}\Bigg[\log\frac{1-\hat\bP_{1}(\omega)}{1-\hat\bP_{2}(\omega)}-\log\frac{1-\bP_{1}(\omega)}{1-\bP_{2}(\omega)}\Bigg]+\tilde\rho^{1/3}I^*\Bigg\},
\end{align*}
where $\bar\bX(\omega):=\bX(\omega)-\bP_{1}(\omega)$.  The first term is the oracle term. \ZL{Let $s_*\in[0,1]$ be the optimizer in \eqref{eq:def-I-star-mmsbm}. Applying Markov's inequality at $s_*$, as in the proof of Lemma~\ref{lem:optimal-rate-chernoff}, gives}
\begin{align}\label{eq:optimal-rate-first-term}
&\PP\Bigg\{\sum_{\omega\in\calJ_d}\bar\bX(\omega)
\log\frac{\bP_{2}(\omega)[1-\bP_{1}(\omega)]}
{\bP_{1}(\omega)[1-\bP_{2}(\omega)]}\nonumber\\
&\hspace{3em}\ZL{\ge}\sum_{\omega\in\calJ_d}
D_{KL}\left(\bP_{1}(\omega)||\bP_{2}(\omega)\right)-\tilde\rho^{1/3}I^*\Bigg\}\nonumber\\
&\quad\le \ZL{\exp\big(-I^*+s_*\tilde\rho^{1/3}I^*\big)\le\exp\big(-(1-\tilde\rho^{1/3})I^*\big)}.
\end{align}

It remains to control the perturbation term.
\ZL{By Jensen's inequality, for every $s\in[0,1]$, we have
\begin{align*}
 -\log\sum_{\bu} f_1(\bu)^{1-s}f_2(\bu)^s
 &=-\log\mathbb E_1\exp\left(s\log\frac{f_2(\bx)}{f_1(\bx)}\right)\\
 &\le-s\mathbb E_1\log\frac{f_2(\bx)}{f_1(\bx)}
 =sD_{\rm KL}(f_1\Vert f_2)
 \le D_{\rm KL}(f_1\Vert f_2),
\end{align*}
where $D_{\rm KL}(f_1\Vert f_2)=\sum_{\bu}f_1(\bu)\log(f_1(\bu)/f_2(\bu))$.
Taking the maximum over $s$, we obtain
$I^*\le D_{\rm KL}(f_1\Vert f_2)$.
For each Bernoulli edge, by
$\log x\le x-1$, we have
\begin{align*}
 D_{\rm KL}(\bP_1(\omega)\Vert\bP_2(\omega))
 &\le\frac{\bP_1(\omega)^2}{\bP_2(\omega)}
 +\frac{[1-\bP_1(\omega)]^2}{1-\bP_2(\omega)}-1\\
 &=\frac{[\bP_1(\omega)-\bP_2(\omega)]^2}
 {\bP_2(\omega)[1-\bP_2(\omega)]}.
\end{align*}
Since the edge probabilities are mutually comparable and bounded away from one,
summing over the independent edges yields}
\[
\frac{I^*}{d(d+1)/2}\lesssim \frac{1}{d(d+1)/2}\sum_{\omega\in\calJ_d}\frac{[\bP_1(\omega)-\bP_2(\omega)]^2}{\bP_1(\omega)\vee \bP_2(\omega)}\lesssim \min_{\omega\in\calJ_d}\min_{m\in\{1,2\}}\bP_m(\omega),
\]
where the last inequality uses the homogeneous parameter space in \eqref{eq:parameter-space}.  Hence, on $\calE$ and for all sufficiently large $n$,
\[
|\hat\bP_m(\omega)-\bP_m(\omega)|\le \frac{1}{2}\bP_m(\omega),\qquad \omega\in\calJ_d,\ m=1,2.
\]
\ZL{The same information bound and $\tilde\rho=o(1)$ give
\[
 \max_{m,\omega}|\hat\bP_m(\omega)-\bP_m(\omega)|=o(1)
 \quad\text{on }\calE.
\]
Since the true probabilities are bounded away from one, for all
sufficiently large $n$, on $\calE$, we obtain
\begin{equation}\label{eq:testing-bernoulli-denominators}
 \min\{\hat\bP_m(\omega),\bP_m(\omega)\}
 \ge\tfrac12\bP_m(\omega),\qquad
 1-\max\{\hat\bP_m(\omega),\bP_m(\omega)\}\ge c>0,
\end{equation}
for a fixed constant $c>0$.}
The mean-value theorem and \eqref{eq:testing-bernoulli-denominators}
therefore imply
\begin{align*}
 \left|\log\frac{\hat\bP_m(\omega)}{\bP_m(\omega)}\right|
 &\le\frac{|\hat\bP_m(\omega)-\bP_m(\omega)|}
 {\min\{\hat\bP_m(\omega),\bP_m(\omega)\}},\\
 \left|\log\frac{1-\hat\bP_m(\omega)}{1-\bP_m(\omega)}\right|
 &\le\frac{|\hat\bP_m(\omega)-\bP_m(\omega)|}
 {1-\max\{\hat\bP_m(\omega),\bP_m(\omega)\}}.
\end{align*}
On $\calE$, using $|\log(1-x)-\log(1-y)|\lesssim |x-y|$ \ZL{for $x,y$ uniformly bounded away from one}, we have
\[
\Bigg|\sum_{\omega\in\calJ_d}\Bigg[\log\frac{1-\hat\bP_{1}(\omega)}{1-\hat\bP_{2}(\omega)}-\log\frac{1-\bP_{1}(\omega)}{1-\bP_{2}(\omega)}\Bigg]\Bigg|\lesssim \sum_{m\in\{1,2\}}\sum_{\omega\in\calJ_d}|\hat\bP_m(\omega)-\bP_m(\omega)|\lesssim \tilde\rho I^*\le \frac{\tilde\rho^{1/3} I^*}{8}
\]
for all sufficiently large $n$.
Therefore $\hat\bP_m(\omega)\asymp \bP_m(\omega)$ uniformly in $\omega$, and for every such $\omega$, the two mean-value bounds give
\begin{align}\label{eq:testing-bernoulli-log-envelope}
 &\left|\log\frac{\hat\bP_m(\omega)}{\bP_m(\omega)}\right|
 +\left|\log\frac{1-\hat\bP_m(\omega)}{1-\bP_m(\omega)}\right|\nonumber\\
 &\qquad\le C\tilde\rho\cdot
 \frac{I^*}{d(d+1)\bP_1(\omega)},\qquad m=1,2,
\end{align}
after using $\bP_1(\omega)\asymp \bP_2(\omega)$ uniformly over
$\omega$ under \eqref{eq:parameter-space}. Consequently, on $\calE$, the second indicator above is bounded by
\[
\II\left[\sum_{\omega\in\calJ_d}\bX(\omega)\frac{C\tilde\rho I^*}{d(d+1)\bP_1(\omega)}\ZL{\ge}\frac{\tilde\rho^{1/3} I^*}{4}\right].
\]
It suffices to bound the probability of the last event, whose coefficients
are deterministic. By the moment generating function of the Bernoulli distribution,
we have
\begin{align*}
 &\log\EE\exp\left(
 \bX(\omega)\frac{C\tilde\rho^{1/3}I^*}
 {d(d+1)\bP_1(\omega)}\right)\\
 &\quad=\log\left(1-\bP_1(\omega)+\bP_1(\omega)
 \exp\left(\frac{C\tilde\rho^{1/3}I^*}
 {d(d+1)\bP_1(\omega)}\right)\right)\\
 &\quad\le\frac{2C\tilde\rho^{1/3}I^*}{d(d+1)}.
\end{align*}
where we used $\log(1+x)\le x$, $e^x-1\le2x$ for small $x>0$ and
$\tilde\rho^{1/3}I^*/[d(d+1)\bP_1(\omega)]=o(1)$.
By Chernoff's inequality and independence, we obtain
\begin{align*}
&\PP\left[\sum_{\omega\in\calJ_d}\bX(\omega)
 \frac{C\tilde\rho I^*}{d(d+1)\bP_1(\omega)}
 \ZL{\ge}\frac{\tilde\rho^{1/3}I^*}4\right]\\
&\quad\le \exp\left(-\frac{\tilde\rho^{-1/3}I^*}4\right)
 \prod_{\omega\in\calJ_d}\EE\exp\left(
 \bX(\omega)\frac{C\tilde\rho^{1/3}I^*}
 {d(d+1)\bP_1(\omega)}\right)\\
&\quad\le \exp\left(-\frac{\tilde\rho^{-1/3}I^*}4
                         +C\tilde\rho^{1/3}I^*\right)\\
&\quad\le \exp(-c\tilde\rho^{-1/3}I^*)=o(e^{-I^*})
\end{align*}
for all sufficiently large $n$, since $\tilde\rho=o(1)$ and $I^*\to\infty$.  Combining this observation with \eqref{eq:optimal-rate-first-term}, we obtain
\[
\PP(\hat z=2,\calE)\leq \ZL{\exp\big(-[1-o(1)]I^*\big)}.
\]
The same argument applies when $\bX$ is sampled from $\bP_2$.
\endgroup

\subsection{Proof of Lemma~\ref{lem:MMPBM-optimal-rate-dev}}
\begingroup
Suppose the true label is $1$. Recall that $\hat z$ is defined in
\eqref{eq:MMPBM-MLE-emp}. On the event $\calE$ in the lemma, choose $0<\tilde\rho=o(1)$ such that
\[
 \max_{m\in\{1,2\}}\max_{\omega\in\calJ_d}
 |\hat\bP_m(\omega)-\bP_m(\omega)|
 \le\tilde\rho\frac{I^*}{d(d+1)/2}.
\]
By definition, on the event $\calE$,
\begin{align*}
\II(\hat z=2)\II(\calE)\leq&\II\left\{\sum_{\omega\in\calJ_d}\bX(\omega)\log\frac{\bP_{2}(\omega)}{\bP_{1}(\omega)}\ZL{\ge}\sum_{\omega\in\calJ_d}\big[\bP_2(\omega)-\bP_1(\omega)\big]-\tilde\rho^{1/3}I^*\right\}\\
&+\II(\calE)\II\Bigg\{\sum_{\omega\in\calJ_d}\bX(\omega)\left[\log\frac{\hat\bP_{2}(\omega)}{\bP_{2}(\omega)}+\log\frac{\bP_{1}(\omega)}{\hat \bP_{1}(\omega)}\right]\\
&\qquad\qquad\qquad\ZL{\ge}\sum_{\omega\in\calJ_d}\Big\{\big[\hat\bP_2(\omega)-\hat\bP_1(\omega)\big]-\big[\bP_2(\omega)-\bP_1(\omega)\big]\Big\}+\tilde\rho^{1/3}I^*\Bigg\},
\end{align*}
The first term is the oracle term. \ZL{Let $s_*\in[0,1]$ be the optimizer in \eqref{eq:def-I-star-mmpbm}. Applying Markov's inequality at $s_*$, as in the proof of Lemma~\ref{lem:oracle-MMPBM}, gives}
\[
\begin{aligned}&\PP\left\{\sum_{\omega\in\calJ_d}\bX(\omega)\log\frac{\bP_{2}(\omega)}{\bP_{1}(\omega)}\ZL{\ge}\sum_{\omega\in\calJ_d}\big[\bP_2(\omega)-\bP_1(\omega)\big]-\tilde\rho^{1/3}I^*\right\}\\
&\quad\leq \ZL{\exp\big(-I^*+s_*\tilde\rho^{1/3}I^*\big)\le\exp\big(-(1-\tilde\rho^{1/3})I^*\big)}.
\end{aligned}
\]

It remains to control the perturbation term.
On $\calE$,
\[
\begin{aligned}&\Bigg|\sum_{\omega\in\calJ_d}\Big\{\big[\hat\bP_2(\omega)-\hat\bP_1(\omega)\big]-\big[\bP_2(\omega)-\bP_1(\omega)\big]\Big\}\Bigg|\\
&\quad\le \sum_{m\in\{1,2\}}\sum_{\omega\in\calJ_d}|\hat\bP_m(\omega)-\bP_m(\omega)|\lesssim \tilde\rho I^*\le \frac{\tilde\rho^{1/3} I^*}{8}
\end{aligned}
\]
for all sufficiently large $n$. Since all intensities in
\eqref{eq:parameter-space-MMPBM} lie between $\theta_0$ and
$\gamma\theta_0$, \ZL{each summand in the Chernoff objective
\eqref{eq:def-I-star-mmpbm} is at most $\gamma\theta_0$ for every
$s\in[0,1]$}. Consequently,
\[
\frac{I^*}{d(d+1)/2}\lesssim \theta_0\asymp \min_{\omega\in\calJ_d}\min_{m\in\{1,2\}}\bP_m(\omega).
\]
Hence, on $\calE$ and for all sufficiently large $n$,
\[
|\hat\bP_m(\omega)-\bP_m(\omega)|\le \frac{1}{2}\bP_m(\omega),\qquad \omega\in\calJ_d,\ m=1,2.
\]
Therefore $\hat\bP_m(\omega)\asymp \bP_m(\omega)$ uniformly in $\omega$, and for every such $\omega$,
\[
\left|\log\frac{\hat\bP_m(\omega)}{\bP_m(\omega)}\right|\le C\tilde\rho\cdot \frac{I^*}{d(d+1)\bP_1(\omega)},\qquad m=1,2,
\]
after using $\bP_1(\omega)\asymp \bP_2(\omega)$ uniformly over $\omega$ under \eqref{eq:parameter-space-MMPBM}.  Consequently, on $\calE$, the second indicator above is bounded by
\[
\II\left[\sum_{\omega\in\calJ_d}\bX(\omega)\frac{C\tilde\rho I^*}{d(d+1)\bP_1(\omega)}\ZL{\ge}\frac{\tilde\rho^{1/3} I^*}{4}\right].
\]
It suffices to bound the probability of the last event. By the moment
generating function of the Poisson distribution, we have
\begin{align*}
 &\log\EE\exp\left(
 \bX(\omega)\frac{C\tilde\rho^{1/3}I^*}
 {d(d+1)\bP_1(\omega)}\right)\\
 &\quad=\bP_1(\omega)\left\{
 \exp\left(\frac{C\tilde\rho^{1/3}I^*}
 {d(d+1)\bP_1(\omega)}\right)-1\right\}\\
 &\quad\le\frac{2C\tilde\rho^{1/3}I^*}{d(d+1)},
\end{align*}
where the last inequality follows from
$\tilde\rho^{1/3}I^*/[d(d+1)\bP_1(\omega)]=o(1)$.
By Chernoff's inequality and independence, we obtain
\begin{align*}
&\PP\left[\sum_{\omega\in\calJ_d}\bX(\omega)
 \frac{C\tilde\rho I^*}{d(d+1)\bP_1(\omega)}
 \ZL{\ge}\frac{\tilde\rho^{1/3}I^*}4\right]\\
&\quad\le \exp\left(-\frac{\tilde\rho^{-1/3}I^*}4\right)
 \prod_{\omega\in\calJ_d}\EE\exp\left(
 \bX(\omega)\frac{C\tilde\rho^{1/3}I^*}
 {d(d+1)\bP_1(\omega)}\right)\\
&\quad\le \exp\left(-\frac{\tilde\rho^{-1/3}I^*}4
                         +C\tilde\rho^{1/3}I^*\right)\\
&\quad\le \exp(-c\tilde\rho^{-1/3}I^*)=o(e^{-I^*})
\end{align*}
for all sufficiently large $n$, since $\tilde\rho=o(1)$ and $I^*\to\infty$.  This yields
\[
\PP(\hat z=2,\calE)\leq \ZL{\exp\big(-[1-o(1)]I^*\big)}.
\]
The same argument applies when $\bX$ is sampled from $\bP_2$.
\endgroup

\ZL{
{\color{black}\subsection{Proofs of Theorems~\ref{thm:minimax-lower-bound} and~\ref{thm:lower-MMPBM}}
\label{sec:repair-lower}}

{\color{black}Fix the two component matrices and let $f_1,f_2$ be their
layer mass functions.} {\color{black}The parameter spaces require
$0<\bP_m(\omega)<1$ for Bernoulli observations and
$\bP_m(\omega)>0$ for Poisson observations, so $f_1,f_2>0$ on
$\mathcal X_d$.} Write
\begin{align}
 {\color{black}I_{1/2}}&\mathrel{\color{black}=}{\color{black}-2\log\sum_{\bu}\sqrt{f_1(\bu)f_2(\bu)},}\label{eq:repair-renyi}\\
 \psi(s)&=\log\sum_{\bu} f_1(\bu)^{1-s}f_2(\bu)^s,
 \qquad I^*=-\min_{0\leq s\leq1}\psi(s).
 \label{eq:repair-chernoff}
\end{align}
Here $I^*$ is the Chernoff information of a complete layer, and $I_{1/2}$ is its R\'enyi-$1/2$ divergence. We have
\begin{equation}\label{eq:repair-info-comparison}
 \frac{I_{1/2}}{2}\leq I^*\leq I_{1/2}.
\end{equation}
The first inequality follows by taking $s=1/2$ in
\eqref{eq:repair-chernoff}. To show the second inequality, let $s_*$ minimize
$\psi$. By convexity and $\psi(0)=\psi(1)=0$, we have
\[
 -\psi(1/2)\ge\frac{I^*}{2\max\{s_*,1-s_*\}}\ge\frac{I^*}{2}.
\]
Together with $I_{1/2}=-2\psi(1/2)$, this proves the second inequality.

\subsubsection{Change of measure}

{\color{black}Define the log likelihood ratio
$L:\mathcal X_d\to\mathbb R$ by $L(\bu)=\log(f_2(\bu)/f_1(\bu))$.}
For $s\in[0,1]$, define
\begin{equation}\label{eq:repair-all-tilts}
 Q_s(\bu)=\frac{f_1(\bu)^{1-s}f_2(\bu)^s}
 {\sum_{\bv\in\mathcal X_d}f_1(\bv)^{1-s}f_2(\bv)^s}
 =e^{-\psi(s)}f_1(\bu)^{1-s}f_2(\bu)^s,
 \qquad \bu\in\mathcal X_d.
\end{equation}
Let $s_*$ minimize $\psi$ on $[0,1]$, and write
\begin{equation}\label{eq:repair-tilt}
 Q(\bu)=Q_{s_*}(\bu)=e^{I^*} f_1(\bu)^{1-s_*}f_2(\bu)^{s_*},
 \qquad V=\operatorname{Var}_Q(L(\bx)).
\end{equation}

For Bernoulli observations, direct calculation gives
\begin{align*}
 \psi'(s)
 &=\frac{\sum_{\bu}f_1(\bu)^{1-s}f_2(\bu)^sL(\bu)}{e^{\psi(s)}},\\
 \psi''(s)
 &=\frac{\sum_{\bu}f_1(\bu)^{1-s}f_2(\bu)^s[L(\bu)-\psi'(s)]^2}
 {e^{\psi(s)}}\ge0.
\end{align*}
For Poisson observations, by \eqref{eq:def-I-star-mmpbm},
we obtain
\begin{align*}
 \psi(s)
 &=\sum_{\omega\in\calJ_d}\left[
 -(1-s)\bP_1(\omega)-s\bP_2(\omega)
 +\bP_1(\omega)^{1-s}\bP_2(\omega)^s\right],\\
 \psi'(s)
 &=\sum_{\omega\in\calJ_d}\left[
 \bP_1(\omega)-\bP_2(\omega)
 +\bP_1(\omega)^{1-s}\bP_2(\omega)^s
 \log\frac{\bP_2(\omega)}{\bP_1(\omega)}\right],\\
 \psi''(s)
 &=\sum_{\omega\in\calJ_d}
 \bP_1(\omega)^{1-s}\bP_2(\omega)^s
 \log^2\frac{\bP_2(\omega)}{\bP_1(\omega)}.
\end{align*}
By \eqref{eq:testing-poisson-tilt}, the edges are independent under
$Q_s$ with means and variances $\bP_1(\omega)^{1-s}\bP_2(\omega)^s$.
Together with
{\color{black}\[
 L(\bx)=\sum_{\omega\in\calJ_d}\left[
 x_\omega\log\frac{\bP_2(\omega)}{\bP_1(\omega)}
 -\bP_2(\omega)+\bP_1(\omega)\right]
\]}
this implies $\psi'(s)=\mathbb E_{Q_s}L$ and
$\psi''(s)=\operatorname{Var}_{Q_s}L$. Moreover,
\begin{align*}
 \psi'(0)&=\sum_{\bu}f_1(\bu)\log(f_2(\bu)/f_1(\bu))
 =-D_{\rm KL}(f_1\Vert f_2)\le0,\\
 \psi'(1)&=D_{\rm KL}(f_2\Vert f_1)\ge0.
\end{align*}
Hence $s_*$ can be taken in $(0,1)$, and
\begin{equation}\label{eq:repair-centered}
 \mathbb E_Q L
 =e^{I^*}\sum_{\bu}f_1(\bu)^{1-s_*}f_2(\bu)^{s_*}L(\bu)
 =\psi'(s_*)=0.
\end{equation}

{\color{black}The equal-prior Bayes error is
\[
 e_B=\inf_{\varphi:\mathcal X_d\to\{1,2\}}
 \frac{\mathbb P_1[\varphi(\bx)=2]+\mathbb P_2[\varphi(\bx)=1]}2.
\]
The infimum is attained by
$\varphi(\bu)=1+\mathbf1[f_2(\bu)>f_1(\bu)]$.}

\begin{lemma}\label{lem:repair-testing}
{\color{black}The equal-prior Bayes error}
\[
 {\color{black}e_B=\frac12\sum_{\bu}\min\{f_1(\bu),f_2(\bu)\}}
\]
satisfies
\begin{equation}\label{eq:repair-testing-bounds}
 \frac38 e^{-I^*-2\sqrt V}\leq e_B\leq e^{-I^*}.
\end{equation}
Moreover, the likelihood-ratio rule with threshold zero has each of
its two error probabilities at most $e^{-I^*}$. In particular, if
$I^*\to\infty$ and $V=o((I^*)^2)$, then
$e_B=\exp(-[1+o(1)]I^*)$.
\end{lemma}

\begin{proof}
The likelihood ratios relative to $Q$ are
\[
 \frac{f_1(\bu)}{Q(\bu)}=e^{-I^*-s_*L(\bu)},\qquad
 \frac{f_2(\bu)}{Q(\bu)}=e^{-I^*+(1-s_*)L(\bu)},\qquad \bu\in\mathcal X_d.
\]
By \eqref{eq:repair-centered} and Chebyshev's inequality,
$Q(|L|\le2\sqrt V)\ge3/4$. Since
$\min\{-s_*L,(1-s_*)L\}\ge-|L|$, it follows that
\begin{align*}
 e_B
 &=\frac12\mathbb E_Q
       \min\{e^{-I^*-s_*L},e^{-I^*+(1-s_*)L}\}\\
 &\geq\frac12e^{-I^*-2\sqrt V}
       Q(|L|\leq2\sqrt V)
 \geq\frac38e^{-I^*-2\sqrt V}.
\end{align*}
{\color{black}For the upper bound, Markov's inequality gives}
\[
 {\color{black}\mathbb P_1(L\mathrel{\ZL{\geq}}0)\leq\mathbb E_1e^{\ZL{s_*}L}=\ZL{e^{-I^*}},
 \qquad
 \mathbb P_2(L\mathrel{\ZL{\leq}}0)\leq\mathbb E_2e^{-(\ZL{1-s_*})L}=\ZL{e^{-I^*}}.}
\]
This proves \eqref{eq:repair-testing-bounds}. The last assertion follows
from $\sqrt V=o(I^*)$ and $I^*\to\infty$.
\end{proof}

\begin{lemma}
\label{lem:repair-variance}
Suppose one of the following holds, with fixed constants
$\Gamma\geq1$ and $p_0<1$:
\begin{enumerate}
 \item Bernoulli edges:
 $\bP_2(\omega)/\bP_1(\omega)\in[\Gamma^{-1},\Gamma]$ and
 $0<\bP_m(\omega)\leq p_0$ for all $m,\omega$;
 \item Poisson edges:
 $\bP_2(\omega)/\bP_1(\omega)
       \in[\Gamma^{-1},\Gamma]$ for all $\omega$.
\end{enumerate}
Then the variance of $L$ under every tilt
$Q_s\propto f_1^{1-s}f_2^s$, $0\leq s\leq1$, is at most $C I_{1/2}$,
where $C$ depends only on the displayed fixed constants.
In particular, $V\leq2C I^*$.
\end{lemma}

\begin{proof}
{\color{black}To verify independence under the tilt, use the product form of
$f_m$ and sum over the coordinates separately:
\begin{align*}
 \sum_{\bu} f_1(\bu)^{1-\ZL{s}}f_2(\bu)^{\ZL{s}}
 &=\prod_{\omega\in\calJ_d}\sum_k
 \mathbb P_1[x_\omega=k]^{1-\ZL{s}}
 \mathbb P_2[x_\omega=k]^{\ZL{s}},\\
 Q_s(\bu)
 &=\prod_{\omega\in\calJ_d}
 \frac{\mathbb P_1[x_\omega=u_\omega]^{1-\ZL{s}}
       \mathbb P_2[x_\omega=u_\omega]^{\ZL{s}}}
 {\sum_k\mathbb P_1[x_\omega=k]^{1-\ZL{s}}
         \mathbb P_2[x_\omega=k]^{\ZL{s}}}.
\end{align*}
Here $k\in\{0,1\}$ for Bernoulli observations and $k\in\mathbb N_0$
for Poisson observations. Thus the coordinates remain independent
under $Q_s$, so the variances of their contributions to $L$ add.}

For Bernoulli observations, the tilted success probability is
\[
 Q_s[x_\omega=1]
 =\frac{\bP_1(\omega)^{1-s}\bP_2(\omega)^s}
 {\bP_1(\omega)^{1-s}\bP_2(\omega)^s
  +[1-\bP_1(\omega)]^{1-s}[1-\bP_2(\omega)]^s}.
\]
Since $\bP_1(\omega)^{1-s}\bP_2(\omega)^s\le\Gamma\bP_1(\omega)$
and $[1-\bP_1(\omega)]^{1-s}[1-\bP_2(\omega)]^s\ge1-p_0$, we have
\[
 Q_s[x_\omega=1]\le\frac{\Gamma\bP_1(\omega)}{1-p_0},
 \qquad s\in[0,1].
\]
Moreover,
{\color{black}\[
 L(\bx)=\sum_{\omega\in\calJ_d}\left\{
 x_\omega\log\frac{\bP_2(\omega)[1-\bP_1(\omega)]}
 {\bP_1(\omega)[1-\bP_2(\omega)]}
 +\log\frac{1-\bP_2(\omega)}{1-\bP_1(\omega)}\right\}.
\]}
The second term is deterministic. By the mean-value theorem applied
to $\log x$ and $\log(1-x)$,
\begin{align*}
 &\left|\log\frac{\bP_2(\omega)[1-\bP_1(\omega)]}
 {\bP_1(\omega)[1-\bP_2(\omega)]}\right|\\
 &\quad\le\frac{|\bP_2(\omega)-\bP_1(\omega)|}
 {\min\{\bP_1(\omega),\bP_2(\omega)\}}
 +\frac{|\bP_2(\omega)-\bP_1(\omega)|}{1-p_0}\\
 &\quad\le\left(\Gamma+\frac1{1-p_0}\right)
 \frac{|\bP_2(\omega)-\bP_1(\omega)|}{\bP_1(\omega)}.
\end{align*}
{\color{black}Since the variance of a Bernoulli variable is at most its success
probability, it follows that}
\begin{align*}
 {\color{black}\operatorname{Var}_{Q_{\ZL{s}}}L}
 &\mathrel{\color{black}=}{\color{black}\sum_{\omega\in\calJ_d}\operatorname{Var}_{Q_{\ZL{s}}}(x_\omega)
 \log^2\frac{\bP_2(\omega)[1-\bP_1(\omega)]}
 {\bP_1(\omega)[1-\bP_2(\omega)]}}\\
 &\le C\sum_{\omega\in\calJ_d}
 \frac{[\bP_2(\omega)-\bP_1(\omega)]^2}{\bP_1(\omega)}.
\end{align*}
By $-\log x\ge1-x$, we have
\begin{align*}
 &-2\log\left(\sqrt{\bP_1(\omega)\bP_2(\omega)}
       +\sqrt{[1-\bP_1(\omega)][1-\bP_2(\omega)]}\right)\\
 &\quad\ge2-2\sqrt{\bP_1(\omega)\bP_2(\omega)}
       -2\sqrt{[1-\bP_1(\omega)][1-\bP_2(\omega)]}\\
 &\quad=[\sqrt{\bP_2(\omega)}-\sqrt{\bP_1(\omega)}]^2
       +[\sqrt{1-\bP_2(\omega)}-\sqrt{1-\bP_1(\omega)}]^2\\
 &\quad\ge\frac{[\bP_2(\omega)-\bP_1(\omega)]^2}
 {[\sqrt{\bP_2(\omega)}+\sqrt{\bP_1(\omega)}]^2}\\
 &\quad\ge\frac{[\bP_2(\omega)-\bP_1(\omega)]^2}
 {(1+\sqrt\Gamma)^2\bP_1(\omega)}.
\end{align*}
Summing over $\omega$ gives
$\operatorname{Var}_{Q_s}L\le CI_{1/2}$, proving the Bernoulli case.

{\color{black}For Poisson observations, the tilted mass of edge $\omega$ is}
\begin{equation}\label{eq:testing-poisson-tilt}
 {\color{black}Q_{\ZL{s}}[x_\omega=k]
 =\frac{\exp(-\ZL{\bP_1(\omega)^{1-s}\bP_2(\omega)^s})
 [\ZL{\bP_1(\omega)^{1-s}\bP_2(\omega)^s}]^k}{k!},
 \qquad k=0,1,\ldots.}
\end{equation}
{\color{black}Indeed, the unnormalized mass is proportional to
$[\ZL{\bP_1(\omega)^{1-s}\bP_2(\omega)^s}]^k/k!$, and summing over
$k$ gives $\exp(\ZL{\bP_1(\omega)^{1-s}\bP_2(\omega)^s})$.
The tilted mean and variance are therefore both
$\ZL{\bP_1(\omega)^{1-s}\bP_2(\omega)^s}$. Moreover,}
{\color{black}\[
 L(\bx)=\sum_{\omega\in\calJ_d}\left[
 x_\omega\log\frac{\bP_2(\omega)}{\bP_1(\omega)}
 -\bP_2(\omega)+\bP_1(\omega)\right],
\]}
{\color{black}so independence gives
\[
 \operatorname{Var}_{Q_{\ZL{s}}}L
 =\sum_{\omega\in\calJ_d}\ZL{\bP_1(\omega)^{1-s}\bP_2(\omega)^s}
 \log^2\frac{\bP_2(\omega)}{\bP_1(\omega)}.
\]}
{\color{black}To bound each summand, note that
$\ZL{\bP_1(\omega)^{1-s}\bP_2(\omega)^s\le\Gamma\bP_1(\omega)}$
and
\[
 \left|\log\frac{\bP_2(\omega)}{\bP_1(\omega)}\right|
 \le\frac{|\bP_2(\omega)-\bP_1(\omega)|}
 {\min\{\bP_1(\omega),\bP_2(\omega)\}}
 \le\Gamma\frac{|\bP_2(\omega)-\bP_1(\omega)|}{\bP_1(\omega)}.
\]
Consequently,
\begin{align*}
 &\ZL{\bP_1(\omega)^{1-s}\bP_2(\omega)^s}
 \log^2\frac{\bP_2(\omega)}{\bP_1(\omega)}\\
 &\quad\le\Gamma^3\frac{[\bP_2(\omega)-\bP_1(\omega)]^2}{\bP_1(\omega)}\\
 &\quad=\Gamma^3
 \frac{[\sqrt{\bP_2(\omega)}+\sqrt{\bP_1(\omega)}]^2}{\bP_1(\omega)}
 [\sqrt{\bP_2(\omega)}-\sqrt{\bP_1(\omega)}]^2\\
 &\quad\le\Gamma^3(1+\sqrt\Gamma)^2
 [\sqrt{\bP_2(\omega)}-\sqrt{\bP_1(\omega)}]^2.
\end{align*}}
Summing over $\omega$ proves the Poisson case. Taking $s=s_*$ and
applying \eqref{eq:repair-info-comparison}, we obtain $V\le2CI^*$.
\end{proof}

\paragraph{Bernoulli probabilities without a bound away from one.}
Let $f_1,f_2$ be product Bernoulli mass functions, with
$p_\omega,q_\omega\in(0,1)$ and
$q_\omega/p_\omega\in[\Gamma^{-1},\Gamma]$ for a fixed $\Gamma$.
Write $I^*=-\min_{s\in[0,1]}\psi(s)$, where
$\psi(s)=\log\sum_{\bu}f_1(\bu)^{1-s}f_2(\bu)^s$.
If $I^*\to\infty$, then the equal-prior Bayes error satisfies
\begin{equation}\label{eq:testing-bernoulli-no-boundary}
 \frac38\exp(-I^*-C(I^*)^{3/4})
 \le e_B\le e^{-I^*}
\end{equation}
for all sufficiently large $I^*$, with $C$ depending only on $\Gamma$.

\begin{proof}
First, for $0<\varepsilon\le1/2$, $s\in[\varepsilon,1-\varepsilon]$
and $u,v>0$, we have
\begin{equation}\label{eq:testing-complement-scalar}
 u^{1-s}v^s\log^2(v/u)
 \le C\varepsilon^{-2}(\sqrt u-\sqrt v)^2.
\end{equation}
Without loss of generality, assume $v\ge u$ and write $v/u=e^x$.
Since $1-e^{-x/2}\ge c\min\{x,1\}$ and
$\sup_{x\ge0}x^2e^{-\varepsilon x}\le C\varepsilon^{-2}$, we have
\begin{align*}
 x^2e^{-(1-s)x}
 &\le C\varepsilon^{-2}\min\{x^2,1\}
 \le C\varepsilon^{-2}(1-e^{-x/2})^2,\\
 u^{1-s}v^s\log^2(v/u)
 &=ue^x x^2e^{-(1-s)x}
 \le C\varepsilon^{-2}u(e^{x/2}-1)^2,
\end{align*}
which proves \eqref{eq:testing-complement-scalar}.

Let $Q_s=e^{-\psi(s)}f_1^{1-s}f_2^s$ and $L=\log(f_2/f_1)$.
Under $Q_s$ the coordinates remain independent. We shall show that
\begin{equation}\label{eq:testing-interior-variance}
 \operatorname{Var}_{Q_s}L\le C\varepsilon^{-2}I_{1/2}
 \le2C\varepsilon^{-2}I^*,
 \qquad s\in[\varepsilon,1-\varepsilon].
\end{equation}
For a fixed edge write $p=p_\omega$ and $q=q_\omega$.
If $\max\{p,q\}\le1/2$, then $Q_s(X=1)\le2\Gamma p$ and
$|\log(q(1-p)/(p(1-q)))|\le C|q-p|/p$. Hence
\begin{align*}
 \operatorname{Var}_{Q_s}\left(X\log\frac qp+(1-X)\log\frac{1-q}{1-p}\right)
 &\le Cp\frac{(q-p)^2}{p^2}
 \le C(\sqrt p-\sqrt q)^2\\
 &\le -2C\log(\sqrt{pq}+\sqrt{(1-p)(1-q)}),
\end{align*}
where the second inequality uses $q/p\in[\Gamma^{-1},\Gamma]$.

If $\max\{p,q\}>1/2$, bounded ratios give
$\min\{p,q\}\ge1/(2\Gamma)$. The denominator of the tilted
Bernoulli probabilities is at least $p^{1-s}q^s\ge1/(2\Gamma)$.
It follows that
\begin{align*}
 \operatorname{Var}_{Q_s}\left(X\log\frac qp+(1-X)\log\frac{1-q}{1-p}\right)
 &=Q_s(X=1)Q_s(X=0)\log^2\frac{q(1-p)}{p(1-q)}\\
 &\le 2\log^2\frac qp
       +4\Gamma(1-p)^{1-s}(1-q)^s\log^2\frac{1-q}{1-p}\\
 &\le C(\sqrt p-\sqrt q)^2
   +C\varepsilon^{-2}(\sqrt{1-p}-\sqrt{1-q})^2.
\end{align*}
where the last inequality follows from the mean-value theorem on
$[1/(2\Gamma),1]$ and \eqref{eq:testing-complement-scalar}. Since
\[
 (\sqrt p-\sqrt q)^2+(\sqrt{1-p}-\sqrt{1-q})^2
 \le-2\log(\sqrt{pq}+\sqrt{(1-p)(1-q)}),
\]
summing over the independent edges proves
\eqref{eq:testing-interior-variance}.

Take $\varepsilon=(I^*)^{-1/4}$ for sufficiently large $I^*$,
and choose $s_*$ to minimize $\psi$ on $[0,1]$.
Put $s=\max\{\varepsilon,\min\{s_*,1-\varepsilon\}\}$.
By convexity and $\psi(0)=\psi(1)=0$, we have
\[
 \psi'(\varepsilon)\le\frac{-\psi(\varepsilon)}{1-\varepsilon}
 \le\frac{I^*}{1-\varepsilon},\qquad
 -\psi'(1-\varepsilon)\le\frac{-\psi(1-\varepsilon)}{1-\varepsilon}
 \le\frac{I^*}{1-\varepsilon}.
\]
Since $\psi'(s_*)=0$, the choice of $s$ gives
\begin{align*}
 |\psi'(s)|&\le\max\{0,\psi'(\varepsilon),-\psi'(1-\varepsilon)\}
 \le\frac{I^*}{1-\varepsilon},\\
 \min\{-s\psi'(s),(1-s)\psi'(s)\}&=-\varepsilon|\psi'(s)|
 \ge-\frac{\varepsilon I^*}{1-\varepsilon}.
\end{align*}
Since $\mathbb E_{Q_s}L=\psi'(s)$, Chebyshev's inequality gives
\[
 Q_s\big(|L-\psi'(s)|\le2\sqrt{\operatorname{Var}_{Q_s}L}\big)\ge3/4.
\]
On this event, we have
\begin{align*}
 \min\{-sL,(1-s)L\}
 &\ge\min\{-s\psi'(s),(1-s)\psi'(s)\}-|L-\psi'(s)|\\
 &\ge-\frac{\varepsilon I^*}{1-\varepsilon}-2\sqrt{\operatorname{Var}_{Q_s}L}.
\end{align*}
The density ratios satisfy
$f_1/Q_s=e^{\psi(s)-sL}$ and $f_2/Q_s=e^{\psi(s)+(1-s)L}$.
Since $\psi(s)\ge-I^*$,
\begin{align*}
 e_B&=\frac12\mathbb E_{Q_s}
       \min\{e^{\psi(s)-sL},e^{\psi(s)+(1-s)L}\}\\
 &\ge\frac38\exp\left(-I^*
       -\frac{\varepsilon I^*}{1-\varepsilon}
       -2\sqrt{2C\varepsilon^{-2}I^*}\right)\\
 &\ge\frac38\exp(-I^*-C(I^*)^{3/4}).
\end{align*}
Together with the Chernoff upper bound, this completes the proof.
\end{proof}

For arbitrary $I^*$, the equal-prior testing error is at least
\[
 \frac14 e^{-I_{1/2}}\ge\frac14 e^{-2I^*}.
\]
Here $\|f_1-f_2\|_{\rm TV}:=\frac12\sum_{\bu}|f_1(\bu)-f_2(\bu)|$
is the total variation distance. By the Cauchy--Schwarz inequality,
\begin{align*}
 \|f_1-f_2\|_{\rm TV}
 &=\frac12\sum_{\bu}|\sqrt{f_1(\bu)}-\sqrt{f_2(\bu)}|
                    [\sqrt{f_1(\bu)}+\sqrt{f_2(\bu)}]\\
 &\le\frac12\left\{\sum_{\bu}[\sqrt{f_1(\bu)}-\sqrt{f_2(\bu)}]^2\right\}^{1/2}
             \left\{\sum_{\bu}[\sqrt{f_1(\bu)}+\sqrt{f_2(\bu)}]^2\right\}^{1/2}\\
 &=\sqrt{1-\left[\sum_{\bu}\sqrt{f_1(\bu)f_2(\bu)}\right]^2}
 =\sqrt{1-e^{-I_{1/2}}}.
\end{align*}
Since $e_B=(1-\|f_1-f_2\|_{\rm TV})/2$, we obtain
\[
 e_B\ge\frac{1-\sqrt{1-e^{-I_{1/2}}}}2
 =\frac{e^{-I_{1/2}}}{2(1+\sqrt{1-e^{-I_{1/2}}})}
 \ge\frac14e^{-I_{1/2}}.
\]
The anchored reduction below therefore also proves necessity of $I^*\to\infty$ for uniform label consistency when $\alpha>1$ is fixed.

\subsubsection{Minimax lower bound for fixed components}

{\color{black}For fixed $f_1,f_2$, define}
\[
 {\color{black}\mathcal Z_\alpha=
 \left\{\bz\in\{1,2\}^n:
 \frac{n}{2\alpha}\leq|\{i\in[n]:z_i=m\}|
 \leq\frac{\alpha n}{2},\ m=1,2\right\},}
\]
{\color{black}where $[n]=\{1,\ldots,n\}$ and $\alpha\ge1$.
Under $\mathbb P_\bz$, the layers are independent with mass
\[
 \mathbb P_\bz(\bx_1=\bu_1,\ldots,\bx_n=\bu_n)
 =\prod_{i=1}^nf_{z_i}(\bu_i).
\]
Write $\mathbb E_\bz$ for the corresponding expectation and
$\mathfrak S_k$ for the permutations of $[k]$. Define
\[
 h(\bu,\bz)=\frac1n\min_{\pi\in\mathfrak S_2}
    \sum_{i=1}^n\mathbf1[u_i\ne\pi(z_i)],
 \qquad
 R_\alpha=\inf_{\widehat\bz:\mathcal X_d^n\to\{1,2\}^n}
    \sup_{\bz\in\mathcal Z_\alpha}\mathbb E_\bz h(\widehat\bz,\bz),
\]
where $[\pi(\bz)]_i=\pi(z_i)$.}

\begin{theorem}
\label{thm:repair-fixed-pair}
Let $n\to\infty$, let $\alpha>1$ be a fixed constant, and let
$I^*\to\infty$. If $V=o((I^*)^2)$ in \eqref{eq:repair-tilt}, then
\begin{equation}\label{eq:repair-minimax}
 R_\alpha=\exp(-[1+o(1)]I^*).
\end{equation}
Moreover, for a constant $c_\alpha>0$ and all sufficiently
large $n$,
\begin{equation}\label{eq:repair-minimax-finite}
 c_\alpha e^{-I^*-2\sqrt V}\leq R_\alpha\leq e^{-I^*}.
\end{equation}
The variance condition holds under either set of assumptions in
Lemma~\ref{lem:repair-variance}. The Bernoulli assumptions include the
homogeneous sparse regime, and the Poisson assumptions hold under the
bounded-intensity-ratio parameter class.
\end{theorem}

\begin{proof}
Choose a fixed constant
$0<c<\min\{1-\alpha^{-1},\alpha-1,1/4\}$, and set
$m=n-2\lfloor(1-c)n/2\rfloor$. Then $cn\le m<cn+2$, so
$m/n\to c$ and $n-m$ is even. {\color{black}Fix
$\ZL{(n-m)/2}$ anchor labels of each type and let the remaining
$m$ labels range freely over $\{1,2\}$. More precisely, define
\[
 \mathcal Z_0=
 \left\{\bz\in\{1,2\}^n:
 \begin{array}{ll}
 z_i=1,&1\le i\le\ZL{(n-m)/2},\\
 z_i=2,&\ZL{(n-m)/2}<i\le n-m
 \end{array}\right\}.
\]
Each vector $(z_{n-m+1},\ldots,z_n)\in\{1,2\}^m$ has exactly one
extension in $\mathcal Z_0$, so $|\mathcal Z_0|=2^m$.}
Since $m/n\to c<\min\{1-\alpha^{-1},\alpha-1\}$, for all
sufficiently large $n$,
\[
 \frac{n}{2\alpha}\le\frac{n-m}{2}
 \le|\{i:z_i=j\}|\le\frac{n+m}{2}\le\frac{\alpha n}{2},
 \qquad \bz\in\mathcal Z_0,\quad j=1,2.
\]
Thus $\mathcal Z_0\subset\mathcal Z_\alpha$.
{\color{black}For two elements $\bu,\bz$ of this cube, all differences occur
on the $m$ free coordinates. The nontrivial permutation in
$\mathfrak S_2$ exchanges the two labels, so
\begin{align*}
 h(\bu,\bz)
 &=\min\left\{\frac1n\sum_i\mathbf1(u_i\ne z_i),
               1-\frac1n\sum_i\mathbf1(u_i\ne z_i)\right\}\\
 &=\frac1n\sum_{i=n-m+1}^n\mathbf1(u_i\ne z_i),
 \qquad \bu,\bz\in\mathcal Z_0,
\end{align*}
where the last equality follows from $m/n<1/2$.}

For any estimator $\widehat\bz$, let
\[
 \widetilde\bz\in\underset{\bw\in\mathcal Z_0}{\arg\min}
 h(\bw,\widehat\bz).
\]
By the triangle inequality and the definition of $\widetilde\bz$, we have
\[
 h(\widetilde \bz,\bz)
 \le h(\widetilde \bz,\widehat \bz)+h(\widehat \bz,\bz)
 \le2h(\widehat \bz,\bz),\qquad \bz\in\mathcal Z_0.
\]

{\color{black}Under the uniform prior on $\mathcal Z_0$, the free labels
are independent and equally likely to be 1 or 2; conditional on
$\bz$, the observations have law $\mathbb P_\bz$. Write $\mathbb P,\mathbb E$
for this joint law. For $i>n-m$ and $a\in\{1,2\}$, independence gives
\begin{align}
 &\mathbb P(z_i=a,\bx_1=\bu_1,\ldots,\bx_n=\bu_n)\nonumber\\
 &\quad=\frac12f_a(\bu_i)
 \prod_{j=1}^{(n-m)/2}f_1(\bu_j)
 \prod_{j=(n-m)/2+1}^{n-m}f_2(\bu_j)
 \prod_{\substack{n-m<j\le n\\j\ne i}}
              \frac{f_1(\bu_j)+f_2(\bu_j)}2.
 \label{eq:testing-free-coordinate-factorization}
\end{align}
Thus
\[
 \mathbb P(z_i=a\mid\bx_1,\ldots,\bx_n)
 =\frac{f_a(\bx_i)}{f_1(\bx_i)+f_2(\bx_i)},
\]
and every coordinate estimator satisfies
\begin{align*}
 2^{-m}\sum_{\bz\in\mathcal Z_0}
 \mathbb P_\bz(\widetilde z_i\ne z_i)
 &\ge\mathbb E\frac{\min\{f_1(\bx_i),f_2(\bx_i)\}}
                         {f_1(\bx_i)+f_2(\bx_i)}\\
 &=\frac12\sum_{\bu}\min\{f_1(\bu),f_2(\bu)\}=e_B.
\end{align*}}
Combining the preceding inequalities, we obtain
\begin{align*}
 {\color{black}\sup_{\bz\in\mathcal Z_\alpha}\mathbb E_\bz h(\widehat \bz,\bz)}
 &\mathrel{\color{black}\ge}{\color{black}2^{-m}\sum_{\bz\in\mathcal Z_0}\mathbb E_\bz h(\widehat \bz,\bz)}\\
 &\ge\frac12\,2^{-m}\sum_{\bz\in\mathcal Z_0}
                       \mathbb E_\bz h(\widetilde \bz,\bz)\\
 &=\frac1{2n}\sum_{i=n-m+1}^n
       2^{-m}\sum_{\bz\in\mathcal Z_0}
                       \mathbb P_\bz(\widetilde z_i\ne z_i)\\
 &\ge\frac{m}{2n}e_B
 \ge\frac{3m}{16n}e^{-I^*-2\sqrt V}.
\end{align*}
Taking the infimum over $\widehat\bz$ and using $m/n\to c>0$, we obtain
the lower bound for any fixed $c_\alpha<3c/16$.

For the upper bound, the component distributions are fixed and known.
Define the coordinatewise likelihood-ratio rule by
\[
 \widehat z_i(\bu_1,\ldots,\bu_n)
 =1+\mathbf1[f_2(\bu_i)>f_1(\bu_i)],\qquad i\in[n].
\]
By Lemma~\ref{lem:repair-testing} and
{\color{black}the definition of $h$},
\[
 \sup_{\bz\in\mathcal Z_\alpha}\mathbb E_\bz h(\widehat \bz,\bz)
 \le\sup_{\bz\in\mathcal Z_\alpha}\frac1n\sum_{i=1}^n
                  \mathbb P_\bz(\widehat z_i\ne z_i)
 \le e^{-I^*}.
\]
By $I^*\to\infty$ and $\sqrt V=o(I^*)$, \eqref{eq:repair-minimax-finite}
implies \eqref{eq:repair-minimax}. This completes the proof.
\end{proof}

\subsubsection{Exactly balanced classes}

\begin{proposition}
\label{prop:repair-exact-balance}
Let $n\to\infty$ be even and let $\alpha=1$. Suppose the edges are
Bernoulli with probabilities in $(0,1)$ or Poisson with positive intensities,
and $\bP_2(\omega)/\bP_1(\omega)\in[\Gamma^{-1},\Gamma]$ for all $\omega$,
where $\Gamma\ge1$ is fixed. Then, for a fixed constant $c>0$,
\begin{equation}\label{eq:repair-exact-bounds}
 c\exp(-I_{1/2}-O(\sqrt{I_{1/2}}))
 \leq R_1
 \leq\min\left\{e^{-I^*},\frac n2e^{-I_{1/2}}\right\}.
\end{equation}
In particular, if $I_{1/2}\to\infty$ and $\log n=o(I_{1/2})$, then
\begin{equation}\label{eq:repair-exact-rate}
 R_1=\exp(-[1+o(1)]I_{1/2}).
\end{equation}
\end{proposition}

\begin{proof}
For the second upper bound, let $\tau\in\mathfrak S_n$ satisfy
$L(\bx_{\tau(1)})\ge\cdots\ge L(\bx_{\tau(n)})$, and define
\[
 \widehat z_{\tau(j)}=
 \begin{cases}
 2,&1\le j\le n/2,\\
 1,&n/2<j\le n.
 \end{cases}
\]
For independent $\bx\sim f_1$ and $\by\sim f_2$, with joint mass
$f_1(\bu)f_2(\bv)$ on $\mathcal X_d^2$, Markov's inequality gives
\begin{align*}
 \mathbb P[L(\bx)\ge L(\by)]
 &\le\mathbb E\exp\big([L(\bx)-L(\by)]/2\big)\\
 &=\mathbb E_1e^{L(\bx)/2}\,\mathbb E_2e^{-L(\by)/2}\\
 &=\left[\sum_{\bu}\sqrt{f_1(\bu)f_2(\bu)}\right]^2
 =e^{-I_{1/2}}.
\end{align*}
Both the true labels and the sorting rule contain exactly $n/2$
observations in each class. Set
$r=|\{i\in[n]:z_i=1,\widehat z_i=2\}|$. Then
\begin{gather*}
 |\{j:z_j=2,\widehat z_j=1\}|=r,\\
 L(\bx_i)\ge L(\bx_j)\quad\text{if }z_i=1,\ \widehat z_i=2,\ z_j=2,\ \widehat z_j=1,
\end{gather*}
since the sorting rule assigns label 2 to the $n/2$ largest scores.
Thus $\sum_{i:z_i=1}\sum_{j:z_j=2}\mathbf1(L(\bx_i)\ge L(\bx_j))\ge r^2\ge r$,
and hence
\[
 h(\widehat \bz,\bz)\le\frac{2r}{n}
 \le\frac2n\sum_{i:z_i=1}\sum_{j:z_j=2}\mathbf1(L(\bx_i)\ge L(\bx_j)).
\]
Taking expectations gives
\[
 \sup_{\bz\in\mathcal Z_1}\mathbb E_\bz h(\widehat \bz,\bz)
 \le\frac2n\cdot\frac{n^2}{4}e^{-I_{1/2}}
 =\frac n2e^{-I_{1/2}}.
\]
The coordinatewise rule used in Theorem~\ref{thm:repair-fixed-pair}
also applies to $\mathcal Z_1$ and gives $R_1\le e^{-I^*}$.

For the lower bound, take $m=\lfloor n/8\rfloor$ and use pairs
$(2j-1,2j)$, $j\in\{1,\ldots,m\}$. For every
$\bb=(b_1,\ldots,b_m)\in\{0,1\}^m$, define a label vector $\bz(\bb)$ by
\[
 z_i(\bb)=
 \begin{cases}
 1+b_j,&i=2j-1,\quad 1\le j\le m,\\
 2-b_j,&i=2j,\quad 1\le j\le m,\\
 1,&2m<i\le n/2+m,\\
 2,&n/2+m<i\le n.
 \end{cases}
\]
Thus the first $2m$ coordinates form the free pairs, and the remaining
coordinates are anchors. Define the finite parameter set
\[
 \mathcal Z_{\rm pair}=\{\bz(\bb):\bb\in\{0,1\}^m\}.
\]
Each $\bz(\bb)$ has $n/2$ labels of each type. Hence
$\mathcal Z_{\rm pair}\subset\mathcal Z_1$ and $|\mathcal Z_{\rm pair}|=2^m$.
Since two such vectors differ on at most $2m<n/2$ coordinates,
\[
 h(\bz(\bb),\bz(\bb'))=\frac2n\sum_{j=1}^m\mathbf1(b_j\ne b'_j),
 \qquad \bb,\bb'\in\{0,1\}^m.
\]
For any estimator $\widehat\bz$, let
\[
 \widetilde\bb\in\underset{\bb\in\{0,1\}^m}{\arg\min}
 h(\bz(\bb),\widehat\bz).
\]
As in the proof of Theorem~\ref{thm:repair-fixed-pair}, for every
$\bb\in\{0,1\}^m$, we have
\[
 h(\bz(\widetilde \bb),\bz(\bb))\le2h(\widehat \bz,\bz(\bb)).
\]
Let $b_1,\ldots,b_m$ be independent Bernoulli$(1/2)$ variables,
and conditionally sample the layers independently with masses
$f_{z_i(\bb)}$. Under this joint law $\mathbb P$, for $a\in\{0,1\}$,
\begin{align*}
 &\mathbb P(b_j=a,\bx_1=\bu_1,\ldots,\bx_n=\bu_n)\\
 &\quad=\frac12f_{1+a}(\bu_{2j-1})f_{2-a}(\bu_{2j})
 \prod_{\substack{1\le k\le m\\k\ne j}}
 \frac{f_1(\bu_{2k-1})f_2(\bu_{2k})+f_2(\bu_{2k-1})f_1(\bu_{2k})}2\\
 &\qquad\cdot\prod_{i=2m+1}^{n/2+m}f_1(\bu_i)
                  \prod_{i=n/2+m+1}^{n}f_2(\bu_i).
\end{align*}
Hence
\begin{equation}\label{eq:testing-pair-posterior}
 \mathbb P(b_j=a\mid \bx_1,\ldots,\bx_n)
 =\frac{f_{1+a}(\bx_{2j-1})f_{2-a}(\bx_{2j})}
 {f_1(\bx_{2j-1})f_2(\bx_{2j})+f_2(\bx_{2j-1})f_1(\bx_{2j})}.
\end{equation}
Averaging the smaller posterior probability gives
\begin{align}
 &2^{-m}\sum_{\bb\in\{0,1\}^m}
       \mathbb P_{\bz(\bb)}(\widetilde b_j\ne b_j)\nonumber\\
 &\qquad\ge\frac12\sum_{\bu,\bv\in\mathcal X_d}
          \min\{f_1(\bu)f_2(\bv),f_2(\bu)f_1(\bv)\}.
 \label{eq:testing-pair-bayes-reduction}
\end{align}
The right-hand side is the Bayes error for testing
$f_1(\bu)f_2(\bv)$ against $f_2(\bu)f_1(\bv)$. Since
\begin{align*}
 &\sum_{\bu,\bv}[f_1(\bu)f_2(\bv)]^{1-s}[f_2(\bu)f_1(\bv)]^s\\
 &\quad=\left[\sum_{\bu} f_1(\bu)^{1-s}f_2(\bu)^s\right]
         \left[\sum_{\bv} f_1(\bv)^sf_2(\bv)^{1-s}\right]
 =e^{\psi(s)+\psi(1-s)}.
\end{align*}
and $\psi(s)+\psi(1-s)\ge2\psi(1/2)$ by convexity, its Chernoff
information is
\[
 -\min_{0\le s\le1}[\psi(s)+\psi(1-s)]
 =-2\psi(1/2)=I_{1/2}.
\]
The tilted joint mass at $s=1/2$ is
\[
 e^{I_{1/2}}\sqrt{f_1(\bu)f_2(\bu)}\sqrt{f_1(\bv)f_2(\bv)}
 =Q_{1/2}(\bu)Q_{1/2}(\bv).
\]
By independence and the common marginal $Q_{1/2}$,
\begin{align*}
 \mathbb E_{Q_{1/2}\otimes Q_{1/2}}[L(\bx)-L(\by)]
 &=\mathbb E_{Q_{1/2}}L-\mathbb E_{Q_{1/2}}L=0,\\
 \operatorname{Var}_{Q_{1/2}\otimes Q_{1/2}}(L(\bx)-L(\by))
 &=2\operatorname{Var}_{Q_{1/2}}L\le2CI_{1/2}
\end{align*}
by \eqref{eq:testing-interior-variance} with $s=\varepsilon=1/2$
in the Bernoulli case, after enlarging $C$, and by
Lemma~\ref{lem:repair-variance} in the Poisson case. By
Lemma~\ref{lem:repair-testing}, we obtain the lower bound
\[
 \frac38\exp(-I_{1/2}-2\sqrt{2CI_{1/2}})
\]
for the pair Bayes error in \eqref{eq:testing-pair-bayes-reduction}.
Combining the preceding inequalities, for any $\widehat\bz$, we obtain
\begin{align*}
 \sup_{\bz\in\mathcal Z_1}\mathbb E_\bz h(\widehat \bz,\bz)
 &\ge\frac12\,2^{-m}\sum_{\bb\in\{0,1\}^m}
       \mathbb E_{\bz(\bb)}h(\bz(\widetilde \bb),\bz(\bb))\\
 &=\frac1n\sum_{j=1}^m2^{-m}\sum_{\bb\in\{0,1\}^m}
       \mathbb P_{\bz(\bb)}(\widetilde b_j\ne b_j)\\
 &\ge\frac{3m}{8n}\exp(-I_{1/2}-2\sqrt{2CI_{1/2}}).
\end{align*}
Taking the infimum over $\widehat\bz$ and using $m/n\to1/8$ proves
the lower bound. If $I_{1/2}\to\infty$, taking logarithms and using
the sorting upper bound gives
\[
 -1-2\sqrt{\frac{2C}{I_{1/2}}}
    +\frac{\log(3m/(8n))}{I_{1/2}}
 \le\frac{\log R_1}{I_{1/2}}
 \le-1+\frac{\log(n/2)}{I_{1/2}}.
\]
By $m/n\to1/8$ and $\log n=o(I_{1/2})$, we obtain
\eqref{eq:repair-exact-rate}.
\end{proof}

It remains to apply these testing bounds to Theorems~\ref{thm:minimax-lower-bound} and
\ref{thm:lower-MMPBM}. Under \eqref{eq:parameter-space}, each diagonal
block probability lies between $p_m$ and $\gamma p_m$, while each
off-diagonal probability lies between $q_m/\gamma$ and $q_m$.
The cross-component ratio conditions imply
\begin{align*}
 \frac{\max_{m,\omega}\bP_m(\omega)}
      {\min_{m,\omega}\bP_m(\omega)}
 &\le\frac{\max\{\gamma p_1,\gamma p_2,q_1,q_2\}}
           {\min\{p_1,p_2,q_1/\gamma,q_2/\gamma\}}\\
 &\le\gamma^2\frac{\max\{p_1,p_2,q_1,q_2\}}
                       {\min\{p_1,p_2,q_1,q_2\}}
 \le\gamma^3.
\end{align*}
Hence $\gamma^{-3}\le\bP_2(\omega)/\bP_1(\omega)\le\gamma^3$ for all
$\omega\in\calJ_d$. Applying \eqref{eq:testing-bernoulli-no-boundary}
and the proof of Theorem~\ref{thm:repair-fixed-pair}, we obtain
\[
 R_\alpha\ge\frac{m}{2n}e_B
 \ge\frac{3m}{16n}\exp(-I^*-C(I^*)^{3/4}),\qquad m/n\to c>0.
\]
Since $(I^*)^{3/4}=o(I^*)$, this proves the Bernoulli lower bound.
Under \eqref{eq:parameter-space-MMPBM}, all intensities lie in
$[\theta_0,\gamma\theta_0]$, giving the Poisson ratio condition with
$\Gamma=\gamma$. Lemma~\ref{lem:repair-variance} then gives
\[
 V\le2CI^*,\qquad V/(I^*)^2\le2C/I^*\to0.
\]
Theorem~\ref{thm:repair-fixed-pair} proves the Poisson lower bound.
In both cases $\alpha>1$ is fixed. Together with the oracle upper
bounds, these inequalities give the same exponential rate.
}
\ZL{
{\color{black}\subsection{Proofs of Theorems~\ref{thm:main} and~\ref{thm:MMPBM-main}}
\label{fs:section}
\label{fs:node-class}}

Suppose the true node and layer classes are balanced, $\eta_z=o(1)$,
and $K\eta_\sigma=o(1)$.
{\color{black}Let $\mathcal J_d=\{(u,v):1\le u\le v\le d\}$. For the true
community vectors $\bsig_m=(\bsig_m(u))_{u\in[d]}\in[K]^d$ and true symmetric block matrices
$\bB_m$, write
\[
 p_{m\omega}=\bB_m[\bsig_m(u),\bsig_m(v)],
 \qquad m\in[2],\quad \omega=(u,v)\in\mathcal J_d.
\]
The observed entries $\bX_i(\omega)$, $i\in[n]$, $\omega\in\mathcal J_d$,
are independent Bernoulli or Poisson variables with means
$p_{z_i^*\omega}$. Put}
\[
 {\color{black}\bar p=\max_{m,\omega}p_{m\omega}.}
\]
{\color{black}Assume $\bar p/G\le p_{m\omega}\le\bar p$ for fixed $G$.
For Bernoulli observations assume $\bar p\le1-c_0$ for a fixed $c_0\in(0,1)$.} {\color{black}Constants $c,C>0$ may change from line to line and depend on the fixed
balance constants, $G$, and, for Bernoulli observations, $c_0$.
In probability bounds with a prescribed exponent $A$, they may also
depend on $A$, but never on $K,n$ or $d$. The maximum edge mean
used here agrees with the main-text sparsity envelope up to a fixed
homogeneity constant.} Write $I^*$ for the complete-layer Chernoff information and $I_{1/2}$ for the auxiliary R\'enyi-$1/2$ divergence. Then
\begin{equation}\label{fs:information}
 I_{1/2}/2\le I^*\le I_{1/2},\qquad
 I^*\asymp\sum_{\omega\in\mathcal J_d}
              (p_{1\omega}-p_{2\omega})^2/\bar p,
 \qquad I^*\lesssim d^2\bar p.
\end{equation}
By \eqref{eq:repair-renyi}--\eqref{eq:repair-chernoff} and concavity,
$I_{1/2}/2\le I^*\le I_{1/2}$.
{\color{black}For Bernoulli means $p,q\in[\bar p/G,\bar p]$, we have
\begin{align*}
 &2\big[1-\sqrt{pq}-\sqrt{(1-p)(1-q)}\big]\\
 &\quad=(\sqrt p-\sqrt q)^2+(\sqrt{1-p}-\sqrt{1-q})^2\\
 &\quad=(p-q)^2\left[\frac1{(\sqrt p+\sqrt q)^2}
             +\frac1{(\sqrt{1-p}+\sqrt{1-q})^2}\right]
 \asymp\frac{(p-q)^2}{\bar p}.
\end{align*}
Moreover,
\[
 \sqrt{pq}+\sqrt{(1-p)(1-q)}\ge c_0,
 \qquad 1-x\le-\log x\le c_0^{-1}(1-x)\quad(c_0\le x\le1).
\] Summing the resulting bounds gives
\[
 I_{1/2}=-2\sum_\omega
 \log\big(\sqrt{p_{1\omega}p_{2\omega}}
     +\sqrt{(1-p_{1\omega})(1-p_{2\omega})}\big)
 \asymp\sum_\omega\frac{(p_{1\omega}-p_{2\omega})^2}{\bar p}.
\]
For Poisson means the same comparison follows directly from
\[
 I_{1/2}=\sum_\omega
 \frac{(p_{1\omega}-p_{2\omega})^2}
      {(\sqrt{p_{1\omega}}+\sqrt{p_{2\omega}})^2}
 \asymp\sum_\omega\frac{(p_{1\omega}-p_{2\omega})^2}{\bar p}.
\]
Finally, each summand on the right is at most $\bar p$ and
$|\mathcal J_d|=d(d+1)/2$, which gives the last bound in
\eqref{fs:information}.}

{\color{black}Let $\eta_z,\eta_\sigma=o(1)$ bound the layer and node error
fractions. All candidate labels below are aligned with the true labels.
Specifically, the candidate assignments are the vectors
$\boldsymbol{\zeta}\in[2]^n$ and $\boldsymbol{\tau}_m=(\boldsymbol{\tau}_m(u))_{u\in[d]}\in[K]^d$, $m\in[2]$, satisfying
\begin{equation}\label{fs:candidate-assignments}
 \frac1n\sum_{i=1}^n\mathbf1(\zeta_i\ne z_i^*)\le\eta_z,
 \qquad
 \max_{m\in[2]}\frac1d\sum_{u=1}^d
              \mathbf1[\boldsymbol{\tau}_m(u)\ne\bsig_m(u)]\le\eta_\sigma.
\end{equation}}
Let $\widetilde \bz\in[2]^n$ and
$\widetilde\bsig_m=(\widetilde\bsig_m(u))_{u\in[d]}\in[K]^d$ denote the initializer before alignment.
A component permutation $\pi\in\mathfrak S_2$ and community
permutations $\rho_m\in\mathfrak S_K$ give aligned assignments
\[
 \zeta_i=\pi^{-1}(\widetilde z_i),\qquad
 \boldsymbol{\tau}_m(u)=\rho_m(\widetilde\bsig_{\pi(m)}(u)).
\]
On the initialization event, choose $(\pi,\rho_1,\rho_2)$ satisfying
\eqref{fs:candidate-assignments}. Let
$(\widehat\bB_1^{\rm raw},\widehat\bB_2^{\rm raw})$ be the block averages
before alignment. The aligned entries are
\[
 \widehat \bB_m(k,\ell)=
 \widehat\bB^{\rm raw}_{\pi(m)}[\rho_m^{-1}(k),\rho_m^{-1}(\ell)].
\]
{\color{black}Using $\binom da\le(ed/a)^a$, the numbers of candidate
layer assignments and pairs of node assignments are bounded by
\[
 \sum_{a=0}^{\lfloor n\eta_z\rfloor}\binom na
 \le\exp(Cn\eta_z\log(e/\eta_z))
\]
and
\[
 \left[\sum_{a=0}^{\lfloor d\eta_\sigma\rfloor}
       \binom da(K-1)^a\right]^2
 \le\exp(Cd\eta_\sigma\log(eK/\eta_\sigma)),
\]
respectively.}
{\color{black}For a candidate node vector $\boldsymbol{\tau}_m$, define its \ZL{unordered} edge cell by
\[
 A_{mk\ell}:=\{(u,v)\in\mathcal J_d:
  \ZL{(\boldsymbol{\tau}_m(u),\boldsymbol{\tau}_m(v))\in\{(k,\ell),(\ell,k)\}}\},
 \qquad m\in[2],\quad 1\le k\le\ell\le K.
\]
The index $j=(m,k,\ell)$ ranges over $m\in[2]$ and $1\le k\le\ell\le K$. Write $A_j=A_{mk\ell}$, $M_j=|A_j|$, and
$B_j=\bB_m(k,\ell)$, and use the same indexing for fitted block entries.
\ZL{Across the two components there are $K(K+1)$ cells.}
In particular, the cell sizes are
\[
 M_{mk\ell}=
 \begin{cases}
 |\boldsymbol{\tau}_m^{-1}(k)|\,|\boldsymbol{\tau}_m^{-1}(\ell)|,&k<\ell,\\[1mm]
 |\boldsymbol{\tau}_m^{-1}(k)|[|\boldsymbol{\tau}_m^{-1}(k)|+1]/2,&k=\ell.
 \end{cases}
\]
Here $\boldsymbol{\tau}_m^{-1}(k)=\{u\in[d]:\boldsymbol{\tau}_m(u)=k\}$. Since $|\boldsymbol{\tau}_m^{-1}(k)|\asymp d/K$, these formulas give
$M_j\asymp d^2/K^2$.}
The cells from different components can overlap.

{\color{black}For a single observation, the exponential-family mass function is
\[
 f_\theta(x)=\exp(x\theta-b(\theta)+c(x)),\qquad \theta\in\mathbb R.
\]
For Bernoulli observations, $x\in\{0,1\}$,
$\theta(p)=\log(p/(1-p))$ for $p\in(0,1)$,
$b(\theta)=\log(1+e^\theta)$, and $c(x)=0$.
For Poisson observations, $x\in\mathbb Z_{\ge0}$,
$\theta(p)=\log p$ for $p>0$, $b(\theta)=e^\theta$, and
$c(x)=-\log(x!)$.}

{\color{black}\subsubsection{Uniform bounds for the block averages}}

{\color{black}For every candidate $(\boldsymbol{\zeta},\boldsymbol{\tau}_1,\boldsymbol{\tau}_2)$ satisfying
\eqref{fs:candidate-assignments}, the empirical block averages are
\[
 \widehat \bB_m(k,\ell)=
 \frac{\sum_{i:\zeta_i=m}\sum_{\omega\in A_{mk\ell}}\bX_i(\omega)}
 {|\{i:\zeta_i=m\}|M_{mk\ell}},
 \qquad m\in[2],\quad 1\le k\le\ell\le K.
\]
When needed, write $\widehat\bB_m(\boldsymbol\zeta,\boldsymbol\tau)$
to indicate the candidate assignments.}

\begin{lemma}\label{fs:oracle-comparison}
Suppose $I^*\to\infty$, $K^2\log n=o(n)$, and
\[
 \eta_z\vee\eta_\sigma
 =o\!\left(\frac{I^*}{K^2d^2\bar p}\right),\qquad
 \frac{d}{n}\eta_\sigma\log(e/\eta_\sigma)
 =o\!\left(\frac{(I^*)^2}{K^2d^2\bar p}\right).
\]
For the true community assignments, write
\[
 A^*_{mk\ell}
 =\{(u,v)\in\mathcal J_d:
 (\bsig_m(u),\bsig_m(v))\in\{(k,\ell),(\ell,k)\}\},
 \qquad m\in[2],\quad 1\le k\le\ell\le K,
\]
and define the oracle block averages by
\[
 \overline\bB_m(k,\ell)
 =\frac{\sum_{i:z_i^*=m}\sum_{\omega\in A^*_{mk\ell}}\bX_i(\omega)}
 {|\{i:z_i^*=m\}|\,|A^*_{mk\ell}|}.
\]
For every fixed $A>0$, with probability at least $1-C n^{-A}e^{-4I^*}$,
simultaneously for every candidate in \eqref{fs:candidate-assignments},
\begin{equation}\label{fs:candidate-oracle-close}
 \max_{m\in[2],\,k\le\ell}
 |\widehat\bB_m(k,\ell)-\overline\bB_m(k,\ell)|
 =o\!\left(\frac{I^*}{d^2}\right).
\end{equation}
The little-$o$ term is deterministic and uniform over these candidates.
On the same event,
\begin{equation}\label{fs:oracle-relative-close}
 \max_{m\in[2],\,k\le\ell}
 \frac{|\overline\bB_m(k,\ell)-\bB_m(k,\ell)|}{\bar p}
 \le C\left[
 \sqrt{\frac{K^2(\log n+I^*)}{nd^2\bar p}}
 +\frac{K^2(\log n+I^*)}{nd^2\bar p}\right]=o(1).
\end{equation}
Consequently both the oracle and candidate block averages are positive
and uniformly within $1+o(1)$ of the corresponding true block entries.
In the Bernoulli case they are strictly below one.
\end{lemma}

\begin{proof}
Now fix a candidate satisfying \eqref{fs:candidate-assignments}. Define
\[
\begin{aligned}
 I_m&=\{i:\zeta_i=m\},\quad I_m^*=\{i:z_i^*=m\},\\
 S_{mk}&=\{u:\boldsymbol\tau_m(u)=k\},\quad
 S_{mk}^*=\{u:\bsig_m(u)=k\}.
\end{aligned}
\]
Then
\[
 \begin{aligned}
 \sum_{i=1}^n|\mathbf1(i\in I_m)-\mathbf1(i\in I_m^*)|&\le n\eta_z,\\
 \sum_{u=1}^d|\mathbf1(u\in S_{mk})-\mathbf1(u\in S_{mk}^*)|
 &\le d\eta_\sigma.
 \end{aligned}
\]
By \eqref{fs:information}, $I^*\lesssim d^2\bar p$. The first accuracy
condition therefore implies
\[
 \eta_z=o(1),\qquad K\eta_\sigma=o(1).
\]
The balance assumptions now give, uniformly over all candidates,
\begin{equation}\label{fs:oracle-comparison-sizes}
 |I_m|\asymp |I_m^*|\asymp n,\qquad
 |S_{mk}|\asymp|S_{mk}^*|\asymp d/K,\qquad
 |A_{mk\ell}|\asymp|A^*_{mk\ell}|\asymp d^2/K^2.
\end{equation}
By the definition of $A_{mk\ell}$, we have
\begin{equation}\label{fs:changed-edge-count}
\begin{split}
 &\sum_{\omega\in\mathcal J_d}
 |\mathbf1(\omega\in A_{mk\ell})-\mathbf1(\omega\in A^*_{mk\ell})|\\
 &\quad\le
 \left[\sum_{a\in\{k,\ell\}}\sum_{u=1}^d
 |\mathbf1(u\in S_{ma})-\mathbf1(u\in S_{ma}^*)|\right]
 \sum_{b\in\{k,\ell\}}(|S_{mb}|+|S_{mb}^*|)\\
 &\quad\le C\frac{d^2\eta_\sigma}{K}.
\end{split}
\end{equation}
Moreover,
\begin{align}
 &\left|\sum_{i\in I_m}\sum_{\omega\in A_{mk\ell}}\bX_i(\omega)
       -\sum_{i\in I_m^*}\sum_{\omega\in A^*_{mk\ell}}\bX_i(\omega)
   \right|\nonumber\\
 &\quad\le
 \sum_{i=1}^n|\mathbf1(i\in I_m)-\mathbf1(i\in I_m^*)|
       \sum_{\omega\in A_{mk\ell}}\bX_i(\omega)\nonumber\\
 &\qquad+
 \sum_{i\in I_m^*}\sum_{\omega\in\mathcal J_d}
 |\mathbf1(\omega\in A_{mk\ell})-\mathbf1(\omega\in A^*_{mk\ell})|
 \bX_i(\omega).
 \label{fs:changed-observation-sums}
\end{align}
For fixed $I_m,S_{mk},S_{m\ell}$, each sum on the right consists of
independent observations. Their means are at most
\[
 C\frac{nd^2\bar p\eta_z}{K^2}
 \quad\text{and}\quad
 C\frac{nd^2\bar p\eta_\sigma}{K},
\]
respectively. The number of possible sets $I_m$ is at most
\[
 \sum_{a=0}^{\lfloor n\eta_z\rfloor}\binom na
 \le\exp(Cn\eta_z\log(e/\eta_z)).
\]
Similarly, the number of possible sets $S_{mk}$ is at most
\[
 \sum_{a=0}^{\lfloor d\eta_\sigma\rfloor}\binom da
 \le\exp(Cd\eta_\sigma\log(e/\eta_\sigma)).
\]
Since $A_{mk\ell}$ is determined by $S_{mk}$ and $S_{m\ell}$,
the first sum in \eqref{fs:changed-observation-sums} has at most
\[
 \exp(Cn\eta_z\log(e/\eta_z)
             +Cd\eta_\sigma\log(e/\eta_\sigma))
\]
possible index sets, and the second has at most
$\exp(Cd\eta_\sigma\log(e/\eta_\sigma))$.
For a sum $S$ of independent Bernoulli or Poisson observations with
mean at most $\mu$, the moment generating function satisfies
\[
 \log\mathbb E e^{\lambda S}\le\mu(e^\lambda-1),\qquad\lambda\ge0.
\]
Taking $\lambda=\log(x/\mu)$ gives
\begin{equation}\label{fs:changed-sum-bennett}
 \mathbb P(S\ge x)
 \le\exp(-x\log(x/\mu)+x-\mu),\qquad x\ge\mu>0.
\end{equation}
For the first sum in \eqref{fs:changed-observation-sums}, the mean
divided by $nI^*/K^2$ is at most $Cd^2\bar p\eta_z/I^*=o(1)$.
Moreover,
\begin{align*}
 &\frac{n\eta_z\log(e/\eta_z)}
 {(nI^*/K^2)\log(I^*/(Cd^2\bar p\eta_z))}\\
 &\quad=
 \frac{K^2\eta_z}{I^*}
 \frac{\log(eI^*/(d^2\bar p\eta_z))
                   +\log(d^2\bar p/I^*)}
      {\log(I^*/(Cd^2\bar p\eta_z))}\\
 &\quad\le
 C\frac{K^2d^2\bar p\eta_z}{I^*}
       \left(\frac1{d^2\bar p}+\frac1{I^*}\right)
 \to0.
\end{align*}
where we used $\log(I^*/(Cd^2\bar p\eta_z))\to\infty$ and
$\max\{0,\log(d^2\bar p/I^*)\}\le d^2\bar p/I^*$.

For the node sets and the probability bound, we have
\begin{align*}
 \frac{d\eta_\sigma\log(e/\eta_\sigma)}{nI^*/K^2}
 &=\frac{K^2(d/n)\eta_\sigma\log(e/\eta_\sigma)}{I^*}
   =o\!\left(\frac{I^*}{d^2\bar p}\right)=o(1),\\
 \frac{\log n+I^*}{nI^*/K^2}
 &=\frac{K^2\log n}{nI^*}+\frac{K^2}{n}=o(1).
\end{align*}
For the second sum in \eqref{fs:changed-observation-sums}, its mean
divided by $nI^*/K^2$ is at most
$CKd^2\bar p\eta_\sigma/I^*=o(1)$. Its exponent in
\eqref{fs:changed-sum-bennett}, divided by $nI^*/K^2$, thus tends to
infinity. Also, $K^2\log n=o(n)$ gives
$\log(K(K+1))\le C\log n$ for all sufficiently large $n$.
By the preceding bounds, there exists $\rho_n\to0$ such that
\[
 \frac{d^2\bar p\eta_z+Kd^2\bar p\eta_\sigma}{\rho_n I^*}=o(1),
 \qquad
 \frac{K^2\eta_z\log(e/\eta_z)}
 {\rho_n I^*\log(\rho_n I^*/(Cd^2\bar p\eta_z))}=o(1),
\]
and
\[
 \frac{K^2[d\eta_\sigma\log(e/\eta_\sigma)
                         +(A+5)\log n+4I^*]}{\rho_n nI^*}=o(1).
\]
Applying \eqref{fs:changed-sum-bennett} and taking a union bound over
both sums, all admissible sets and all blocks, we obtain, with probability at least
$1-C n^{-A}e^{-4I^*}$,
\begin{equation}\label{fs:oracle-numerator-close}
 \max_{m,k\le\ell}\sup_{\text{candidates in }\eqref{fs:candidate-assignments}}
 \left|\sum_{i\in I_m}\sum_{\omega\in A_{mk\ell}}\bX_i(\omega)
       -\sum_{i\in I_m^*}\sum_{\omega\in A^*_{mk\ell}}\bX_i(\omega)
 \right|
 \le2\rho_n\frac{nI^*}{K^2}.
\end{equation}

By Bernstein's inequality with $t=(A+5)\log n+4I^*$ and a union bound
over the $K(K+1)$ oracle blocks, we obtain
\eqref{fs:oracle-relative-close}. Its right side is $o(1)$, since
\[
 \frac{K^2(\log n+I^*)}{nd^2\bar p}
 \lesssim\frac{K^2\log n}{nI^*}+\frac{K^2}{n}=o(1).
\]
In particular, $\overline\bB_m(k,\ell)\le C\bar p$ uniformly.
It remains to compare the denominators. By
\eqref{fs:oracle-comparison-sizes} and
\eqref{fs:changed-edge-count},
\[
 \frac{\big||I_m|\,|A_{mk\ell}|-|I_m^*|\,|A^*_{mk\ell}|\big|}
 {|I_m|\,|A_{mk\ell}|}
 \le C(\eta_z+K\eta_\sigma).
\]
Thus \eqref{fs:oracle-numerator-close} yields
\begin{align*}
 |\widehat\bB_m(k,\ell)-\overline\bB_m(k,\ell)|
 &\le C\rho_n\frac{I^*}{d^2}
       +C\bar p(\eta_z+K\eta_\sigma)\\
 &=o\!\left(\frac{I^*}{d^2}\right),
\end{align*}
uniformly, because
\[
 \frac{d^2\bar p(\eta_z+K\eta_\sigma)}{I^*}
 =o(K^{-2}+K^{-1})=o(1).
\]
This proves \eqref{fs:candidate-oracle-close}.
Since $I^*\lesssim d^2\bar p$ and
$\bB_m(k,\ell)\ge\bar p/G$, the last assertions follow from
\eqref{fs:oracle-relative-close} and
\eqref{fs:candidate-oracle-close}. For Bernoulli observations we also
use $\bar p\le1-c_0$ and relative consistency of the oracle and candidate means.
\end{proof}

The first accuracy condition in Lemma~\ref{fs:oracle-comparison}
follows from the second condition in \eqref{eq:main-init-cond} or
\eqref{eq:MMPBM-main-init-cond}. For Bernoulli observations,
$d^2\bar p\asymp\bar a(1+d/n)$ and $(d/n)/(1+d/n)\le1$, so the first
condition in \eqref{eq:main-init-cond} gives
\[
 \frac{d}{n}\eta_\sigma\log(e/\eta_\sigma)
 =o\!\left(\frac{(I^*)^2}{K^2d^2\bar p}\right).
\]
For Poisson observations, $\bar p\asymp\theta_0$ and this is exactly
the first condition in \eqref{eq:MMPBM-main-init-cond}, up to fixed
constants. Replacing $\log(1/\eta_\sigma)$ by
$\log(e/\eta_\sigma)$ leaves the condition unchanged.

Furthermore, $\eta_\sigma=o(1/K^2)$ implies
\[
 \log(eK/\eta_\sigma)\asymp\log(e/\eta_\sigma).
\]
In the Poisson model, the bound
$x\log(eK/x)\le CK\sqrt{x}$ for $0\le x\le1$ and $K\ge1$ gives
\[
 d\eta_\sigma\log(eK/\eta_\sigma)
 \le CdK\sqrt{\eta_\sigma}
 =o\!\left(\sqrt{I^*/\theta_0}\right)=o(I^*),
\]
using $\theta_0\ge c>0$ and $I^*\to\infty$.

\subsubsection{Likelihood comparison with the oracle block averages}

Write $M^*_{mk\ell}=|A^*_{mk\ell}|$ and
$n_m=|\{i:z_i^*=m\}|$. The oracle matrices $\overline\bB_m$
are those defined in Lemma~\ref{fs:oracle-comparison}.
Their true cell sizes satisfy $M^*_{mk\ell}\asymp d^2/K^2$.

Fix a layer $i$, and write $m=z_i^*$. For $a\in[2]$, put
\[
 \overline\bB_a^{(-i)}(k,\ell):=
 \frac{\sum_{j:z_j^*=a,\,j\ne i}
             \sum_{\omega\in A^*_{ak\ell}}\bX_j(\omega)}
 {[n_a-\mathbf1(a=m)]M^*_{ak\ell}}.
\]
Then $\overline\bB_a^{(-i)}$, $a\in[2]$, are independent of $\bX_i$,
and $\overline\bB_{3-m}^{(-i)}=\overline\bB_{3-m}$. Define
\[
\begin{aligned}
 \mathcal F_i&:=\bigcap_{a,k\le\ell}
 \left\{\left|\overline\bB_a^{(-i)}(k,\ell)-\bB_a(k,\ell)\right|
                  \le\tfrac{c_0}{4} \bB_a(k,\ell)\right\},
\\
 R_i&:=\sum_{a,k\le\ell}M^*_{ak\ell}
 \frac{[\overline\bB_a^{(-i)}(k,\ell)-\bB_a(k,\ell)]^2}{\bar p}.
\end{aligned}
\]
For Poisson observations, use $c_0=1$ in the definition of $\mathcal F_i$ and in its proof below.

\begin{lemma}\label{fs:oracle-block-mgf}
Suppose $K^2\log n=o(n)$ and $d^2\bar p\to\infty$.
For every fixed $A>0$, for all sufficiently large $n$,
\[
 \mathbb P\left(\bigcup_{i=1}^n\mathcal F_i^c\right)\le n^{-A}e^{-4I^*}.
\]
There are constants $c,C>0$, depending only on the fixed balance and
homogeneity constants and, in the Bernoulli case, on $c_0$, such that
\begin{equation}\label{fs:oracle-quadratic-mgf}
 \mathbb E\left[e^{tR_i}\mathbf1(\mathcal F_i)\right]
       \le\exp(CK^2t/n),\qquad 0\le t\le cn.
\end{equation}
\end{lemma}

\begin{proof}
Fix $a,k,\ell$. The sum defining this average contains
$[n_a-\mathbf1(a=z_i^*)]M^*_{ak\ell}$ independent variables, each with
mean $\bB_a(k,\ell)$. Its centered log moment generating function is at most
\[
 [n_a-\mathbf1(a=z_i^*)]M^*_{ak\ell}\bB_a(k,\ell)(e^s-1-s)
\]
for both observation families.
Consequently, Bernstein's inequality gives
\[
 \mathbb P\left(\left|\overline\bB_a^{(-i)}(k,\ell)-\bB_a(k,\ell)\right|
                       \ge x\right)
 \le2\exp\left(-\frac{[n_a-\mathbf1(a=z_i^*)]M^*_{ak\ell}x^2}
 {2\bB_a(k,\ell)+2x/3}\right),
 \qquad x>0.
\]
For $0<x\le c_0\bB_a(k,\ell)/4$, $n_a-\mathbf1(a=z_i^*)\asymp n$ and
$\bB_a(k,\ell)\le\bar p$ show that the exponent is at most
$-cnM^*_{ak\ell}x^2/\bar p$. Thus, writing
\[
 U_{ak\ell}:=M^*_{ak\ell}
 \frac{[\overline\bB_a^{(-i)}(k,\ell)-\bB_a(k,\ell)]^2}{\bar p},
 \quad
 F_{ak\ell}:=\left\{
 |\overline\bB_a^{(-i)}(k,\ell)-\bB_a(k,\ell)|\le c_0\bB_a(k,\ell)/4
 \right\},
\]
we have
\[
 \mathbb P(U_{ak\ell}>u,F_{ak\ell})\le2e^{-cnu},\qquad u\ge0.
\]
For $0\le t<cn$, it follows that
\begin{align*}
 \mathbb E[e^{tU_{ak\ell}}\mathbf1(F_{ak\ell})]
 &=\mathbb P(F_{ak\ell})+
   \int_0^\infty te^{tu}\mathbb P(U_{ak\ell}>u,F_{ak\ell})\,du\\
 &\le1+\frac{2t}{cn-t}.
\end{align*}
The oracle averages are independent because their observation sets
are disjoint. Thus, for $0\le t\le cn/2$,
\[
 \mathbb E[e^{tR_i}\mathbf1(\mathcal F_i)]
 \le\left(1+\frac{2t}{cn-t}\right)^{K(K+1)}
 \le\exp(CK^2t/n),
\]
which proves \eqref{fs:oracle-quadratic-mgf} after reducing $c$.

At $x=c_0\bB_a(k,\ell)/4$, the same Bernstein bound and
$\bB_a(k,\ell)\ge\bar p/G$ give
\[
 \mathbb P(\mathcal F_i^c)
 \le2K(K+1)\exp(-cn d^2\bar p/K^2).
\]
By $K^2\log n=o(n)$, $d^2\bar p\to\infty$ and
$I^*\lesssim d^2\bar p$, we have
\[
 \frac{nd^2\bar p}{K^2\log n}\to\infty,
 \qquad\frac{nd^2\bar p}{K^2I^*}\to\infty.
\]
A union bound over $i\in[n]$ proves the probability bound for every
fixed $A>0$.
\end{proof}

For the target layer $\bx=\bX_i=(x_\omega)_{\omega\in\mathcal J_d}$, write its log likelihood at a block matrix and a
community vector as
\[
 \ell(\bx;\bH,\boldsymbol\tau)
 :=\sum_{\omega=(u,v)\in\mathcal J_d}
 \log f_{\theta(\bH[\boldsymbol\tau(u),\boldsymbol\tau(v)])}(x_\omega).
\]
Terms independent of $(\bH,\boldsymbol\tau)$ may be omitted.
For $m=z_i^*$ and a fixed true cell $A^*_{mk\ell}$, its log likelihood
as a function of the mean is
\[
 p\longmapsto
 \begin{cases}
 \displaystyle\sum_{\omega\in A^*_{mk\ell}}
 [x_\omega\log p+(1-x_\omega)\log(1-p)],&\text{Bernoulli},\\[2mm]
 \displaystyle\sum_{\omega\in A^*_{mk\ell}}(x_\omega\log p-p),&\text{Poisson}.
 \end{cases}
\]
This is concave and maximized at
$(M^*_{mk\ell})^{-1}\sum_{\omega\in A^*_{mk\ell}}x_\omega$. Since
\[
 \overline\bB_m(k,\ell)
 =\left(1-\frac1{n_m}\right)\overline\bB_m^{(-i)}(k,\ell)
       +\frac1{n_m}\frac{\sum_{\omega\in A^*_{mk\ell}}x_\omega}
                              {M^*_{mk\ell}},
\]
we obtain
\begin{equation}\label{fs:oracle-own-class-monotonicity}
 \ell(\bx;\overline\bB_m,\bsig_m)
 \ge\ell(\bx;\overline\bB_m^{(-i)},\bsig_m).
\end{equation}
Since $\overline\bB_{3-m}=\overline\bB_{3-m}^{(-i)}$, we have, for every $\bx$,
\begin{align*}
 &\ell(\bx;\overline\bB_{3-m},\bsig_{3-m})
             -\ell(\bx;\overline\bB_m,\bsig_m)\\
 &\qquad\le\ell(\bx;\overline\bB_{3-m}^{(-i)},\bsig_{3-m})
             -\ell(\bx;\overline\bB_m^{(-i)},\bsig_m).
\end{align*}

\begin{theorem}\label{fs:main}\label{fs:oracle-refinement}
Suppose $K^2\log n=o(n)$, $I^*\to\infty$, and
\[
 d^2\bar p\eta_\sigma=o(I^*),\qquad
 d\eta_\sigma\log(eK/\eta_\sigma)=o(I^*).
\]
For a fixed $A>0$, suppose that, outside an event of probability $C n^{-A}$, the
full-sample oracle means lie between $\bar p/(2G)$ and $2\bar p$,
and that there is a deterministic $\delta_n=o(I^*/d^2)$ such that
\begin{equation}\label{fs:candidate-oracle-entry}
 \max_{m,k\le\ell}
 |\widehat\bB_m(k,\ell)-\overline\bB_m(k,\ell)|\le\delta_n
\end{equation}
simultaneously for all aligned candidate assignments in
\eqref{fs:candidate-assignments}. In the Bernoulli case suppose
$\bar p\le1-c_0$, and suppose that the full-sample oracle means are at most $1-3c_0/4$ on the same event. If the jointly aligned initialization event in
\eqref{eq:init-cond} has probability at least $1-Cn^{-2}$, then,
for a deterministic $o(1)$,
\begin{align}
 \mathbb E h(\widehat\bz,\bz^*)
 &\le \exp(-[1-o(1)]I^*)
       +\tfrac12(Cn^{-2}+C n^{-A}),\label{fs:risk}\\
 \mathbb P\!\left\{h(\widehat\bz,\bz^*)>
 e^{(I^*)^{1-\epsilon}-[1-o(1)]I^*}\right\}
 &\le e^{-(I^*)^{1-\epsilon}}+Cn^{-2}+C n^{-A},
 \qquad 0<\epsilon<1.\label{fs:tail}
\end{align}
Exact recovery has failure probability at most
$n\exp(-[1-o(1)]I^*)+Cn^{-2}+C n^{-A}$.
\end{theorem}

\begin{proof}
Let $\mathcal G$ be the intersection of the initialization event,
the event in \eqref{fs:candidate-oracle-entry}, the stated event for
the full-sample oracle means, and $\bigcap_i\mathcal F_i$.
Lemma~\ref{fs:oracle-block-mgf} gives
$\mathbb P(\mathcal G^c)\le Cn^{-2}+C n^{-A}$ after enlarging $C$.
Since $I^*\lesssim d^2\bar p$, $\delta_n/\bar p=o(1)$.
All candidate and oracle means on $\mathcal G$ are therefore
between $\bar p/(3G)$ and $3\bar p$ for all sufficiently large $n$.
For Bernoulli observations they are also bounded above by $1-c_0/2$.
Without loss of generality, all labels below are aligned with the true labels.

Fix an aligned candidate pair $\boldsymbol\tau=(\boldsymbol\tau_1,
\boldsymbol\tau_2)$, and define
\[
 E_{\boldsymbol\tau}:=
 \bigcup_{a=1}^2\{(u,v)\in\mathcal J_d:
      \boldsymbol\tau_a(u)\ne\bsig_a(u)\text{ or }\boldsymbol\tau_a(v)\ne\bsig_a(v)\}.
\]
Since at most $d$ unordered edges meet any vertex,
\[
 |E_{\boldsymbol\tau}|
 \le d\sum_{a=1}^2\sum_{u=1}^d
      \mathbf1[\boldsymbol\tau_a(u)\ne\bsig_a(u)]
 \le2d^2\eta_\sigma.
\]
If $p,q\in[\bar p/(3G),3\bar p]$ (and $p,q\le1-c_0/2$ in the Bernoulli case), the canonical links satisfy
$|\theta(p)-\theta(q)|\le C|p-q|/\bar p$ and
$|\theta(p)-\theta(q)|\le C$. The Bernoulli normalizing term
$\log(1-p)$ has derivative in absolute value at most $2/c_0$ on this interval;
the Poisson normalizing term is $-p$.
Applying the first bound to
$\widehat\bB_a[\boldsymbol\tau_a(u),\boldsymbol\tau_a(v)]-
\overline\bB_a[\boldsymbol\tau_a(u),\boldsymbol\tau_a(v)]$, and the second to the
change from $\boldsymbol\tau_a$ to $\bsig_a$, gives
\begin{align}\label{fs:positive-count-envelope}
 &\sum_{a=1}^2\left|
 \ell(\bx;\widehat\bB_a,\boldsymbol\tau_a)
       -\ell(\bx;\overline\bB_a,\bsig_a)\right|\nonumber\\
 &\quad\le C\frac{\delta_n}{\bar p}\sum_\omega x_\omega
       +C\sum_{\omega\in E_{\boldsymbol\tau}}x_\omega
       +Cd^2\delta_n+Cd^2\bar p\eta_\sigma.
\end{align}
Now fix $i\in[n]$ and let $m=z_i^*$. By
\eqref{fs:oracle-own-class-monotonicity}, $\{\widehat z_i\ne z_i^*,\mathcal G\}$
for this candidate pair is contained in
\begin{equation}\label{fs:oracle-error-inclusion}
\begin{split}
 \bigg\{&\ell(\bX_i;\overline\bB_{3-m}^{(-i)},\bsig_{3-m})
       -\ell(\bX_i;\overline\bB_m^{(-i)},\bsig_m)\\
 &+C\frac{\delta_n}{\bar p}\sum_\omega \bX_i(\omega)
  +C\sum_{\omega\in E_{\boldsymbol\tau}}\bX_i(\omega)
  +Cd^2\delta_n+Cd^2\bar p\eta_\sigma\ge0,\quad\mathcal F_i
 \bigg\}.
\end{split}
\end{equation}
Condition on all layers other than $i$.
For $a\in[2]$ and $\omega\in\mathcal J_d$, put
\[
 e_{a\omega}:=
 \theta(\overline\bB_a^{(-i)}[\bsig_a(u),\bsig_a(v)])
             -\theta(p_{a\omega}),\qquad \omega=(u,v),
\]
and put
\[
 w_\omega:=C\delta_n/\bar p+C\mathbf1(\omega\in E_{\boldsymbol\tau}).
\]
On $\mathcal F_i$, $|e_{a\omega}|\le C$ and
$\sum_{a,\omega}\bar p e_{a\omega}^2\le CR_i$.
Also $|w_\omega|\le C$ and
\begin{equation}\label{fs:positive-envelope-size}
 \sum_\omega\bar p w_\omega+
 \sum_\omega\bar p w_\omega^2
 \le Cd^2\delta_n+Cd^2\bar p\eta_\sigma.
\end{equation}
Here the bound for the squared terms follows from $\delta_n/\bar p=o(1)$.

Suppose $m=1$, and choose
\[
 s_*\in\underset{s\in[0,1]}{\arg\min}\,
 \sum_\omega\left\{
 b((1-s)\theta(p_{1\omega})+s\theta(p_{2\omega}))
 -(1-s)b(\theta(p_{1\omega}))-s b(\theta(p_{2\omega}))
 \right\}.
\]
The minimum is $-I^*$. Set
$\theta_{*\omega}=(1-s_*)\theta(p_{1\omega})+
s_*\theta(p_{2\omega})$. Their means
$q_\omega=b'(\theta_{*\omega})$ lie between $p_{1\omega}$ and
$p_{2\omega}$. Conditional on the other layers, the logarithm of
the moment generating function of the random terms inside
\eqref{fs:oracle-error-inclusion}, at $s_*$, equals
\begin{align*}
 -I^*&+\sum_\omega
 \left\{b(\theta_{*\omega}+s_*(e_{2\omega}-e_{1\omega}+w_\omega))
                       -b(\theta_{*\omega})\right\}\\
 &-s_*\sum_\omega\left\{
 b(\theta(p_{2\omega})+e_{2\omega})-b(\theta(p_{2\omega}))
 -b(\theta(p_{1\omega})+e_{1\omega})+b(\theta(p_{1\omega}))
 \right\}.
\end{align*}
On the line segments in this display, $b''\le C\bar p$ because
the canonical perturbations are bounded. By Taylor's formula, the
preceding display is at most
\begin{align*}
 -I^*&+s_*\sum_\omega[
 (q_\omega-p_{2\omega})e_{2\omega}
 -(q_\omega-p_{1\omega})e_{1\omega}+q_\omega w_\omega]\\
 &+C\sum_{a,\omega}\bar p e_{a\omega}^2
   +C\sum_\omega\bar p w_\omega^2\\
 &\le-I^*+C\sqrt{I^*R_i}+CR_i
             +Cd^2\delta_n+Cd^2\bar p\eta_\sigma.
\end{align*}
The last step uses \eqref{fs:positive-envelope-size} and
\[
 \sum_{a,\omega}(q_\omega-p_{a\omega})^2/\bar p\le CI^*.
\]
The same bound holds for $m=2$ with $1-s_*$ in place of $s_*$.
By Markov's inequality and \eqref{fs:oracle-error-inclusion}, we obtain
\begin{align}\label{fs:oracle-integrated-bound}
 &\mathbb P(\text{event in \eqref{fs:oracle-error-inclusion}})
 \nonumber\\
 &\quad\le
 e^{-I^*+Cd^2\delta_n+Cd^2\bar p\eta_\sigma}
 \mathbb E\left[
 e^{C\sqrt{I^*R_i}+CR_i}\mathbf1(\mathcal F_i)\right].
\end{align}

By Young's inequality,
\[
 C\sqrt{I^*R_i}
 \le\frac{I^*+K\sqrt{I^*}}{\sqrt n}
       +\frac{C\sqrt n}{1+K/\sqrt{I^*}}R_i.
\]
Applying \eqref{fs:oracle-quadratic-mgf} with
$t=C[1+\sqrt n/(1+K/\sqrt{I^*})]\le cn$ gives
\begin{align*}
 \mathbb E\left[e^{C\sqrt{I^*R_i}+CR_i}\mathbf1(\mathcal F_i)\right]
 &\le\exp\left(\frac{I^*+K\sqrt{I^*}}{\sqrt n}
       +\frac{CK^2}{n}
          \left(1+\frac{\sqrt n}{1+K/\sqrt{I^*}}\right)\right)\\
 &\le\exp\left(C\left(\frac{I^*+K\sqrt{I^*}}{\sqrt n}
                         +\frac{K^2}{n}\right)\right)
 =\exp(o(I^*)),
\end{align*}
since $K^2\log n=o(n)$ and $I^*\to\infty$.
Combining this bound with $d^2\delta_n=o(I^*)$ and
$d^2\bar p\eta_\sigma=o(I^*)$, we bound the probability
in \eqref{fs:oracle-integrated-bound} by
$\exp(-[1-o(1)]I^*)$, uniformly over the candidate pair.

Taking a union bound over at most
$\exp(C d\eta_\sigma\log(eK/\eta_\sigma))$ candidate node pairs gives
\[
 \mathbb P(\widehat z_i\ne z_i^*,\mathcal G)
 \le\exp(-[1-o(1)]I^*).
\]
{\color{black}Averaging over $i$ gives}
\[
 {\color{black}\mathbb E[h(\widehat\bz,\bz^*)\mathbf1(\mathcal G)]
 \le\frac1n\sum_{i=1}^n\mathbb P(\widehat z_i\ne z_i^*,\mathcal G)}
 \le\exp(-[1-o(1)]I^*).
\]
Since $h\le1/2$,
\begin{align*}
 {\color{black}\mathbb E h(\widehat\bz,\bz^*)}
 &\le\exp(-[1-o(1)]I^*)+\tfrac12\mathbb P(\mathcal G^c)\\
 &\le\exp(-[1-o(1)]I^*)+\tfrac12(Cn^{-2}+Cn^{-A}).
\end{align*}
Markov's inequality and a union bound yield
\begin{align*}
 \mathbb P\!\left\{h(\widehat\bz,\bz^*)>
   e^{(I^*)^{1-\epsilon}-[1-o(1)]I^*}\right\}
 &\le e^{-(I^*)^{1-\epsilon}}+Cn^{-2}+Cn^{-A},\\
 \mathbb P[h(\widehat\bz,\bz^*)>0]
 &\le n\exp(-[1-o(1)]I^*)+Cn^{-2}+Cn^{-A}.
\end{align*}
These are \eqref{fs:risk}, \eqref{fs:tail}, and the exact-recovery bound.
Under the assumptions of Lemma~\ref{fs:oracle-comparison}, let
$\mathcal H$ be the intersection of the block-average comparison
event, the relative event for the full-sample oracle means, and
$\bigcap_i\mathcal F_i$. Lemmas~\ref{fs:oracle-comparison} and~\ref{fs:oracle-block-mgf}
give
\[
 \mathbb P(\mathcal H^c)\le C n^{-A}e^{-4I^*}.
\]
In the Bernoulli case, the relative oracle bound is $o(1)$ and hence
ensures $\overline\bB_a(k,\ell)\le1-3c_0/4$ for all sufficiently
large $n$. The preceding likelihood argument on
$\mathcal A_0\cap\mathcal H$ therefore gives
\begin{align*}
 \mathbb E[h(\widehat\bz,\bz^*)\mathbf1(\mathcal A_0)]
 &\le \mathbb E[h(\widehat\bz,\bz^*)\mathbf1(\mathcal A_0\cap\mathcal H)]
       +\tfrac12\mathbb P(\mathcal H^c)\\
 &\le \exp(-[1-o(1)]I^*).
\end{align*}
For any fixed $\epsilon\in(0,1)$, define
\[
 \mathcal A_1:=
 \left\{h(\widehat\bz,\bz^*)\mathbf1(\mathcal A_0)
       \le \exp((I^*)^{1-\epsilon}-[1-o(1)]I^*)\right\}.
\]
Markov's inequality yields
\[
 \mathbb P(\mathcal A_1^c)\le\exp(-(I^*)^{1-\epsilon}).
\]
Since $(I^*)^{1-\epsilon}=o(I^*)$, the claimed bound holds on
$\mathcal A_0\cap\mathcal A_1$. The Poisson argument gives
$\mathcal A_1^{\rm P}$ in the same way. Adding
$\mathbb P(\mathcal A_0^c)/2$ or
$\mathbb P[(\mathcal A_0^{\rm P})^c]/2$ gives the respective
unconditional expectation bound.
\end{proof}

Theorems~\ref{thm:main} and~\ref{thm:MMPBM-main} follow from
the preceding proof. Their second accuracy conditions give
$d^2\bar p\eta_\sigma=o(I^*)$, and the entropy conditions were verified
above. The initialization failure probabilities give the respective
$O(n^{-2})$ and $O((n\vee d)^{-2})$ terms.

}
\ZL{
\subsection{Refinement with exact community estimates}
\label{fs:exact-community-section}

We consider Bernoulli observations and use the notation of
Section~\ref{fs:section}. In this subsection, take
$\boldsymbol\tau_m=\bsig_m$, $m\in[2]$, in the definition of $A_j$.
Thus, for $j=(m,k,\ell)$,
\[
 A_j=\{(u,v)\in\mathcal J_d:
 (\bsig_m(u),\bsig_m(v))\in\{(k,\ell),(\ell,k)\}\},
 \qquad M_j=|A_j|,\qquad B_j=\bB_m(k,\ell).
\]
Set
\[
 \mu_{\min}=\min_j M_jB_j.
\]
For fixed $A>0$, define
\begin{equation}\label{fs:exact-refinement-radius}
 \begin{split}
 q&=\eta_z\log\frac e{\eta_z}
       +\frac{A\log n+4I^*+\log(2K(K+1))}{n},\\
 R&=\eta_z^2 I^*+K(K+1)q\left(1+\frac q{\mu_{\min}}\right),
 \end{split}
\end{equation}
\begin{proposition}\label{fs:exact-community-refinement}
Suppose the observations are Bernoulli and
\[
 \begin{gathered}
 |\{i:z_i^*=m\}|\ge\frac n{2\alpha},\qquad
 0<p_{m\omega}\le1-c_0,\\
 G^{-1}\le\frac{p_{m\omega}}{p_{m'\omega'}}\le G,
 \qquad m,m'\in[2],\quad\omega,\omega'\in\mathcal J_d,
 \end{gathered}
\]
for fixed $\alpha,G\ge1$ and $c_0\in(0,1)$, and suppose $M_j>0$ for every cell.
Let $\calA_0$ be the event in \eqref{eq:init-cond} with
$\eta_\sigma=0$, and suppose $\PP(\calA_0^c)\le\delta_{\rm init}$.
There exist constants $c,C>0$, depending only on $\alpha,G,c_0$,
such that if
\[
 \eta_z\le c,\qquad R\le c\mu_{\min},
\]
then the output of Algorithm~\ref{algorithm} satisfies
\begin{equation}\label{fs:exact-refinement-expectation}
 \begin{split}
 \EE h(\widehat\bz,\bz^*)
 \le{}&2\exp\left(-I^*+C\sqrt{I^*R}+CR
   +C\sqrt R\left(\sqrt{K(K+1)}
           +\frac{K(K+1)}{\sqrt{\mu_{\min}}}\right)\right)\\
 &+\delta_{\rm init}+n^{-A}\exp(-4I^*).
 \end{split}
\end{equation}
For every $u>0$,
\begin{equation}\label{fs:exact-refinement-probability}
 \begin{split}
 \PP\Bigg(&h(\widehat\bz,\bz^*)>
 2\exp\Bigg(-I^*+C\sqrt{I^*R}+CR\\
 &\qquad+C\sqrt R\left(\sqrt{K(K+1)}
            +\frac{K(K+1)}{\sqrt{\mu_{\min}}}\right)+u\Bigg)\Bigg)\\
 &\le \exp(-u)+\delta_{\rm init}+n^{-A}\exp(-4I^*).
 \end{split}
\end{equation}
The initial estimates and the block averages may use the full sample.
\end{proposition}

\begin{proof}
Without loss of generality, assume all labels are aligned as in
Section~\ref{fs:section}. On $\calA_0$, the estimated community vectors
equal $\bsig_1,\bsig_2$, and $\widetilde\bz$ belongs to
\[
 \mathcal Z_{\eta_z}
 =\left\{\boldsymbol\zeta\in[2]^n:
       \sum_{i=1}^n\mathbb I(\zeta_i\ne z_i^*)\le n\eta_z\right\}.
\]
For a fixed $\boldsymbol\zeta\in\mathcal Z_{\eta_z}$, write
$N_m(\boldsymbol\zeta)=|\{i:\zeta_i=m\}|$. By layer balance,
\[
 N_m(\boldsymbol\zeta)
 \ge\frac n{2\alpha}-n\eta_z\ge\frac n{4\alpha}
\]
when $c\le1/(4\alpha)$. For $j=(m,k,\ell)$, the block average is
\begin{equation}\label{fs:exact-candidate-average}
 \widehat B_j(\boldsymbol\zeta)
 =\frac{\sum_{i:\zeta_i=m}\sum_{\omega\in A_j}\bX_i(\omega)}
       {N_m(\boldsymbol\zeta)M_j}.
\end{equation}

By \eqref{fs:exact-candidate-average}, for any $j=(m,k,\ell)$, we have
\begin{equation}\label{fs:exact-mixture-bias}
 \EE\widehat B_j(\boldsymbol\zeta)-B_j
 =\frac{|\{i:\zeta_i=m,\ z_i^*\ne m\}|}{N_m(\boldsymbol\zeta)}
       \frac1{M_j}\sum_{\omega\in A_j}(p_{3-m,\omega}-p_{m\omega}).
\end{equation}
In \eqref{fs:exact-mixture-bias},
\[
 \frac{|\{i:\zeta_i=m,\ z_i^*\ne m\}|}{N_m(\boldsymbol\zeta)}
 \le\frac{n\eta_z}{n/(4\alpha)}=4\alpha\eta_z.
\]
Since $p_{m\omega}=B_j$ on $A_j$, Jensen's inequality and
\eqref{fs:information} imply
\begin{align}
 &\sum_{j=(m,k,\ell)}\frac{M_j}{B_j}
 \left[\frac1{M_j}\sum_{\omega\in A_j}
               (p_{3-m,\omega}-p_{m\omega})\right]^2\notag\\
 &\qquad\le\sum_{\omega\in\mathcal J_d}
       \frac{(p_{3-m,\omega}-p_{m\omega})^2}{p_{m\omega}}
 \le CI^*.
 \label{fs:exact-projected-contrast}
\end{align}
Here the first sum ranges over $1\le k\le\ell\le K$ for the indicated
component $m$. Summing over $m\in[2]$ in
\eqref{fs:exact-projected-contrast} yields
\begin{equation}\label{fs:exact-bias-radius}
 \sum_j\frac{M_j[\EE\widehat B_j(\boldsymbol\zeta)-B_j]^2}{B_j}
 \le C\eta_z^2 I^*.
\end{equation}
Also, by homogeneity,
\[
 \max_j\frac{|\EE\widehat B_j(\boldsymbol\zeta)-B_j|}{B_j}
 \le C\eta_z.
\]

For a fixed candidate, the numerator of
\eqref{fs:exact-candidate-average} is a sum of independent Bernoulli
variables, with variance at most $CN_m(\boldsymbol\zeta)M_jB_j$.
By Bernstein's inequality, for each $t>0$,
\[
 \left|\widehat B_j(\boldsymbol\zeta)
                 -\EE\widehat B_j(\boldsymbol\zeta)\right|
 \le C\left(\sqrt{\frac{B_jt}{nM_j}}+\frac t{nM_j}\right)
\]
with probability at least $1-2\exp(-t)$. Moreover, for $\eta_z\le1/2$,
\[
 |\mathcal Z_{\eta_z}|=\sum_{a=0}^{\lfloor n\eta_z\rfloor}\binom na
 \le\exp\left(n\eta_z\log\frac e{\eta_z}\right).
\]
Taking $t=nq$ and a union bound over these candidates and the $K(K+1)$
cells gives an event $\calE$ such that
\begin{equation}\label{fs:exact-concentration-probability}
 \PP(\calE^c)\le n^{-A}\exp(-4I^*)
\end{equation}
and, on $\calE$, simultaneously for all candidates and cells,
\begin{equation}\label{fs:exact-centered-average}
 \left|\widehat B_j(\boldsymbol\zeta)
                -\EE\widehat B_j(\boldsymbol\zeta)\right|
 \le C\left(\sqrt{\frac{B_jq}{M_j}}+\frac q{M_j}\right).
\end{equation}
Consequently,
\begin{align*}
 \sum_j\frac{M_j[\widehat B_j(\boldsymbol\zeta)
                     -\EE\widehat B_j(\boldsymbol\zeta)]^2}{B_j}
 &\le C\sum_j\left(q+\frac{q^2}{M_jB_j}\right)\\
 &\le CK(K+1)q\left(1+\frac q{\mu_{\min}}\right).
\end{align*}
Combining this inequality with \eqref{fs:exact-bias-radius} gives
\begin{equation}\label{fs:exact-estimation-radius}
 \sum_j\frac{M_j[\widehat B_j(\boldsymbol\zeta)-B_j]^2}{B_j}
 \le CR.
\end{equation}
Dividing \eqref{fs:exact-centered-average} by $B_j$ also gives
\begin{equation}\label{fs:exact-relative-average}
 \max_j\frac{|\widehat B_j(\boldsymbol\zeta)-B_j|}{B_j}
 \le C\left(\eta_z+\sqrt{\frac q{\mu_{\min}}}
                       +\frac q{\mu_{\min}}\right).
\end{equation}
Both bounds hold uniformly over $\mathcal Z_{\eta_z}$, hence for
$\widetilde\bz$ on $\calA_0\cap\calE$.

Fix fitted entries $r_j$ satisfying \eqref{fs:exact-estimation-radius}
and \eqref{fs:exact-relative-average}. For sufficiently small $c$,
\[
 \frac q{\mu_{\min}}
 \le\frac R{K(K+1)\mu_{\min}}\le\frac c{K(K+1)},
 \qquad
 \max_j\frac{|r_j-B_j|}{B_j}\le\min\{1/2,c_0/2\}.
\]
Consequently,
\[
 B_j/2\le r_j\le3B_j/2,
 \qquad 0<r_j\le1-c_0/2.
\]
Put $u_j=\theta(r_j)-\theta(B_j)$, with $\theta,b$ as in
Section~\ref{fs:section}.
Since $\theta'(p)=1/[p(1-p)]$, the mean value theorem gives
\begin{equation}\label{fs:exact-canonical-radius}
 \sum_jM_jB_ju_j^2
 \le C\sum_j\frac{M_j(r_j-B_j)^2}{B_j}\le CR.
\end{equation}
Along the segment between $\theta(B_j)$ and $\theta(r_j)$,
$b''(\theta)\le CB_j$. Taylor's expansion therefore gives
\begin{equation}\label{fs:exact-canonical-remainder}
 \left|b(\theta(B_j)+u_j)-b(\theta(B_j))-B_ju_j\right|
 \le CB_ju_j^2.
\end{equation}

Choose $s_*\in[0,1]$ attaining the maximum in the definition of $I^*$,
and define
\[
 q_{*\omega}
 =\frac{p_{1\omega}^{1-s_*}p_{2\omega}^{s_*}}
 {p_{1\omega}^{1-s_*}p_{2\omega}^{s_*}
 +(1-p_{1\omega})^{1-s_*}(1-p_{2\omega})^{s_*}}.
\]
Let $Q_*$ be the product Bernoulli law with probabilities $q_{*\omega}$.
Since
\[
 \theta(q_{*\omega})=(1-s_*)\theta(p_{1\omega})+s_*\theta(p_{2\omega}),
\]
$q_{*\omega}$ lies between $p_{1\omega}$ and $p_{2\omega}$. By
\eqref{fs:information},
\begin{equation}\label{fs:exact-tilted-contrast}
 \sum_{m=1}^2\sum_{\omega\in\mathcal J_d}
 \frac{(q_{*\omega}-p_{m\omega})^2}{p_{m\omega}}\le CI^*.
\end{equation}
For $j=(m,k,\ell)$ and $\omega\in A_j$, write
$r_{m\omega}=r_j$ and $u_{m\omega}=u_j$. For a layer $\bX$, define
\begin{align*}
 L_p(\bX)&=\sum_{\omega\in\mathcal J_d}
 \left(\bX(\omega)[\theta(p_{2\omega})-\theta(p_{1\omega})]
           -b(\theta(p_{2\omega}))+b(\theta(p_{1\omega}))\right),\\
 L_r(\bX)&=\sum_{\omega\in\mathcal J_d}
 \left(\bX(\omega)[\theta(r_{2\omega})-\theta(r_{1\omega})]
           -b(\theta(r_{2\omega}))+b(\theta(r_{1\omega}))\right).
\end{align*}
Define $\mathbf g(\bX)\in\mathbb R^{K(K+1)}$ by
\[
 g_j(\bX)=\frac{\sum_{\omega\in A_j}[\bX(\omega)-q_{*\omega}]}
                    {\sqrt{M_jB_j}}.
\]
Subtracting the two likelihoods and applying
\eqref{fs:exact-canonical-remainder}, we obtain
\begin{align*}
 &\left|L_r(\bX)-L_p(\bX)
      -\sum_\omega[\bX(\omega)-q_{*\omega}](u_{2\omega}-u_{1\omega})
      -\sum_{m=1}^2(-1)^m\sum_\omega
                   (q_{*\omega}-p_{m\omega})u_{m\omega}\right|\\
 &\qquad\le C\sum_{m=1}^2\sum_\omega p_{m\omega}u_{m\omega}^2
 =C\sum_jM_jB_ju_j^2\le CR.
\end{align*}
By Cauchy--Schwarz and \eqref{fs:exact-canonical-radius},
\begin{align*}
 \left|\sum_\omega[\bX(\omega)-q_{*\omega}]
                         (u_{2\omega}-u_{1\omega})\right|
 &\le\left(\sum_jM_jB_ju_j^2\right)^{1/2}\|\mathbf g(\bX)\|_2\\
 &\le C\sqrt R\,\|\mathbf g(\bX)\|_2.
\end{align*}
Similarly, by \eqref{fs:exact-tilted-contrast},
\[
 \left|\sum_{m=1}^2(-1)^m\sum_\omega
             (q_{*\omega}-p_{m\omega})u_{m\omega}\right|
 \le\left[\sum_{m,\omega}\frac{(q_{*\omega}-p_{m\omega})^2}{p_{m\omega}}\right]^{1/2}
       \left(\sum_jM_jB_ju_j^2\right)^{1/2}
 \le C\sqrt{I^*R}.
\]
Uniformly over fitted entries satisfying
\eqref{fs:exact-estimation-radius} and \eqref{fs:exact-relative-average},
\begin{equation}\label{fs:exact-likelihood-comparison}
 |L_r(\bX)-L_p(\bX)|
 \le C\sqrt{I^*R}+CR+C\sqrt R\,\|\mathbf g(\bX)\|_2.
\end{equation}
For any unit vector $\mathbf v\in\mathbb R^{K(K+1)}$,
\[
 \mathbf v^\top\mathbf g(\bX)
 =\sum_\omega\left(\sum_j\frac{v_j\mathbf1(\omega\in A_j)}{\sqrt{M_jB_j}}\right)
                   [\bX(\omega)-q_{*\omega}].
\]
Since each edge belongs to one cell in each component,
\[
 \max_\omega\left|\sum_j\frac{v_j\mathbf1(\omega\in A_j)}{\sqrt{M_jB_j}}\right|
 \le\sqrt{\frac2{\mu_{\min}}},
\]
and the variance under $Q_*$ is at most
\begin{align*}
 \sum_\omega q_{*\omega}
   \left(\sum_j\frac{v_j\mathbf1(\omega\in A_j)}{\sqrt{M_jB_j}}\right)^2
 &\le2\sum_j\frac{v_j^2\sum_{\omega\in A_j}q_{*\omega}}{M_jB_j}\\
 &\le2G\sum_jv_j^2=2G.
\end{align*}
Bernstein's inequality gives, for $t>0$,
\[
 Q_*\left(\mathbf v^\top\mathbf g(\bX)
                 >C\sqrt t+\frac{Ct}{\sqrt{\mu_{\min}}}\right)\le \exp(-t).
\]
Take a $1/2$-net $\mathcal N$ of the Euclidean unit sphere
in $\mathbb R^{K(K+1)}$, with $|\mathcal N|\le5^{K(K+1)}$.
Since $\|\mathbf g\|_2\le2\max_{\mathbf v\in\mathcal N}\mathbf v^\top\mathbf g$,
a union bound with $t=K(K+1)\log5+x$ yields
\begin{equation}\label{fs:exact-score-tail}
 Q_*\left(\|\mathbf g(\bX)\|_2>
 C\left[\sqrt{K(K+1)+x}+\frac{K(K+1)+x}{\sqrt{\mu_{\min}}}\right]\right)
 \le \exp(-x),\qquad x\ge0.
\end{equation}
Integrating \eqref{fs:exact-score-tail} and using
$\sqrt{K(K+1)+x}\le\sqrt{K(K+1)}+\sqrt x$, we obtain, for $a\ge0$,
\begin{align*}
 \EE_{Q_*}\exp(a\|\mathbf g(\bX)\|_2)
 &\le\exp\left(Ca\left(\sqrt{K(K+1)}
              +\frac{K(K+1)}{\sqrt{\mu_{\min}}}\right)\right)\\
 &\qquad\times\int_0^\infty
         \exp\left(Ca\sqrt x+\frac{Cax}{\sqrt{\mu_{\min}}}-x\right)\,dx.
\end{align*}
If $Ca/\sqrt{\mu_{\min}}\le1/4$, then
$Ca\sqrt x\le x/4+C^2a^2$, and hence
\[
 \int_0^\infty
    \exp\left(Ca\sqrt x+\frac{Cax}{\sqrt{\mu_{\min}}}-x\right)\,dx
 \le \exp(C^2a^2)\int_0^\infty \exp(-x/2)\,dx=2\exp(C^2a^2).
\]
Increasing $C$ if necessary, we obtain
\begin{equation}\label{fs:exact-score-moment}
 \EE_{Q_*}\exp(a\|\mathbf g(\bX)\|_2)
 \le2\exp\left(Ca\left(\sqrt{K(K+1)}
           +\frac{K(K+1)}{\sqrt{\mu_{\min}}}\right)+Ca^2\right)
\end{equation}
whenever $Ca/\sqrt{\mu_{\min}}\le1/4$.

Finally, fix $i\in[n]$ and suppose $z_i^*=1$. On $\calA_0\cap\calE$, let
$r_j$ be the aligned fitted entries in Algorithm~\ref{algorithm}.
By $\{\widehat z_i\ne z_i^*\}\subseteq\{L_r(\bX_i)\ge0\}$ and
\eqref{fs:exact-likelihood-comparison}, we have
\begin{align*}
 &\mathbb I(\widehat z_i\ne z_i^*)\mathbb I(\calA_0\cap\calE)\\
 &\quad\le\exp\left(s_*L_p(\bX_i)+C\sqrt{I^*R}+CR
                      +C\sqrt R\,\|\mathbf g(\bX_i)\|_2\right).
\end{align*}
Let $f_m$ be the probability mass function of a layer from component
$m$. Since
\[
 \exp(s_*L_p(\bX))f_1(\bX)=\exp(-I^*)Q_*(\bX),
\]
taking expectations and applying \eqref{fs:exact-score-moment}
with $a=C\sqrt R$ gives
\begin{equation}\label{fs:exact-individual-error}
 \begin{split}
 &\PP(\widehat z_i\ne z_i^*,\calA_0\cap\calE)\\
 &\quad\le2\exp\left(-I^*+C\sqrt{I^*R}+CR
       +C\sqrt R\left(\sqrt{K(K+1)}
             +\frac{K(K+1)}{\sqrt{\mu_{\min}}}\right)\right),
 \end{split}
\end{equation}
where $R\le c\mu_{\min}$ ensures $Ca/\sqrt{\mu_{\min}}\le1/4$
for sufficiently small $c$. For $z_i^*=2$, apply the same argument to
$-L_p$, $-L_r$ and $1-s_*$. The identity
\[
 \exp(-(1-s_*)L_p(\bX))f_2(\bX)=\exp(-I^*)Q_*(\bX)
\]
gives \eqref{fs:exact-individual-error} with the same bound.
By definition of the Hamming loss,
\[
 h(\widehat\bz,\bz^*)\mathbb I(\calA_0\cap\calE)
 \le\frac1n\sum_{i=1}^n
          \mathbb I(\widehat z_i\ne z_i^*,\calA_0\cap\calE).
\]
By \eqref{fs:exact-individual-error} and $0\le h\le1$,
\begin{align*}
 \EE h(\widehat\bz,\bz^*)
 &\le\frac1n\sum_{i=1}^n
        \PP(\widehat z_i\ne z_i^*,\calA_0\cap\calE)
       +\PP(\calA_0^c)+\PP(\calE^c),\\
 \PP(\calA_0^c)+\PP(\calE^c)
 &\le\delta_{\rm init}+n^{-A}\exp(-4I^*).
\end{align*}
This proves \eqref{fs:exact-refinement-expectation}. For $u>0$, applying
Markov's inequality to
$h(\widehat\bz,\bz^*)\mathbb I(\calA_0\cap\calE)$ and adding the same
failure probabilities proves \eqref{fs:exact-refinement-probability}.
\end{proof}
}

\section{Initialization and concentration arguments}\label{sec:lemma-proofs}
\ZL{
{\color{black}\subsection{Proof of Theorem~\ref{thm:init}}
Throughout, matrix norms without a subscript are operator norms, and
$\|\bA\|_{2,\infty}=\|\bA\|_{2\to\infty}=\max_j\|\bA_{j,:}\|_2$.
Put $\bar a=nd\bar p$ and suppose $n=O(d)$.}
Since $\bar p=(n\wedge d)(p_1\vee p_2)/n$ and $n=O(d)$,
the parameter space and \eqref{eq:pmqm} imply
\[
 c\bar p\le\bP_m(u,v)\le C\bar p\qquad(m\in[2],\ u,v\in[d]).
\]
{\color{black}Throughout this proof, condition on an independent split satisfying
Theorem~\ref{thm:init}'s assumptions. Constants may depend on the fixed
homogeneity and balance constants, the fixed incoherence bound,
the fixed $K$-means approximation factor, and the bound on $n/d$,
but not on $\kappa,K,n,d$. The clipping constant $C_0$ is sufficiently large.} {\color{black}For two label vectors
$\bx,\by$ of equal length, write
\[
 h(\bx,\by)=\min_{\pi}\frac{1}{|\bx|}\sum_j
 \boldsymbol{1}[x_j\ne\pi(y_j)],
\]
where $\pi$ ranges over label permutations.
Write $\mathfrak S_q$ for the permutations of $[q]$ and
$\mathbb O_q=\{\bO\in\mathbb R^{q\times q}:\bO^\top \bO=\bI_q\}$ for the
orthogonal group; minima over $\bO$ use the appropriate dimension.
Let $\be_j$ be the $j$-th coordinate vector.
The true labels and block matrices are fixed throughout.
$\operatorname{SVD}_q(\bA)$ is a matrix of
orthonormal eigenvectors of $\bA\bA^\top$ for its $q$ largest eigenvalues.
}

{\color{black}\subsubsection{The population core and incoherence}}

{\color{black}For $m\in[2]$, write
\[
 \bZ_m(j,a)=\boldsymbol 1[\bsig_m(j)=a],\qquad
 \bP_m=\bZ_m\bB_m\bZ_m^\top,\qquad
 n_m=|\{i\in[n]:z_i^*=m\}|.
\]
Thus $\bZ_m\in\{0,1\}^{d\times K}$ and
$S_a(\bsig_m)=\{j\in[d]:\bsig_m(j)=a\}$.
Let $\bZ=[\bZ_1,\bZ_2]=\bU\bD\bR^\top$ be a compact singular value decomposition,
where $\bU\in\mathbb R^{d\times r}$ and $\bR\in\mathbb R^{2K\times r}$
have orthonormal columns, and $\bD\in\mathbb R^{r\times r}$ is diagonal
with its positive singular values in decreasing order.}
Since both membership matrices have rank $K$, we have
\[
 K\le r:=\operatorname{rank}(\bZ)\le 2K-1,
 \qquad \bZ_1\boldsymbol{1}_K=\bZ_2\boldsymbol{1}_K.
\]
{\color{black}The block tensor $\bcalB\in\mathbb R^{2K\times2K\times2}$ has slices
\[
 \bcalB(:,:,1)=\begin{pmatrix}\bB_1&0\\0&0\end{pmatrix},\qquad
 \bcalB(:,:,2)=\begin{pmatrix}0&0\\0&\bB_2\end{pmatrix},
 \qquad \bcalB_0=\bcalB/\bar p.
\]}
Define
\[
 \overline{\bcalB}_0
 =\bcalB_0\times_1\bR^\top\times_2\bR^\top.
\]
By symmetry of the first two modes,
\[
 \sigma_{\min}(\overline{\bcalB}_0)
 =\min\big\{\sigma_r(\scrM_1(\overline{\bcalB}_0)),
              \sigma_2(\scrM_3(\overline{\bcalB}_0))\big\}.
\]

{\color{black}The population tensor has the orthonormal Tucker representation with
vertex factors $\bU$, normalized layer-membership factor
\[
 \bW_0(i,m)=\boldsymbol1(z_i^*=m)/\sqrt{n_m},
 \qquad i\in[n],\ m\in[2],
\]
and core
\[
 \bcalC=\bar p\,\overline{\bcalB}_0
       \times_1\bD\times_2\bD\times_3\bD_N,
 \qquad \bD_N=\operatorname{diag}(\sqrt{n_1},\sqrt{n_2}).
\]
This is the core denoted by $\bar\bcalC$ in the main text.
In particular, $\bU^\top \bU=\bI_r$, $\bW_0^\top \bW_0=\bI_2$, and
\[
 \mathbb E\bcalX
 =\bcalC\times_1\bU\times_2\bU\times_3\bW_0,
 \qquad
 (\mathbb E\bcalX)(:,:,i)=\bP_{z_i^*}.
\]
Since each vertex belongs to one community in each component,
$\|\bZ\|_F^2=2d$. Consequently,
\[
 2d=\sum_{j=1}^r\sigma_j(\bD)^2
 \le r\sigma_1(\bD)^2\le r\kappa^2\sigma_r(\bD)^2.
\]
Also $\bU=\bZ\bR\bD^{-1}$ and $\|\bZ_{j,:}\|_2=\sqrt2$, so
\[
 \|\bU\|_{2,\infty}
 \le\frac{\sqrt2}{\sigma_r(\bD)}
 \le\kappa\sqrt{r/d}.
\]
}
By orthonormality of the population factors,
\begin{equation}\label{eq:init-layer-singular-value}
 \sigma_2(\scrM_3(\bcalC))
 =\sigma_2(\scrM_3(\mathbb E\bcalX))>0.
\end{equation}
where we have used condition (A1).

{\color{black}For the split in Algorithm~\ref{algorithm:init}, write
$V_k:=\mathcal V^{[k]}$ and $N_k:=\mathcal N^{[k]}$, $k\in\{0,1\}$.
Thus $V_0\mathbin{\dot\cup}V_1=[d]$ and
$N_0\mathbin{\dot\cup}N_1=[n]$. When the splits are random, they
are independent of the observations, with each first half a uniform
subset of the prescribed size. Set
\[
 \bX_i^{(k)}=\bX_i[V_k,V_k],\qquad \bP_m^{(k)}=\bP_m[V_k,V_k],\qquad
 \bz^{*[k]}=\bz^*|_{N_k},\qquad \bsig_m^{[k]}=\bsig_m|_{V_k}.
\]
For $N\subseteq[n]$ and $V\subseteq[d]$, the tensor
$\bcalX_{N,V}$ has slices $\bX_i[V,V]$, $i\in N$.
We use column vectorization, so row $i$ of $\scrM_3(\bcalX)$
is $\operatorname{vec}(\bX_i)^\top$, with column $u+(v-1)d$
corresponding to $\bX_i(u,v)$.

For a $K$-means call with $J$ centers on the rows of
$\bA\in\mathbb R^{q\times p}$, the returned labels
$\boldsymbol{\ell}:[q]\to[J]$ and centers $\bc_1,\ldots,\bc_J\in\mathbb R^p$ satisfy
\[
 \boldsymbol{\ell}(i)\in\underset{a\in[J]}{\arg\min}\,\|\bA_{i,:}-\bc_a\|_2,
 \qquad
 \sum_{i=1}^q\|\bA_{i,:}-\bc_{\boldsymbol{\ell}(i)}\|_2^2
 \le C
 \min_{\substack{\bg:[q]\to[J]\\\ba_1,\ldots,\ba_J\in\mathbb R^p}}
       \sum_{i=1}^q\|\bA_{i,:}-\ba_{\bg(i)}\|_2^2,
\]
where the approximation factor $C\ge1$ is fixed. We use $J=2$ for layer labels
and $J=K$ for community labels. The cross-fitting step is}
\begin{align}
 \widehat \bU_k
 &=\operatorname{SVD}_r\left(\sum_{i\in N_{1-k}}\bX_i^{(k)}\right),
 \label{ir:crossfit-u}\\
 \widehat \bW_k
 &=\operatorname{SVD}_2\left(
 \scrM_3(\bcalX_{N_k,V_k})
 (\widetilde \bU_k\otimes\widetilde \bU_k)\right),
 \label{ir:crossfit-w}
\end{align}
{\color{black}Here
$\widetilde \bU_k=\mathcal P_{C_0\sqrt{r/|V_k|}}(\widehat \bU_k)$,
where $C_0$ is the algorithm's fixed constant covering the population
row bound. For $\bA\in\mathbb R^{q\times r}$ with $q\ge r$ and
$\delta>0$, the regularization operator is
\[
 \mathcal P_\delta(\bA)=\operatorname{SVD}_r
 \left(\left(\bA_{j,:}\min\{1,\delta/\|\bA_{j,:}\|_2\}\right)_{j}\right).
\]
The raw labels
$\widetilde \bz^{[k]}:N_k\to[2]$ are the returned labels of the
$2$-means call on $\widehat \bW_k$.}
{\color{black}Recover communities on $V_{1-k}$ using}
{\color{black}\[
 \sum_{i\in N_k:\widetilde z_i^{[k]}=m}\bX_i^{(1-k)}.
\]
Apply $K$-means to the rows of the matrix of its leading $K$ left
singular vectors to obtain
$\widetilde\bsig_m^{[1-k]}:V_{1-k}\to[K]$ before alignment.
Conditional on the split, $\{\bX_i^{(1-k)}:i\in N_k\}$ is independent
of $\widetilde\bz^{[k]}$, which uses only edges within $V_k$.}

{\color{black}\subsubsection{Sample-splitting properties}}

In this subsection, probabilities refer to the random node and layer
splits, with the population parameters fixed. Expectations of restricted
observations are conditional on the split.

\begin{lemma}
\label{ir:vertex-split}
Suppose $\bU^\top \bU=\bI_r$ and
$\max_j\|\bU_{j,:}\|_2^2\le\mu^2r/d$.  Let $V$ be a uniform subset of
$s\in\{\lfloor d/2\rfloor,\lceil d/2\rceil\}$ vertices.
For $0<\eta<1$, with probability at least
\[
 1-2r\exp(-c\eta^2d/(\mu^2r)),
\]
the matrix $\bH_V=\bU_V^\top \bU_V$ obeys
\[
 (1-\eta)(s/d)\bI_r\preceq \bH_V
 \preceq(1+\eta)(s/d)\bI_r.
\]
On this event, if $\bS=\bU\bT\bU^\top$ has rank $r$, then
\[
 \sigma_r(\bS[V,V])\ge(1-\eta)(s/d)\sigma_r(\bS),\qquad
 \|\bS[V,V]\|\le(1+\eta)(s/d)\|\bS\|.
\]
Moreover, $\operatorname{rank}(\bZ_V)=r$ and
\[
 \kappa(\bZ_V)\le\sqrt{(1+\eta)/(1-\eta)}\,\kappa(\bD).
\]
The projected normalized block core has exactly the same matricization
singular values as before vertex restriction.
\end{lemma}

\begin{proof}
Write $\bU_V=\bU[V,:]$ and $\bZ_V=\bZ[V,:]$. The identities
\[
 \sum_{j=1}^d\bU_{j,:}^\top\bU_{j,:}=\bI_r,\qquad
 0\preceq\bU_{j,:}^\top\bU_{j,:}\preceq\frac{\mu^2r}{d}\bI_r,
 \qquad \mathbb E\bH_V=\frac sd\bI_r
\]
and the matrix Chernoff bound for sampling without replacement give
\begin{align*}
 &\mathbb P[\lambda_{\min}(\bH_V)<(1-\eta)s/d]
 +\mathbb P[\lambda_{\max}(\bH_V)>(1+\eta)s/d]\\
 &\qquad\le2r\exp(-c\eta^2d/(\mu^2r)).
\end{align*}
Here we use $s\asymp d$ and the matrix exponential bound of
Gross and Nesme, \emph{Note on sampling without replacing from a finite
collection of matrices}, arXiv:1001.2738.
On the event in the lemma, set $\bQ_V=\bU_V\bH_V^{-1/2}$, so that
$\bQ_V^\top\bQ_V=\bI_r$. Then
\[
 \bS[V,V]=\bQ_V\bH_V^{1/2}\bT\bH_V^{1/2}\bQ_V^\top,\qquad
 \bZ_V=\bQ_V\bH_V^{1/2}\bD\bR^\top.
\]
Since $\sigma_r(\bT)=\sigma_r(\bS)$ and $\|\bT\|=\|\bS\|$,
the singular-value product inequality proves the bounds for $\bS[V,V]$.
It also gives $\operatorname{rank}(\bZ_V)=r$ and
\[
 \kappa(\bZ_V)\le
 \sqrt{\frac{\lambda_{\max}(\bH_V)}{\lambda_{\min}(\bH_V)}}\,
 \kappa(\bD)
 \le\sqrt{\frac{1+\eta}{1-\eta}}\,\kappa(\bD).
\]
The row spaces of $\bZ_V$ and $\bZ$ coincide. Thus their right
singular-vector matrices satisfy $\bR_V=\bR\bO$ for some
$\bO\in\mathbb O_r$, and
\[
 \bcalB_0\times_1\bR_V^\top\times_2\bR_V^\top
 =\overline{\bcalB}_0\times_1\bO^\top\times_2\bO^\top.
\]
Orthogonal transformations preserve the singular values of every
matricization, proving the assertion about the projected core.
\end{proof}

The same event preserves the componentwise conditions. Since
$\bZ_m\be_l\in\operatorname{col}(\bU)$,
\[
 (1-\eta)\frac sd|S_l(\bsig_m)|
 \le |S_l(\bsig_m)\cap V|
 \le(1+\eta)\frac sd|S_l(\bsig_m)|.
\]
Also, $\bP_m[V,V]=\bQ_V\bH_V^{1/2}
(\bU^\top\bP_m\bU)\bH_V^{1/2}\bQ_V^\top$ has rank $K$, and
\[
 \sigma_K(\bP_m[V,V])\ge(1-\eta)\frac sd\sigma_K(\bP_m),
 \qquad
 \sigma_1(\bP_m[V,V])\le(1+\eta)\frac sd\sigma_1(\bP_m).
\]
Thus the restricted communities satisfy the balance condition with
constant $\beta/(1-\eta)$, and the component condition-number bound
in (A2) is at most $[(1+\eta)/(1-\eta)]\kappa$.

\begin{lemma}
\label{ir:layer-split}
Let $\bS=\sum_{m=1}^2n_m\bP_m$ have rank $r$, and suppose the symmetric
$\bP_m$ are entrywise nonnegative with column spaces contained in
$\operatorname{col}(\bU)$.  If
\[
 n_m'=\tfrac12n_m(1+e_m),\qquad \max_m|e_m|\le\tau,
 \qquad \tau\kappa(\bS)\le\tfrac12,
\]
then $\bS'=\sum_mn_m'\bP_m$ has rank $r$ and
\[
 \sigma_r(\bS')\ge\tfrac14\sigma_r(\bS),\qquad
 \|\bS'\|\le\tfrac{1+\tau}{2}\|\bS\|.
\]
{\color{black}If layer classes are balanced and layers are divided uniformly into equal
halves, the count condition with
$\tau=C\log n/\sqrt n$ fails with probability at most}
$C\exp(-c\log^2n)$.  {\color{black}The aggregate condition number assumption in
the manuscript gives $\tau\kappa(\bS)\le\ZL{1/2}$ when its
constant $c$ is sufficiently small.}
\end{lemma}

\begin{proof}
Entrywise, $|\bS'-\bS/2|\le(\tau/2)\bS$, and nonnegativity gives
\[
 \|\bS'-\bS/2\|\le\frac\tau2\|\bS\|.
\]
Weyl's inequality and the triangle inequality yield
\[
 \sigma_r(\bS')\ge\frac12\sigma_r(\bS)
 [1-\tau\kappa(\bS)]\ge\frac14\sigma_r(\bS),\qquad
 \|\bS'\|\le\frac{1+\tau}{2}\|\bS\|.
\]
Together with $\operatorname{col}(\bS')\subseteq\operatorname{col}(\bU)$,
these bounds give $\operatorname{rank}(\bS')=r$ and
$\kappa(\bS')\le2(1+\tau)\kappa(\bS)$.

For a uniform layer half $N_k$, the counts
$n_m'=|\{i\in N_k:z_i^*=m\}|$ are hypergeometric, with
$\mathbb E n_m'=n_m|N_k|/n$ and
$|\mathbb E n_m'-n_m/2|\le1$. Since $n_m\asymp n$, for
$\tau=C\log n/\sqrt n$ the sampling-without-replacement bound gives
\[
 \mathbb P\!\left(\left|\frac{2n_m'}{n_m}-1\right|>\tau\right)
 \le2\exp\left(-\frac{2(n_m\tau/2-1)^2}{|N_k|}\right)
 \le2\exp(-c\log^2 n)
\]
for all sufficiently large $n$. A union bound over both components
and both halves proves the count assertion. Finally,
$\kappa(\bS)\le c\sqrt n/\log n$ gives
$\tau\kappa(\bS)\le Cc\le1/2$ for sufficiently small $c$.
\end{proof}

\begin{proposition}\label{ir:split-preservation}
Suppose the conditions in and preceding Theorem~\ref{thm:init}
hold for the full model, with sufficiently large $C,C_0$ and
sufficiently small $c>0$.
Under uniform node and layer splitting independent of the observations,
with probability at least $1-\delta_{\rm split}$, where
\begin{equation}
 \delta_{\rm split}:=
 Cr\exp(-cd/(\kappa^2r))+C\exp(-c\log^2n),
 \label{ir:split-failure}
\end{equation}
all four restricted models satisfy the same conditions with
\[
 (\alpha,\beta,\kappa)\mapsto(2\alpha,2\beta,3\kappa)
\]
and fixed changes in the remaining constants. Moreover,
\begin{equation}\label{ir:split-eta-bound}
 \eta_z\le
 \frac{Cd\bar p}
 {\sigma_2^2(\scrM_3(\mathbb E\bcalX))}.
\end{equation}
\end{proposition}

\begin{proof}
Write $\bS=n_1\bP_1+n_2\bP_2$ and
$\|\bU\|_{2,\infty}\le\mu\sqrt{r/d}$ for the fixed incoherence
constant $\mu$. Apply Lemma~\ref{ir:vertex-split} with $\eta=1/2$
to both vertex halves, and Lemma~\ref{ir:layer-split} to both layer
halves. A union bound gives, except with probability
$\delta_{\rm split}$, simultaneously for every $(N,V)$,
\begin{gather*}
 \frac{s}{2d}\bI_r\preceq\bH_V:=\bU_V^\top\bU_V
       \preceq\frac{3s}{2d}\bI_r,\qquad s=|V|,\\
 n_m':=|\{i\in N:z_i^*=m\}|=\frac{n_m}{2}(1+e_m),\qquad
 \max_m|e_m|\le\tau:=C\log n/\sqrt n,\qquad
 \tau\kappa(\bS)\le\frac12.
\end{gather*}
For large $n$, $\tau\le1/4$, so the count bounds give balance
constants $2\alpha$ and $2\beta$.

Lemma~\ref{ir:vertex-split} gives $\operatorname{rank}(\bZ_V)=r$,
$\kappa(\bD_V)\le\sqrt3\kappa$,
$\operatorname{rank}(\bP_m[V,V])=K$, and
$\kappa(\bP_m[V,V])\le3\kappa$.
It also preserves (A1). With
$\bQ_V=\bU_V\bH_V^{-1/2}$, we have
\[
 \|\bQ_V\|_{2,\infty}
 \le\|\bU\|_{2,\infty}\|\bH_V^{-1/2}\|
 \le\sqrt2\mu\sqrt{r/s}.
\]
By orthogonal invariance, the same row bound holds for the left
singular vectors of $\bZ_V$.
Applying the layer lemma and then the vertex lemma to
$\bS_{N,V}:=\sum_m n_m'\bP_m[V,V]$ yields
\[
 \sigma_r(\bS_{N,V})\ge\frac{s}{8d}\sigma_r(\bS),
 \qquad
 \kappa(\bS_{N,V})\le9\kappa(\bS).
\]
Thus (A3) holds on each restriction after adjusting the constants.

For the mode-$3$ singular value, define
\[
 \bF(m,:)=\sqrt{n_m}\,
       \operatorname{vec}(\bU^\top\bP_m\bU)^\top,\qquad m\in[2].
\]
Since $\bU$ has orthonormal columns, the nonzero singular values of
$\scrM_3(\mathbb E\bcalX)$ equal those of $\bF$.
Using $\bQ_V^\top\bQ_V=\bI_r$, the corresponding singular values
after restriction equal those of
\[
 \bF_{N,V}
 =\operatorname{diag}\!\left(\sqrt{n_1'/n_1},\sqrt{n_2'/n_2}\right)
   \bF(\bH_V^{1/2}\otimes\bH_V^{1/2}).
\]
Consequently,
\begin{equation}\label{ir:split-mode3-gap}
 \sigma_2(\scrM_3(\mathbb E\bcalX_{N,V}))
 \ge\frac{s}{4d}\sigma_2(\scrM_3(\mathbb E\bcalX)).
\end{equation}
Let $G(\kappa)$ be the left-hand side of \eqref{eq:init-variable-gap},
and let $G_{N,V}(\kappa)$ denote its restricted version.
The preceding bounds and $|N|\asymp n$, $s\asymp d$ give
\[
 G(3\kappa)\vee\max_{N,V}G_{N,V}(3\kappa)\le CG(\kappa),
 \qquad |N|s\bar p\asymp\bar a.
\]
Thus the gap and growth conditions hold for the full model and all
four restrictions after adjusting the constants.
Finally, \eqref{ir:split-mode3-gap} gives
\eqref{ir:split-eta-bound}.
\end{proof}

The growth conditions give $d/(\kappa^2r)\gtrsim r^2\sqrt{n\log n}$
and hence $\delta_{\rm split}=o(1)$.
Under the additional assumptions of Corollary~\ref{cor:main},
with $\eta_z\to0$ replaced by
\[
 \frac{d\bar p}{\sigma_2^2(\scrM_3(\mathbb E\bcalX))}\to0,
\]
averaging over the split yields, for any $\epsilon\in(0,1)$,
\begin{align*}
 \PP\!\left\{h(\widehat\bz,\bz^*)>e^{-[1-o(1)]I^*}\right\}
 &\le e^{-(I^*)^{1-\epsilon}}+C(n\vee d)^{-3}+\delta_{\rm split},\\
 \EE h(\widehat\bz,\bz^*)
 &\le e^{-[1-o(1)]I^*}+C(n\vee d)^{-3}+\frac{\delta_{\rm split}}2,
\end{align*}
where the second inequality uses $h\le1/2$.

{\color{black}\subsubsection{A matrix concentration inequality}
\begin{lemma}[\citet{brailovskaya2022universality}]
\label{lem:concentration-square}
Let $\bY=\sum_j \bZ_j$ be a sum of independent centered random matrices
with $q_1$ rows and $q_2$ columns. Define
\begin{align*}
 \sigma(\bY)^2&=\max\{\|\EE \bY\bY^\top\|,\|\EE \bY^\top \bY\|\},\\
 \sigma_*(\bY)^2&=\sup_{\|\bu\|=\|\bv\|=1}\ZL{\EE|\bu^\top \bY\bv|^2},\\
 v(\bY)^2&=\|\operatorname{Cov}(\operatorname{vec}\bY)\|,\qquad
 R(\bY)=\max_j\operatorname{ess\,sup}\|\bZ_j\|.
\end{align*}
For an absolute constant $C$ and any $t\ge0$, with probability at least $1-\ZL{C}(q_1+q_2)\exp(-t)$,
\[
 \begin{aligned}
\|\bY\|&\le \ZL{C}\big\{\sigma(\bY)+v(\bY)^{1/2}\sigma(\bY)^{1/2}
 \log^{3/4}(e(q_1+q_2))+\sigma_*(\bY)\sqrt t\\
 &\quad +R(\bY)^{1/3}\sigma(\bY)^{2/3}t^{2/3}+R(\bY)t\big\}.
\end{aligned}
\]
\end{lemma}}
Assume $q_1\le q_2$ and pad $\bY$ and its summands with zero rows:
\[
 \check\bY=\begin{bmatrix}\bY\\0\end{bmatrix},\qquad
 \check\bZ_j=\begin{bmatrix}\bZ_j\\0\end{bmatrix}
 \in\mathbb R^{q_2\times q_2}.
\]
Zero padding leaves all five quantities unchanged:
\[
 \|\check\bY\|=\|\bY\|,\quad
 \sigma(\check\bY)=\sigma(\bY),\quad
 \sigma_*(\check\bY)=\sigma_*(\bY),\quad
 v(\check\bY)=v(\bY),\quad R(\check\bY)=R(\bY).
\]
The square-matrix bound in
\citet[Corollary 2.17 of arXiv:2201.05142v3]{brailovskaya2022universality}
therefore yields the stated bound, since $q_2\le q_1+q_2$.
For $q_1>q_2$, apply the same argument to $\bY^\top$.

{\color{black}For an aggregate over a deterministic set $T$ with $|T|\asymp n$,
decompose $\bX_i=\bX_i^u+\bX_i^l$ for $i\in T$, where
$\bX_i^u(u,v)=\boldsymbol1(u\le v)\bX_i(u,v)$ for $u,v\in[d]$.
Then we have
\[
 \left\|\sum_{i\in T}(\bX_i-\mathbb E\bX_i)\right\|
 \le\left\|\sum_{i\in T}(\bX_i^u-\mathbb E\bX_i^u)\right\|
   +\left\|\sum_{i\in T}(\bX_i^l-\mathbb E\bX_i^l)\right\|.
\]
We first bound the upper-triangular term. Write
$\bar\bX_i^u=\bX_i^u-\mathbb E\bX_i^u$ for $i\in T$. Notice that
\[
 \bY=\sum_{i\in T}\bar\bX_i^u
 =\sum_{i\in T}\sum_{u\le v}\bar\bX_i^u(u,v)\be_u\be_v^\top
 =\sum_{i\in T}\sum_{u\le v}\bZ_{iuv},
\]
where $\bZ_{iuv}=\bar\bX_i^u(u,v)\be_u\be_v^\top$. Put
$D_{uv}=\sum_{i\in T}\operatorname{Var}(\bX_i(u,v))\le C|T|\bar p$ for $u\le v$.
Independence gives
\begin{align*}
 \sigma(\bY)^2
 &=\max\left\{\max_v\sum_{u\le v}D_{uv},
              \max_u\sum_{v\ge u}D_{uv}\right\}\le Cnd\bar p,\\
 v(\bY)^2&=\max_{u\le v}D_{uv}\le Cn\bar p,\\
 \sigma_*(\bY)&\le v(\bY),\qquad R(\bY)\le1.
\end{align*}
The bound for the strict lower-triangular part is almost identical
except that the diagonal entries are zero.
Taking $t=5\log(e(n\vee d))$ in the lemma yields
\begin{equation}\label{eq:aggregate-concentration-new}
 \left\|\sum_{i\in T}(\bX_i-\EE \bX_i)\right\|
 \le C\sqrt{nd\bar p}
\end{equation}
with failure at most $C(n\vee d)^{-3}$ whenever
$d\ge C\log^3(e(n\vee d))$ and $nd\bar p\ge C\log^4(e(n\vee d))$.
Indeed the four additional terms, divided by $\sqrt{nd\bar p}$,
are bounded respectively by constants times
$d^{-1/4}\log^{3/4}(e(n\vee d))$,
$d^{-1/2}\log^{1/2}(e(n\vee d))$,
$(nd\bar p)^{-1/6}\log^{2/3}(e(n\vee d))$, and
$(nd\bar p)^{-1/2}\log(e(n\vee d))$.}
{\color{black}\subsubsection{Concentration of projected noise}}

\begin{lemma}
\label{ir:covariance}
{\color{black}Let $\bA$ be an $n$-by-$d^2$ matrix of independent centered upper-triangular
Bernoulli entries, with zeros in the unused positions and variances at
most $\bar p$.  Let $\bU\in\mathbb R^{d\times r}$ have orthonormal columns. Set $\bB=\bU\otimes \bU$ and $\bY=\bA\bB$.}
Then the parameters in Lemma~\ref{lem:concentration-square} satisfy
\[
 v(\bY)\le\sqrt{\bar p},\qquad
 \sigma_*(\bY)\le\sqrt{\bar p},\qquad
 \sigma(\bY)\le\sqrt{\max\{n,r^2\}\bar p},\qquad
 R(\bY)\le\|\bU\|_{2,\infty}^2.
\]
\end{lemma}

\begin{proof}
{\color{black}For independent Bernoulli variables $\xi_{iuv}$, write
\[
 \bA_{i,u+(v-1)d}=\xi_{iuv}-\mathbb E\xi_{iuv}\quad(u\le v),
 \qquad \bA_{i,u+(v-1)d}=0\quad(u>v).
\]
Then}
\[
 \bY=\sum_{i=1}^n\sum_{\ell=1}^{d^2}\bZ_{i\ell},\qquad
 \bZ_{i\ell}=\bA_{i\ell}\be_i \bb_\ell^\top,\qquad
 \bb_\ell=\bB_{\ell,:}^\top\in\mathbb R^{r^2}.
\]
{\color{black}Since $\bB^\top\bB=\bI_{r^2}$,
\[
 \|\mathbb E\bY^\top \bY\|
 =\|\bB^\top\mathbb E(\bA^\top \bA)\bB\|
 \le\|\mathbb E(\bA^\top \bA)\|\le n\bar p,
\]
where $\mathbb E(\bA^\top\bA)$ is diagonal with entries at most
$n\bar p$. By row independence and centering,
$\mathbb E\bY\bY^\top$ is also diagonal, with}
\[
 \mathbb E(\bY\bY^\top)_{ii}
 =\sum_{\ell=1}^{d^2}\operatorname{Var}(\bA_{i\ell})\|\bb_\ell\|_2^2
 \le\bar p\|\bB\|_F^2=r^2\bar p.
\]
{\color{black}For $i,j\in[n]$ and $a,b\in[r^2]$,
\[
 \operatorname{Cov}(\bY_{ia},\bY_{jb})
 =\boldsymbol1(i=j)(\bB^\top \operatorname{diag}(\operatorname{Var}\bA_{i,:})\bB)_{ab}.
\]
Ordering coordinates by row makes the covariance block diagonal, so}
\[
 v(\bY)^2=\max_{i\in[n]}\|\bB^\top \operatorname{diag}(\operatorname{Var}\bA_{i,:})\bB\|\le\bar p.
\]
{\color{black}For unit vectors $\bu\in\mathbb R^n$ and $\bv\in\mathbb R^{r^2}$,
\[
 \mathbb E|\bu^\top\bY\bv|^2
 =(\bv\otimes\bu)^\top\operatorname{Cov}(\operatorname{vec}\bY)(\bv\otimes\bu)
 \le v(\bY)^2,
\]
so $\sigma_*(\bY)\le v(\bY)$. Finally,}
if $\ell=(a-1)d+b$ for $a,b\in[d]$, then
$\bb_\ell=\bU_{a,:}^\top\otimes \bU_{b,:}^\top$ and
\[
 \|\bZ_{i\ell}\|=|\bA_{i\ell}|\|\bb_\ell\|_2
 \le\|\bU_{a,:}\|_2\|\bU_{b,:}\|_2
 \le\|\bU\|_{2,\infty}^2.
\]
\end{proof}

By Lemma~\ref{lem:concentration-square}, for $n\ge r^2$, with probability at least $1-Cn\exp(-t)$,
\begin{equation}
 \|\bY\|\le C\left[
 \sqrt{n\bar p}
 +n^{1/4}\sqrt{\bar p}(\log n)^{3/4}
 +\sqrt{\bar p\,t}
 +\|\bU\|_{2,\infty}^{2/3}(n\bar p)^{1/3}t^{2/3}
 +\|\bU\|_{2,\infty}^2t\right].
 \label{ir:projected-full-bound}
\end{equation}
{\color{black}The lower-triangular part satisfies the same bound;
a union bound combines the two parts.}

Take $t=4\log(e(n\vee d))$ and
$\|\bU\|_{2,\infty}^2\le C\kappa^2r/d$.
Suppose $n\ge r^2$, $n\ge C\log^3n$, and
$\sqrt{nd^2\bar p}\ge C\kappa^2r\log^2(e(n\vee d))$.
The second, fourth and fifth terms of \eqref{ir:projected-full-bound},
divided by $\sqrt{n\bar p}$, are bounded by
\[
 C\left(\frac{\log^3n}{n}\right)^{1/4},\qquad
 C\left(\frac{\kappa^2r t^2}{\sqrt{nd^2\bar p}}\right)^{1/3},
 \qquad \frac{C\kappa^2r t}{\sqrt{nd^2\bar p}},
\]
respectively, and are bounded by constants under the stated conditions.
For the remaining term, $n=O(d)$ and
$d\ge C\log^3(e(n\vee d))$ imply
\[
 \sqrt{n\bar p}\le C\sqrt{d\bar p},\qquad
 \sqrt{\bar p t}\le C\sqrt{d\bar p}.
\]
Consequently
\begin{equation}\label{ir:projected-failure}
 \|\bY\|\le C\sqrt{d\bar p}
 \quad\text{with probability at least}\quad
 1-C(n\vee d)^{-3}.
\end{equation}
Here $Cn\exp(-t)\le C(n\vee d)^{-3}$.
The growth assumptions of Theorem~\ref{thm:init} give
$n\ge r^2$, $n\ge C\log^3n$, and
\[
 \sqrt{nd^2\bar p}=\sqrt{d\bar a}
 \ge c\kappa^4r^2\sqrt d\log(e\kappa r)
            \log^3(e(n\vee d))
 \ge C\kappa^2r\log^2(n\vee d)
\]
for $C_0$ sufficiently large. The same bounds hold on each half
after adjusting the constants.
{\color{black}Let $\bcalE=\bcalX-\mathbb E\bcalX$, with slices
$\bE_i=\bX_i-\bP_{z_i^*}$.
For $k=1,2,3$, write
\[
 \|\bcalE\|_{k,\delta}
 =\sup_{\|\bu_j\|_2\le1,\,\|\bu_k\|_\infty\le\delta}
 |\langle\bcalE,\bu_1\otimes \bu_2\otimes \bu_3\rangle|,
\]
where $\bu_1,\bu_2\in\mathbb R^d$, $\bu_3\in\mathbb R^n$, and
\[
 \langle\bcalE,\bu_1\otimes \bu_2\otimes \bu_3\rangle
 =\sum_{i=1}^n\sum_{u,v=1}^d \bE_i(u,v)\bu_1(u)\bu_2(v)\bu_3(i).
\]}

{\color{black}\subsubsection{Uniform concentration for reference recovery}}

\begin{lemma}
\label{ra:uniform}
Suppose $nd\bar p\ge C\log^6(e(n\vee d))$, and set
$\delta_U=C\kappa\sqrt{r/d}$ and
$\delta_N=C/\sqrt n$, where the fixed constants are chosen large enough that $\delta_U^2d\ge2$ and $\delta_N^2n\ge2$.
Suppose these thresholds do not exceed one.  With probability at least
$1-C(n\vee d)^{-3}$,
\begin{align}
 \|\bcalE\|_{1,\delta_U}
 &\le C\kappa\sqrt{rd\bar p}\log(e\kappa r)\log^2(e(n\vee d)),
 \label{ra:tensor-u}\\
 \|\bcalE\|_{3,\delta_N}
 &\le C\sqrt{d\bar p}\log^2(e(n\vee d)).
 \label{ra:tensor-n}
\end{align}
Moreover, uniformly over subsets $T\subseteq[n]$ of size at least
$\delta_N^{-2}$,
\begin{equation}
 \left\|\frac1{|T|}\sum_{i\in T}\bE_i\right\|
 \le C\sqrt{\frac{d\bar p}{n}}\log^2(e(n\vee d)).
 \label{ra:whole-class}
\end{equation}
In particular, $T$ may depend on the observations.
\end{lemma}

\begin{proof}
We apply Theorem 1 of \citet{jing2021community}. Since $n=O(d)$,
we may add zero rows and columns to each adjacency matrix so that its
vertex dimension is $n\vee d\asymp d$. This does not change the tensor
norm. The density condition in that theorem follows from
$nd\bar p\ge C\log^6(e(n\vee d))$.
For mode one, its threshold is a constant times
\[
 \max\{1,\delta\sqrt{n\vee d}\log(n\vee d)\}
 \sqrt{(n\vee d)\bar p}
 \log(\delta^2(n\vee d))\log(n\vee d).
\]
For mode three, replace $\delta\sqrt{n\vee d}$ and
$\delta^2(n\vee d)$ by $\delta\sqrt n$ and $\delta^2n$,
respectively. For $\delta=\delta_U$, we have
\[
 \delta_U\sqrt{n\vee d}\le C\kappa\sqrt r,
 \qquad
 \log(\delta_U^2(n\vee d))\le C\log(e\kappa r).
\]
For $\delta=\delta_N$, $\delta_N\sqrt n=C$ and
$\log(\delta_N^2n)=\log(C^2)$. Thus we may take
\[
 t=\begin{cases}
 C\kappa\sqrt{rd\bar p}\log(e\kappa r)\log^2(e(n\vee d)),
                   &k=1,\ \delta=\delta_U,\\
 C\sqrt{d\bar p}\log^2(e(n\vee d)),
                   &k=3,\ \delta=\delta_N.
 \end{cases}
\]
For the two choices of $t$, we have
\[
 \frac t{\delta_U}\ge c d\sqrt{\bar p}\log^2(e(n\vee d)),
 \qquad
 \frac t{\delta_N}\ge c\sqrt{nd\bar p}\log^2(e(n\vee d)).
\]
Since $n=O(d)$ and $\bar p$ is bounded, both choices satisfy
\begin{align*}
 \frac t\delta
 &\ge c\sqrt{nd\bar p}\log^2(e(n\vee d))
 \ge c\log^5(e(n\vee d)),\\
 \frac{t^2}{\bar p}
 &\ge cd\log^4(e(n\vee d))
 \ge c\sqrt{nd\bar p}\log^2(e(n\vee d))
 \ge c\log^5(e(n\vee d)).
\end{align*}
Consequently, the exponential term in the failure probability is bounded by
\[
 C\log^3(e(n\vee d))\exp(-c\log^5(e(n\vee d)))
 \le C(n\vee d)^{-3}.
\]
For the counting bounds in \citet[proof of Theorem~1]{jing2021community},
separate the upper and strict lower triangular entries. Each count is
bounded by $CS$, where $S$ is a sum of independent Bernoulli variables.
For $\mathbb ES\le\mu$,
\[
 \mathbb E\exp(S)\le\exp((e-1)\mu),\qquad
 \mathbb P\big(S\ge e^2\mu+11\log(e(n\vee d))\big)
 \le[e(n\vee d)]^{-11}.
\]
Replace $3\log(n\vee d)$ in the support counts by
$11\log(e(n\vee d))$. For each support family $\mathcal H$,
\[
 |\mathcal H|\exp\big(-\log|\mathcal H|-11\log(e(n\vee d))\big)
 =[e(n\vee d)]^{-11}.
\]
Taking a union bound over at most $3(n\vee d)^2$ fibers and
$C\log^3(e(n\vee d))$ dyadic levels and support families gives
\[
 \big[3(n\vee d)^2+C\log^3(e(n\vee d))\big]
       [e(n\vee d)]^{-11}\le C(n\vee d)^{-3}.
\]
The thresholds change only by constant factors, yielding
\eqref{ra:tensor-u} and \eqref{ra:tensor-n} with probability at least
$1-C(n\vee d)^{-3}$.

{\color{black}For any allowed $T$, put
$w_i=\boldsymbol1(i\in T)/\sqrt{|T|}$, $i\in[n]$.
Then $\|\bw\|_2=1$ and $\|\bw\|_\infty\le\delta_N$.  Therefore,
simultaneously for all such $T$,
\[
 \left\|\frac1{|T|}\sum_{i\in T}\bE_i\right\|
 =\frac1{\sqrt{|T|}}\sup_{\|\bu\|_2=\|\bv\|_2=1}
 |\langle\bcalE,\bu\otimes \bv\otimes \bw\rangle|
 \le\frac{\|\bcalE\|_{3,\delta_N}}{\sqrt{|T|}},
\]
which proves \eqref{ra:whole-class}.}
\end{proof}

{\color{black}\subsubsection{Regularization}}

\begin{lemma}
\label{ra:clipping}
{\color{black}Let $\bU,\bV\in\mathbb R^{d\times r}$ have orthonormal columns, with
$\|\bU\|_{2,\infty}\le\delta$.  Let $\bA$ be obtained by projecting
each row of $\bV$ onto the Euclidean ball of radius $\delta$, and let
$\bQ=\operatorname{SVD}_r(\bA)$.  If}
\[
 \min_{\bO\in\mathbb O_r}\|\bV-\bU\bO\|\le\varepsilon,
 \qquad \sqrt r\varepsilon\le\tfrac18,
\]
then
\[
 \|\bQ\|_{2,\infty}\le\tfrac87\delta,
 \qquad \sigma_{\min}(\bU^\top \bQ)\ge\tfrac12,
 \qquad \min_{\bO\in\mathbb O_r}\|\bQ-\bU\bO\|\le C\sqrt r\varepsilon.
\]
\end{lemma}

\begin{proof}
Choose $\bO\in\arg\min_{\bO'\in\mathbb O_r}\|\bV-\bU\bO'\|$.
By the contraction property of the row projection, we have
\[
 \|\bA-\bU\bO\|\le\|\bA-\bU\bO\|_F
 \le\|\bV-\bU\bO\|_F\le\sqrt r\varepsilon\le\tfrac18.
\]
Weyl's inequality gives $\sigma_r(\bA)\ge1-\|\bA-\bU\bO\|\ge7/8$.
The orthonormal basis $\bA(\bA^\top\bA)^{-1/2}$ therefore satisfies,
for each row,
\[
 \|\bA_{j,:}(\bA^\top \bA)^{-1/2}\|_2
 \le\|\bA_{j,:}\|_2\|(\bA^\top \bA)^{-1/2}\|
 \le\frac{\delta}{\sigma_r(\bA)}\le\frac87\delta.
\]
The same row bound holds for $\bQ$ by orthogonal invariance. Also,
\begin{align*}
 \|(\bI-\bU\bU^\top)\bQ\|
 &\le\frac{\|\bA-\bU\bO\|}{\sigma_r(\bA)}
 \le\tfrac87\sqrt r\varepsilon\le\tfrac17,\\
 \sigma_{\min}(\bU^\top\bQ)^2
 &=1-\|(\bI-\bU\bU^\top)\bQ\|^2\ge1-\tfrac1{49}>\tfrac14.
\end{align*}
Let $\bU^\top\bQ=\bL\operatorname{diag}(\cos\theta_1,\ldots,
\cos\theta_r)\bR^\top$ be a singular value decomposition, where
$\theta_j\in[0,\pi/2]$, and take $\bO=\bL\bR^\top$. Then
\[
 \|\bQ-\bU\bO\|=2\max_j\sin(\theta_j/2)
 \le\sqrt2\max_j\sin\theta_j
 =\sqrt2\|(\bI-\bU\bU^\top)\bQ\|
 \le C\sqrt r\varepsilon.
\]
This completes the proof.
\end{proof}

The aggregate satisfies the preceding accuracy condition if
\begin{equation}
 \sigma_r(\bS)\ge C\sqrt{rnd\bar p},
 \qquad \bS=\sum_i\bP_{z_i^*}.
 \label{ra:aggregate-safe}
\end{equation}
Indeed, on the aggregate-concentration event, the singular-subspace
perturbation bound gives
\[
 \min_\bO\left\|\operatorname{SVD}_r\left(\sum_i\bX_i\right)-\bU\bO\right\|
 \le\frac{C\|\sum_i\bE_i\|}{\sigma_r(\bS)}
 \le\frac1{8\sqrt r},
\]
where the constant in \eqref{ra:aggregate-safe} is sufficiently large.
This gives the condition in Lemma~\ref{ra:clipping}.
{\color{black}\subsubsection{Full-sample reference accuracy and alignment}
With $\delta_1=C_0\sqrt{r/d}$ as in RSpec, its full-sample projection
and layer embedding are
\[
 \bQ=\mathcal P_{\delta_1}\left(\operatorname{SVD}_r\left(\sum_{i=1}^n\bX_i\right)\right),
 \qquad
 \widehat \bW=\operatorname{SVD}_2(\scrM_3(\bcalX)(\bQ\otimes \bQ)).
\]
The reference layer labels $\check \bz:[n]\to[2]$ are the $2$-means
labels from $\widehat \bW$. For each $m$, the reference community labels
$\check\bsig_m:[d]\to[K]$ are the $K$-means labels from
$\operatorname{SVD}_K(\sum_{i:\check z_i=m}\bX_i)$.}

\begin{theorem}
\label{ra:reference}
{\color{black}Assume the balance conditions, $n=O(d)$,
$d\ge C\log^3(e(n\vee d))$, $n/\log n\ge C\kappa^4r^6$, and
\begin{equation}
 nd\bar p\ge C_0\kappa^8r^4\log^2(e\kappa r)
 \log^6(e(n\vee d)).
 \label{ra:original-growth}
\end{equation}}
Suppose (A1)--(A3) and \eqref{eq:init-variable-gap} hold for the full
model. Use exact or fixed-factor approximate $K$-means for RSpec.
Then, with probability at least $1-C(n\vee d)^{-3}$,
\begin{align}
 h(\check\bz,\bz^*)
 &\le\frac{Cd\bar p}{\sigma_2^2(\scrM_3(\bcalC))}
 \left[1+\frac{\kappa^2r^3(nd\bar p)\log^2(e\kappa r)
                      \log^4(e(n\vee d))}
 {\sigma_r^2(n_1\bP_1+n_2\bP_2)}\right],
 \label{ra:reference-layer}\\
 \max_m h(\check\bsig_{\pi(m)},\bsig_m)
 &\le\max_m\frac{C}{\sigma_K(\bP_m)^2}
 \left[\frac{d\bar p}{n}\log^4(e(n\vee d))
             +d^2\bar p^2h(\check\bz,\bz^*)^2\right].
 \label{ra:reference-node}
\end{align}
Here $\pi$ aligns the reference layer labels with the truth and is
also used for the community labels.

The reference labels can be aligned with the local estimates whenever
their combined layer and community errors satisfy
\begin{equation}
 h(\check\bz,\bz^*)+h(\widetilde\bz,\bz^*)\le c/\alpha,
 \qquad
 \max_m\big[h(\check\bsig_{\pi(m)},\bsig_m)
 +h(\widetilde\bsig_{\pi'(m)},\bsig_m)\big]\le c/(\beta K),
 \label{ra:alignment-condition}
\end{equation}
where $c>0$ is sufficiently small and $\pi'$ aligns the local labels
with the truth on the corresponding half. Restricting the reference
labels to that half increases their error bounds by at most a fixed factor.
\end{theorem}

\begin{proof}
Write $\bS=n_1\bP_1+n_2\bP_2$. By
\eqref{eq:aggregate-concentration-new} and Lemma~\ref{ra:clipping},
there exists $\bO\in\mathbb O_r$ such that
\begin{equation}\label{eq:init-reference-projection-error}
 \|\bQ-\bU\bO\|\le\frac{C\sqrt{rnd\bar p}}{\sigma_r(\bS)},
 \qquad \sigma_{\min}(\bU^\top\bQ)\ge\tfrac12,
 \qquad \|\bQ\|_{2,\infty}\le C\kappa\sqrt{r/d},
\end{equation}
with probability at least $1-C(n\vee d)^{-3}$. Since the population
matrix $\bU$ is fixed,
\eqref{ir:projected-failure} and Lemma~\ref{ra:uniform} give
\begin{gather*}
 \|\scrM_3(\bcalE)(\bU\otimes\bU)\|\le C\sqrt{d\bar p},\\
 \|\bcalE\|_{1,\delta_U}
 \le C\kappa\sqrt{rd\bar p}\log(e\kappa r)\log^2(e(n\vee d)),
 \qquad \delta_U=C\kappa\sqrt{r/d},
\end{gather*}
with probability at least $1-C(n\vee d)^{-3}$. Take the constant in $\delta_U$ large enough to cover the row
bounds of both $\bQ$ and $\bU$. The growth conditions imply
$\delta_U\le1$, since
\[
 \frac{\kappa^2r}{d}
 \le\frac{\sqrt{n/\log n}}{\sqrt{C_0}\,r^2d}
 \le\frac{C}{\sqrt{C_0}\,r^2\sqrt{n\log n}}.
\]
The three bounds hold simultaneously with probability at least
$1-C(n\vee d)^{-3}$. On this event, let
$\boldsymbol\Delta=\bQ-\bU\bO$. We have
\[
 \bQ\otimes\bQ-\bU\bO\otimes\bU\bO
 =\boldsymbol\Delta\otimes\bQ+\bU\bO\otimes\boldsymbol\Delta.
\]
By matricization and the bound on matrix rank, we obtain
\begin{align*}
 &\|\scrM_3(\bcalE)(\boldsymbol\Delta\otimes\bQ)\|\\
 &\quad=\sup_{\|\bw\|_2=1}\|\bQ^\top\left(\sum_iw_i\bE_i\right)\boldsymbol\Delta\|_F\\
 &\quad\le\sqrt r\|\boldsymbol\Delta\|
       \sup_{\|\bw\|_2=1}\|\bQ^\top\left(\sum_iw_i\bE_i\right)\|
 \le\sqrt r\|\boldsymbol\Delta\|\,\|\bcalE\|_{1,\delta_U}.
\end{align*}
For the last inequality, if $\bx\in\mathbb R^r$ and $\by\in\mathbb R^d$
are unit vectors, then $\|\bQ\bx\|_2=1$ and
$\|\bQ\bx\|_\infty\le\delta_U$, whence
\[
 |\bx^\top\bQ^\top\left(\sum_iw_i\bE_i\right)\by|
 =|\langle\bcalE,\bQ\bx\otimes\by\otimes\bw\rangle|
 \le\|\bcalE\|_{1,\delta_U}.
\]
Since each $\bE_i$ is symmetric and
$\|\bU\bO\|_{2,\infty}=\|\bU\|_{2,\infty}\le\delta_U$, the other term satisfies
\begin{align*}
 \|\scrM_3(\bcalE)(\bU\bO\otimes\boldsymbol\Delta)\|
 &=\sup_{\|\bw\|_2=1}\|\boldsymbol\Delta^\top\left(\sum_iw_i\bE_i\right)\bU\bO\|_F\\
 &\le\sqrt r\|\boldsymbol\Delta\|
       \sup_{\|\bw\|_2=1}\|\bO^\top\bU^\top\left(\sum_iw_i\bE_i\right)\|
 \le\sqrt r\|\boldsymbol\Delta\|\,\|\bcalE\|_{1,\delta_U}.
\end{align*}
Also, $\bO\otimes\bO$ is orthogonal for every realization of $\bO$, so
\[
 \|\scrM_3(\bcalE)(\bU\bO\otimes\bU\bO)\|
 =\|\scrM_3(\bcalE)(\bU\otimes\bU)(\bO\otimes\bO)\|
 =\|\scrM_3(\bcalE)(\bU\otimes\bU)\|.
\]
Combining these inequalities with
\eqref{eq:init-reference-projection-error} gives
\begin{equation}\label{eq:init-reference-noise}
 \begin{split}
 \|\scrM_3(\bcalE)(\bQ\otimes\bQ)\|
 &\le C\sqrt{d\bar p}\left[1+
 \frac{\kappa r^{3/2}\sqrt{nd\bar p}\log(e\kappa r)
                   \log^2(e(n\vee d))}{\sigma_r(\bS)}\right].
 \end{split}
\end{equation}
Here the tensor norm bounds the two mixed terms uniformly over
$\bQ$ and $\boldsymbol\Delta$ satisfying the row bounds.
{\color{black}For the population term, the orthonormal Tucker
representation implies
\[
 \scrM_3(\mathbb E\bcalX)(\bQ\otimes \bQ)
 =\bW_0\scrM_3(\bcalC)(\bU^\top \bQ\otimes \bU^\top \bQ).
\]
Since $\bW_0$ has orthonormal columns and
$\sigma_{\min}(\bU^\top \bQ)\ge1/2$,
\begin{equation}\label{eq:init-reference-signal}
 \sigma_2(\scrM_3(\mathbb E\bcalX)(\bQ\otimes \bQ))
 \ge\sigma_2(\scrM_3(\bcalC))
       \sigma_{\min}(\bU^\top \bQ)^2
 \ge\frac14\sigma_2(\scrM_3(\bcalC)).
\end{equation}
Thus its left singular space is exactly $\operatorname{col}(\bW_0)$.
For the empirical embedding $\widehat \bW$ defined above,
the Frobenius singular-subspace bound yields
\[
 \min_{\bO\in\mathbb O_2}\|\widehat \bW-\bW_0\bO\|_F^2
 \le C\frac{\|\scrM_3(\bcalE)(\bQ\otimes \bQ)\|^2}
              {\sigma_2(\scrM_3(\bcalC))^2}.
\]
}
By \eqref{eq:init-reference-noise}, \eqref{eq:init-reference-signal},
$(1+x)^2\le2(1+x^2)$, and \eqref{eq:init-variable-gap},
\[
 \frac{\|\scrM_3(\bcalE)(\bQ\otimes\bQ)\|^2}
 {\sigma_2^2(\scrM_3(\mathbb E\bcalX)(\bQ\otimes\bQ))}
 \le\frac{Cd\bar p}{\sigma_2^2(\scrM_3(\bcalC))}
 \left[1+\frac{\kappa^2r^3(nd\bar p)\log^2(e\kappa r)
                         \log^4(e(n\vee d))}{\sigma_r^2(\bS)}\right]
 \le\frac{Cc}{\kappa\sqrt K}.
\]
For sufficiently small $c$, it follows that
\[
 \|\scrM_3(\bcalE)(\bQ\otimes\bQ)\|
 \le\frac12\sigma_2(\scrM_3(\mathbb E\bcalX)(\bQ\otimes\bQ)).
\]

{\color{black}Let $\widehat\bc_1,\widehat\bc_2$ be the fitted centers and choose
$\bO\in\arg\min_{\bO'\in\mathbb O_2}\|\widehat\bW-\bW_0\bO'\|_F$.
The two distinct rows of $\bW_0\bO$ have squared distance $1/n_1+1/n_2$.
Using them as candidate centers in the $K$-means objective gives
\begin{align}\label{eq:init-layer-centers}
 \sum_i\|\widehat \bc_{\check z_i}-(\bW_0\bO)_{i,:}\|_2^2
 &\le2\sum_i\|\widehat \bc_{\check z_i}-\widehat \bW_{i,:}\|_2^2
       +2\|\widehat \bW-\bW_0\bO\|_F^2\notag\\
 &\le C\|\widehat \bW-\bW_0\bO\|_F^2.
\end{align}
Set
\[
 \Delta^2=1/n_1+1/n_2,\qquad
 \mathcal I=\{i\in[n]:\|\widehat\bc_{\check z_i}-(\bW_0\bO)_{i,:}\|_2\ge\Delta/4\}.
\]
By \eqref{eq:init-layer-centers} and \eqref{eq:init-variable-gap},
for sufficiently small $c$,
\[
 |\mathcal I|\le C\|\widehat\bW-\bW_0\bO\|_F^2/\Delta^2
 <\min_m n_m.
\]
Each true class therefore has a row outside $\mathcal I$. Its fitted
center lies within $\Delta/4$ of the true center. These two balls are
disjoint, so the fitted centers define a permutation under which all
rows outside $\mathcal I$ are correctly labelled. Hence
\[
 h(\check\bz,\bz^*)\le\frac{|\mathcal I|}{n}
 \le\frac{C\|\widehat\bW-\bW_0\bO\|_F^2}{n\Delta^2}
 \le C\|\widehat\bW-\bW_0\bO\|_F^2.
\]}
Substituting \eqref{eq:init-reference-noise} and
\eqref{eq:init-reference-signal} into the preceding bound proves
\eqref{ra:reference-layer}.

{\color{black}For the remainder of this componentwise calculation, choose
\[
 \pi\in\underset{\pi'\in\mathfrak S_2}{\arg\min}
             \sum_{i=1}^n\boldsymbol1[\pi'(\check z_i)\ne z_i^*],
\]
and relabel $\check z_i$ as $\pi(\check z_i)$ and
$\check\bsig_m$ as the original $\check\bsig_{\pi^{-1}(m)}$.
Set
$T_m=\{i\in[n]:\check z_i=m\}$. The number of selected layers satisfies
\[
 |T_m|\ge n_m-\sum_i\mathbf1(\check z_i\ne z_i^*)
 \ge\frac{n}{2\alpha}-nh(\check \bz,\bz^*)\ge\frac{n}{4\alpha},
\]
where the last step uses the sufficiently small layer error.}
Choose the constant in $\delta_N$ sufficiently large that
$\delta_N^2n\ge4\alpha$.
The simultaneous bound \eqref{ra:whole-class} then applies to $T_m$.
{\color{black}Moreover,
\begin{align*}
 \frac1{|T_m|}\sum_{i\in T_m}\bP_{z_i^*}-\bP_m
 &=\frac{|\{i\in T_m:z_i^*\ne m\}|}{|T_m|}(\bP_{3-m}-\bP_m),\\
 \left\|\frac1{|T_m|}\sum_{i\in T_m}\bP_{z_i^*}-\bP_m\right\|
 &\le\frac{nh(\check \bz,\bz^*)}{|T_m|}\|\bP_{3-m}-\bP_m\|
 \le Cd\bar p\,h(\check \bz,\bz^*).
\end{align*}
The last inequality follows from the entrywise probability bound and
$|T_m|\asymp n$. Adding the centered observations therefore gives}
\[
 \left\|\frac1{|T_m|}\sum_{i\in T_m}\bX_i-\bP_m\right\|
 \le C\left[\sqrt{\frac{d\bar p}{n}}\log^2(e(n\vee d))
                  +d\bar p\,h(\check\bz,\bz^*)\right].
\]
{\color{black}The empirical community embedding is
\[
 \widehat \bU_m=\operatorname{SVD}_K\left(\sum_{i\in T_m}\bX_i\right)
 \in\mathbb R^{d\times K}.
\]
An orthonormal basis of the population singular space is
$\bZ_m(\bZ_m^\top \bZ_m)^{-1/2}$, and
\begin{align*}
 &\min_\bO\|\widehat \bU_m-\bZ_m(\bZ_m^\top \bZ_m)^{-1/2}\bO\|_F^2\\
 &\qquad\le\frac{CK}{\sigma_K(\bP_m)^2}
       \left\|\frac1{|T_m|}\sum_{i\in T_m}\bX_i-\bP_m\right\|^2.
\end{align*}
where we have used the singular-subspace perturbation bound.

Let $\widehat \bc_{m,1},\ldots,\widehat \bc_{m,K}\in\mathbb R^K$ be
the centers returned by $K$-means on $\widehat \bU_m$, with labels
$\check\bsig_m$. Choose
\[
 \bO_m\in\underset{\bO\in\mathbb O_K}{\arg\min}
       \|\widehat \bU_m-\bZ_m(\bZ_m^\top \bZ_m)^{-1/2}\bO\|_F.
\]
The population row in community $a$ is
$\be_a^\top(\bZ_m^\top \bZ_m)^{-1/2}\bO_m$.
Distinct community rows have squared separation
$|S_k(\bsig_m)|^{-1}+|S_l(\bsig_m)|^{-1}\ge2K/(\beta d)$.
Define the exceptional set by
\[
 \mathcal I_m=\left\{j\in[d]:
 \left\|\widehat \bc_{m,\check\bsig_m(j)}
       -[\bZ_m(\bZ_m^\top \bZ_m)^{-1/2}\bO_m]_{j,:}\right\|_2
                        \ge\frac14\sqrt{\frac{2K}{\beta d}}\right\}.
\]
As in \eqref{eq:init-layer-centers}, compare the fitted centers with
the $K$ population rows. The $K$-means objective gives
\begin{align*}
 &\sum_{j=1}^d
 \left\|\widehat \bc_{m,\check\bsig_m(j)}
       -[\bZ_m(\bZ_m^\top \bZ_m)^{-1/2}\bO_m]_{j,:}\right\|_2^2\\
 &\qquad\le2\sum_{j=1}^d
       \|\widehat \bc_{m,\check\bsig_m(j)}-(\widehat \bU_m)_{j,:}\|_2^2
       +2\|\widehat \bU_m-\bZ_m(\bZ_m^\top \bZ_m)^{-1/2}\bO_m\|_F^2\\
 &\qquad\le C
       \|\widehat \bU_m-\bZ_m(\bZ_m^\top \bZ_m)^{-1/2}\bO_m\|_F^2.
\end{align*}
Every row of $\mathcal I_m$ contributes at least $K/(8\beta d)$
to the sum on the left. Therefore
\begin{equation}\label{eq:init-community-exceptional-count}
 |\mathcal I_m|
 \le\frac{Cd}{K}
 \min_\bO\|\widehat \bU_m-\bZ_m(\bZ_m^\top \bZ_m)^{-1/2}\bO\|_F^2.
\end{equation}
\ZL{
For $a\in[K]$, define
\[
 \mathcal C_{m,a}=\{b\in[K]:\exists j\in[d]\setminus\mathcal I_m,
          \bsig_m(j)=a,\ \check\bsig_m(j)=b\},\qquad
 \mathcal R_m=\{a:\mathcal C_{m,a}\ne\varnothing\}.
\]
The separation of the true centers makes the sets $\mathcal C_{m,a}$
disjoint: a fitted center cannot be within $\sqrt{2K/(\beta d)}/4$
of two true centers. For $a\notin\mathcal R_m$,
$S_a(\bsig_m)\subseteq\mathcal I_m$. Thus, by balance,
\[
 K-|\mathcal R_m|\le\frac{\beta K}{d}|\mathcal I_m|.
\]
Choose $b_m(a)\in\mathcal C_{m,a}$ for $a\in\mathcal R_m$ and put
\[
 \mathcal D_m=\bigcup_{a\in\mathcal R_m}
       (\mathcal C_{m,a}\setminus\{b_m(a)\}),\qquad
 |\mathcal D_m|\le K-|\mathcal R_m|
                 \le\frac{\beta K}{d}|\mathcal I_m|.
\]
Since $|\{a\in\mathcal R_m:|\mathcal C_{m,a}|>1\}|\le|\mathcal D_m|$,
the upper balance bound gives
\[
 |\{j\notin\mathcal I_m:\check\bsig_m(j)\in\mathcal D_m\}|
 \le\frac{\beta d}{K}|\mathcal D_m|\le\beta^2|\mathcal I_m|.
\]
Extend $\pi_m(b_m(a))=a$ to a permutation $\pi_m\in\mathfrak S_K$.
Every row outside $\mathcal I_m$ whose fitted label is not in
$\mathcal D_m$ is then correctly labelled, so
\[
 \{j:\pi_m(\check\bsig_m(j))\ne\bsig_m(j)\}
 \subseteq\mathcal I_m\cup
 \{j\notin\mathcal I_m:\check\bsig_m(j)\in\mathcal D_m\}.
\]
}
Combining this inclusion with
\eqref{eq:init-community-exceptional-count} gives
\begin{align}\label{eq:init-community-rounding}
 h(\check\bsig_m,\bsig_m)
 &\le\frac{(1+\beta^2)|\mathcal I_m|}{d}\notag\\
 &\le\frac{C}{K}
       \min_\bO\|\widehat \bU_m-\bZ_m(\bZ_m^\top \bZ_m)^{-1/2}\bO\|_F^2\notag\\
 &\le\frac{C}{\sigma_K(\bP_m)^2}
       \left\|\frac1{|T_m|}\sum_{i\in T_m}\bX_i-\bP_m\right\|^2.
\end{align}
Combining this with the preceding matrix bound proves
\eqref{ra:reference-node}.}

{\color{black}Finally, let $\bx,\by:D\to[q]$ be two label vectors on
one vertex or layer index set $D$, with truth $\bz:D\to[q]$.
Choose their theoretical permutations as
\[
 \pi_x\in\underset{\pi\in\mathfrak S_q}{\arg\min}
              \sum_{j\in D}\boldsymbol1[\pi(x_j)\ne z_j],\qquad
 \pi_y\in\underset{\pi\in\mathfrak S_q}{\arg\min}
              \sum_{j\in D}\boldsymbol1[\pi(y_j)\ne z_j].
\]
For this comparison write $\bx,\by$ for $\pi_x(\bx),\pi_y(\by)$.
Put
\[
 e=\sum_{j\in D}\mathbf1(x_j\ne z_j)
      +\sum_{j\in D}\mathbf1(y_j\ne z_j),\qquad
 T_{ab}=|\{j\in D:x_j=a,y_j=b\}|.
\]
Then
\[
 T_{aa}\ge|\{j\in D:z_j=a\}|-e,\qquad T_{ab}\le e\quad(a\ne b).
\]
If $e<\frac14\min_a|\{j\in D:z_j=a\}|$, each diagonal entry
strictly exceeds the other entries in its row. The identity is therefore
the unique maximum-overlap permutation.
In the original label coordinates, the maximizing map from $\by$ to $\bx$
is $\pi_x^{-1}\circ\pi_y$. Apply this first to the layer labels and then,}
using the same component permutation,
to the community labels.
Restriction to a half increases normalized errors by at most a fixed
factor, so \eqref{ra:alignment-condition} ensures the required comparisons.
{\color{black}A union bound over the projection and tensor-concentration
events gives the stated probability.}
\end{proof}

\paragraph{Componentwise condition numbers.}
Under the component condition in (A2), for each $m\in[2]$,
\[
 \operatorname{rank}(\bP_m)=K,\qquad
 \sigma_1(\bP_m)/\sigma_K(\bP_m)\le\kappa.
\]
Homogeneity and the entrywise bound imply
\begin{equation}\label{eq:init-component-gap}
 cd\bar p
 \le\frac{\boldsymbol1_d^\top\bP_m\boldsymbol1_d}{d}
 \le\sigma_1(\bP_m)\le\|\bP_m\|_F\le Cd\bar p,
 \qquad
 \sigma_K(\bP_m)\ge\frac{cd\bar p}{\kappa}.
\end{equation}
Consequently,
\begin{equation}\label{eq:init-K-condition}
 \frac{cK d^2\bar p^2}{\kappa^2}
 \le K\sigma_K(\bP_m)^2\le\|\bP_m\|_F^2
 \le Cd^2\bar p^2,
 \qquad K\le C\kappa^2,\quad K\le r\le2K-1.
\end{equation}
The rank bounds were established above. Substituting
\eqref{eq:init-component-gap} into \eqref{ra:reference-node} gives
\begin{equation}\label{eq:init-reference-kappa}
 \max_m h(\check\bsig_{\pi(m)},\bsig_m)
 \le C\kappa^2\left[
 \frac{\log^4(e(n\vee d))}{\bar a}+h(\check\bz,\bz^*)^2\right].
\end{equation}
By \eqref{eq:init-variable-gap} and \eqref{ra:reference-layer},
\[
 \kappa\sqrt K\,h(\check\bz,\bz^*)\le Cc.
\]
Moreover, \eqref{ra:original-growth} and $K\le r$ imply
\[
 \frac{K\kappa^2\log^4(e(n\vee d))}{\bar a}
 \le\frac{CK}{C_0\kappa^6r^4\log^2(e\kappa r)
                              \log^2(e(n\vee d))}\le C/C_0.
\]
Consequently
\begin{equation}\label{eq:init-finite-reference-match}
 K\max_mh(\check\bsig_{\pi(m)},\bsig_m)
 \le C(C_0^{-1}+c^2).
\end{equation}

{\color{black}\subsubsection{Accuracy of the cross-fitted initializer}
For each $k\in\{0,1\}$, take the compact SVD
\[
 \bZ[V_k,:]=\bU_k\bD_k\bR_k^\top,
 \qquad \bU_k\in\mathbb R^{|V_k|\times r},\quad
 \bD_k\in\mathbb R^{r\times r},\quad \bR_k\in\mathbb R^{2K\times r}.
\]
On an admissible split this matrix has rank $r$; $\bU_k,\bR_k$ have
orthonormal columns and $\bD_k$ is positive diagonal. The normalized
layer factor and local population core are
\begin{align*}
 \bW_{0,k}(i,m)
 &=\frac{\boldsymbol1(z_i^*=m)}{\sqrt{|\{j\in N_k:z_j^*=m\}|}},
       \qquad i\in N_k,\ m\in[2],\\
 \bcalC_k
 &=\bar p\bcalB_0\times_1(\bD_k\bR_k^\top)\times_2(\bD_k\bR_k^\top)\\
 &\quad\times_3\operatorname{diag}\left(
           \sqrt{|\{i\in N_k:z_i^*=1\}|},
           \sqrt{|\{i\in N_k:z_i^*=2\}|}\right).
\end{align*}
Thus $\bW_{0,k}\in\mathbb R^{|N_k|\times2}$,
$\bcalC_k\in\mathbb R^{r\times r\times2}$, and
\[
 \mathbb E\bcalX_{N_k,V_k}
 =\bcalC_k\times_1\bU_k\times_2\bU_k\times_3\bW_{0,k}.
\]}
\begin{lemma}
\label{ir:positive}
Under the assumptions of Theorem~\ref{thm:init}, the raw layer labels
in \eqref{ir:crossfit-w} satisfy, simultaneously for $k\in\{0,1\}$,
\begin{equation}\label{ir:positive-z}
 h(\widetilde\bz^{[k]},\bz^{*[k]})
 \le\frac{C|V_k|\bar p}{\sigma_2^2(\scrM_3(\bcalC_k))}
\end{equation}
with probability at least $1-C(n\vee d)^{-3}$, conditional on the
admissible independent split.
\end{lemma}

\begin{proof}
{\color{black}Condition on the split. For each $k\in\{0,1\}$,
the aggregation event is
\[
 \left\{\left\|\sum_{i\in N_{1-k}}
               [\bX_i^{(k)}-\bP_{z_i^*}^{(k)}]\right\|
                    \le C\sqrt{nd\bar p}\right\}.
\]
By \eqref{eq:aggregate-concentration-new} and the singular-subspace
perturbation bound, we obtain
\[
 \min_{\bO\in\mathbb O_r}\|\widehat \bU_k-\bU_k\bO\|
 \le\frac{C\|\sum_{i\in N_{1-k}}[\bX_i^{(k)}-\bP_{z_i^*}^{(k)}]\|}
               {\sigma_r(\sum_{i\in N_{1-k}}\bP_{z_i^*}^{(k)})}
 \le\frac1{8\ZL{\sqrt r}}.
\]
Together with Lemma~\ref{ra:clipping}, this implies that
\[
 \mathcal A_k=\{\sigma_{\min}(\bU_k^\top\widetilde \bU_k)\ge1/2,\quad
           \|\widetilde \bU_k\|_{2,\infty}\le C\kappa\sqrt{r/d}\}
\]
holds with failure probability $C(n\vee d)^{-3}$. For each $k$,
condition on $\bcalX_{N_{1-k},V_k}$. Then $\widetilde\bU_k$ is fixed
and $\bcalX_{N_k,V_k}$ is independent of the conditioning observations.
On $\mathcal A_k$, Lemma~\ref{ir:covariance} and
\eqref{ir:projected-failure}, with the half-sample dimensions, give
\begin{equation}\label{eq:init-conditional-noise}
 \|\scrM_3(\bcalX_{N_k,V_k}-\mathbb E\bcalX_{N_k,V_k})
       (\widetilde \bU_k\otimes\widetilde \bU_k)\|
 \le C\sqrt{d\bar p}.
\end{equation}
On $\mathcal A_k$, the conditional failure probability in
\eqref{eq:init-conditional-noise} is at most $C(n\vee d)^{-3}$.
Averaging over $\bcalX_{N_{1-k},V_k}$ and adding
$\mathbb P(\mathcal A_k^c)$ gives the same total failure bound.

On the other hand, we have
\[
 \scrM_3(\mathbb E\bcalX_{N_k,V_k})
                 (\widetilde \bU_k\otimes\widetilde \bU_k)
 =\bW_{0,k}\scrM_3(\bcalC_k)
                  (\bU_k^\top\widetilde \bU_k\otimes \bU_k^\top\widetilde \bU_k).
\]
It gives
\begin{align*}
 &\sigma_{\ZL{2}}(\scrM_3(\mathbb E\bcalX_{N_k,V_k})
             (\widetilde \bU_k\otimes\widetilde \bU_k))\\
 &\qquad\ge\sigma_2(\scrM_3(\bcalC_k))
              \sigma_{\min}(\bU_k^\top\widetilde \bU_k)^2
 \ge\frac14\sigma_2(\scrM_3(\bcalC_k)).
\end{align*}
Wedin's theorem therefore implies
\[
 \min_\bO\|\widehat \bW_k-\bW_{0,k}\bO\|_F^2
 \le\frac{Cd\bar p}{\sigma_2(\scrM_3(\bcalC_k))^2}.
\]
The two distinct rows of $\bW_{0,k}$ have squared separation equal to
the sum of the reciprocals of the two local class counts, which is
of order $1/n$. Applying the center-counting argument proved above gives}
\begin{equation}\label{eq:init-raw-layer-rate}
 h(\widetilde\bz^{[k]},\bz^{*[k]})
 \le\frac{C|V_k|\bar p}{\sigma_2^2(\scrM_3(\bcalC_k))}
 =\frac{C|V_k|\bar p}
 {\sigma_2^2(\scrM_3(\mathbb E\bcalX_{N_k,V_k}))}.
\end{equation}
The last equality follows from orthonormality of the local Tucker
factors. By \eqref{eq:init-variable-gap} on this restricted model,
\begin{equation}\label{eq:init-raw-smallness}
 \kappa\sqrt K\,h(\widetilde\bz^{[k]},\bz^{*[k]})\le Cc.
\end{equation}
The same condition, with $c$ sufficiently small, ensures the
perturbation and balance requirements in the subspace and two-means
bounds. This proves the lemma.
\end{proof}

}

\ZL{
\subsection{Exact recovery of community memberships}
\label{fs:exact-local}

We use a $K$-means solution within a fixed factor of the optimum,
assigning each row to its nearest fitted center.

Fix one vertex half and relabel its vertices as $[D]$, where $D\asymp d$.
In this subsection, $\bX_i$ and $\bP_m$ denote the $D\times D$ restrictions
of the adjacency and component probability matrices to this half.
The map $\bsig_m:[D]\to[K]$ is the corresponding restricted community
assignment, and its membership matrix is
\[
 \bZ_m(u,l)=\mathbb I[\bsig_m(u)=l],\qquad u\in[D],\ l\in[K].
\]
For a fixed component $m$, write $D_l=|\{u\in[D]:\bsig_m(u)=l\}|$ for its community sizes.
Write $N\asymp n$ for the number of selected layers. By balance and homogeneity,
\[
 cD/K\le D_l\le CD/K\quad(l\in[K]),\qquad
 c\bar p\le\bP_m(u,v)\le C\bar p\quad(u,v\in[D]).
\]
Moreover,
\[
 \operatorname{rank}(\bP_m)=K,\qquad
 \frac{\sigma_1(\bP_m)}{\sigma_K(\bP_m)}\le C\kappa,
 \qquad \sigma_K(\bP_m)\ge cD\bar p/\kappa.
\]
The above bounds follow from (A2) and \eqref{eq:init-component-gap} after restricting to the vertex half. Constants may depend on the fixed balance and homogeneity bounds, but not on $\kappa,K,n,d$.

\subsubsection{Aggregated adjacency matrices}

{\color{black}Fix $m\in\{1,2\}$ and a selected set $T\subseteq[n]$ with $|T|=N$,
of which at least $(1-\epsilon)N$ layers come from component $m$. Define
\[
 N_{-m}=|\{i\in T:z_i^*\ne m\}|.
\]
The set $T$ may be random but is independent of the observations on
this vertex half. In the next two lemmas, condition on $T$ and keep
$\bz^*$ fixed. The matrix
\[
 \bA=\sum_{i\in T}\bX_i
\]
has independent entries on and above the diagonal, each being a sum
of independent Bernoulli variables. Its conditional mean is
\[
 \bA_* = \EE(\bA\mid T)=N \bP_m + \bB,\qquad
 \bB=N_{-m}(\bP_{3-m}-\bP_m),\qquad N_{-m}\le\epsilon N.
\]}
Since $|\bB_{uv}|\le C\epsilon N\bar p$ and
$|(\bA_*)_{uv}|\le CN\bar p$, we have
\begin{equation}\label{fs:bias-norms}
 {\color{black}\|\bB\|\le C\epsilon ND\bar p,}\qquad
 \|\bB\|_{2\to\infty}\le C\epsilon ND\bar p/\sqrt D,
 \qquad \|\bA_*\|_{2\to\infty}\le CND\bar p/\sqrt D.
\end{equation}

{\color{black}Let $\bU_0=\bZ_m(\bZ_m^\top \bZ_m)^{-1/2}$.
Its columns form an orthonormal basis of the nonzero eigenspace of $N\bP_m$, and it has}
\[
 \|\bU_0\|_{2\to\infty}=\max_{l\in[K]}D_l^{-1/2}
 \le C\sqrt{K/D}.
\]
Order the eigenvalues of $\bA_*$ so that
$|\lambda_1(\bA_*)|\ge\cdots\ge|\lambda_D(\bA_*)|$, and choose associated
orthonormal eigenvectors $\bv_1,\ldots,\bv_D\in\mathbb R^D$.
Set
\[
 \bV=(\bv_1,\ldots,\bv_K)\in\mathbb R^{D\times K},\qquad
 \boldsymbol{\Lambda}=\operatorname{diag}(\lambda_1(\bA_*),\ldots,\lambda_K(\bA_*)).
\]
Then $\bV^\top\bV=\bI_K$ and $\bA_*\bV=\bV\boldsymbol{\Lambda}$. By Weyl's inequality, if $\kappa\epsilon\le c$, we have
\[
 \min_{j\le K}|\lambda_j(\bA_*)|\ge cND\bar p/\kappa,
 \qquad \max_{j>K}|\lambda_j(\bA_*)|\le C\epsilon ND\bar p.
\]
The absolute values of the leading eigenvalues are also bounded by $CND\bar p$.

\begin{lemma}\label{fs:det-bias}
There are constants $c,C>0$, independent of $\kappa,K,N,D$, such
that if $\kappa\epsilon\le c$, then
\[
 \min_{\bO\in\mathbb O_K}\|\bV\bO-\bU_0\|_{2\to\infty}
 \le C\kappa\epsilon\sqrt{K/D}.
\]
\end{lemma}
\begin{proof}
Write $\bA_*\bV=\bV\boldsymbol{\Lambda}$. Since $(\bI-\bU_0\bU_0^\top)N\bP_m=0$,
\[
 (\bI-\bU_0\bU_0^\top)\bV=(\bI-\bU_0\bU_0^\top)\bB\bV\boldsymbol{\Lambda}^{-1}.
\]
By \eqref{fs:bias-norms},
\[
 \|(\bI-\bU_0\bU_0^\top)\bB\|_{2\to\infty}
 \le\|\bB\|_{2\to\infty}
       +\|\bU_0\|_{2\to\infty}\|\bU_0^\top \bB\|
 \le C\epsilon ND\bar p\sqrt{K/D}.
\]
Since $\bV$ has orthonormal columns and $\|\boldsymbol{\Lambda}^{-1}\|\le C\kappa/(ND\bar p)$,
\[
 \|(\bI-\bU_0\bU_0^\top)\bV\|_{2\to\infty}
 \le \|(\bI-\bU_0\bU_0^\top)\bB\|_{2\to\infty}\|\boldsymbol{\Lambda}^{-1}\|
 \le C\kappa\epsilon\sqrt{K/D}.
\]
The same identity in operator norm gives
\[
 \|(\bI-\bU_0\bU_0^\top)\bV\|
 \le \|\bB\|\|\boldsymbol{\Lambda}^{-1}\|\le C\kappa\epsilon.
\]
Since $\sigma_{\min}^2(\bU_0^\top\bV)\ge1-C^2\kappa^2\epsilon^2$, we can take
\[
 \bO=(\bU_0^\top \bV)^\top
       [\bU_0^\top \bV\bV^\top \bU_0]^{-1/2}\in\mathbb O_K,
\]
so that
\[
 0\prec\bU_0^\top\bV\bO
 =[\bU_0^\top\bV\bV^\top\bU_0]^{1/2}\preceq\bI_K.
\]
Since
\[
 \bV^\top \bU_0\bU_0^\top \bV
 =\bI-\bV^\top(\bI-\bU_0\bU_0^\top)\bV.
\]
It follows that
\begin{align*}
 \|\bU_0^\top \bV\bO-\bI\|
 &=1-\sigma_{\min}(\bU_0^\top \bV)\\
 &=\frac{1-\sigma_{\min}^2(\bU_0^\top \bV)}
          {1+\sigma_{\min}(\bU_0^\top \bV)}
 \le \|(\bI-\bU_0\bU_0^\top)\bV\|^2\le C\kappa^2\epsilon^2.
\end{align*}
Combining the above bounds for the two orthogonal projections of $\bV\bO-\bU_0$, we obtain
\begin{align*}
 \|\bV\bO-\bU_0\|_{2\to\infty}
 &\le \|(\bI-\bU_0\bU_0^\top)\bV\bO\|_{2\to\infty}
       +\|\bU_0\|_{2\to\infty}\|\bU_0^\top \bV\bO-\bI\|\\
 &\le C\sqrt{\frac KD}(\kappa\epsilon+\kappa^2\epsilon^2)
 \le C\kappa\epsilon\sqrt{\frac KD},
\end{align*}
where the last inequality uses $\kappa\epsilon\le c$. This proves the lemma.
\end{proof}

\subsubsection{Rowwise random error for Bernoulli-sum entries}

\begin{lemma}\label{fs:random-row}
Let $\widehat\bU\in\mathbb R^{D\times K}$ contain orthonormal eigenvectors
corresponding to the $K$ eigenvalues of $\bA$ with largest absolute values.
There are constants $c,C>0$ such that if
\begin{equation}\label{fs:rowwise-smallness}
 \kappa\epsilon\le c,\qquad
 \kappa\left[\sqrt{\frac{\log(e(n\vee d))}{ND\bar p}}
       +\frac{\sqrt K\log(e(n\vee d))}{ND\bar p}\right]\le c,
\end{equation}
then, with conditional probability at least
$1-C(n\vee d)^{-3}$,
\begin{equation}\label{fs:rowwise-rate}
\begin{split}
 &\min_{\bO\in\mathbb O_K}\|\widehat\bU\bO-\bU_0\|_{2\to\infty}\\
 &\quad\le C\sqrt{\frac KD}\bigg[\kappa\epsilon
 +\left(\kappa+\frac{\kappa^2}{\sqrt K}\right)\\
 &\qquad\qquad\times\left(\sqrt{\frac{\log(e(n\vee d))}{ND\bar p}}
       +\frac{\sqrt K\log(e(n\vee d))}{ND\bar p}\right)\bigg].
\end{split}
\end{equation}
The constants are uniform over all sets $T$ with the stated size and contamination fraction.
\end{lemma}
\begin{proof}
Let $\bE=\bA-\bA_*$. For an edge $u<v$ in layer $i$, its contribution to $\bE$
is
\[
 [\bX_i(u,v)-\EE \bX_i(u,v)](\be_u\be_v^\top+\be_v\be_u^\top),
\]
where $\be_u$ is the $u$th coordinate vector in $\mathbb R^D$. For $u=v$,
the corresponding contribution is $[\bX_i(u,u)-\EE \bX_i(u,u)]\be_u\be_u^\top$.
The summands are independent, centered, and bounded in norm by one.
Their variance matrix is diagonal, with norm at most
\[
 \max_{u\in[D]}\sum_{i\in T}\sum_{v=1}^D
       \operatorname{Var}(\bX_i(u,v))
 \le CND\bar p.
\]
Matrix Bernstein's inequality therefore yields, for every $t>0$,
\[
 \PP(\|\bE\|>t\mid T)
 \le 2D\exp\left(-\frac{t^2}{CND\bar p+Ct}\right).
\]
Taking $t=C\{\sqrt{ND\bar p\log(e(n\vee d))}+\log(e(n\vee d))\}$
with a sufficiently large $C$, and using $D\le n\vee d$, gives
\begin{equation}\label{fs:op-bernstein}
 \|\bE\|\le C\{\sqrt{ND\bar p\log(e(n\vee d))}+\log(e(n\vee d))\}
\end{equation}
except on an event of probability at most $(n\vee d)^{-4}$.

For fixed $\bW\in\mathbb R^{D\times q}$, $q\le K$, and $j\in[D]$, $l\in[q]$,
\[
 (\bE_{j,:}\bW)_l
 =\sum_{i\in T}\sum_{v=1}^D
       [\bX_i(j,v)-\EE \bX_i(j,v)]\bW_{vl},
 \qquad
 \sum_{i\in T}\sum_{v=1}^D
       \operatorname{Var}(\bX_i(j,v))\bW_{vl}^2
 \le CN\bar p\|\bW_{:,l}\|_2^2.
\]
Bernstein's inequality and a union bound over $l\in[q]$ give
\begin{equation}\label{fs:row-bernstein}
 \|\bE_{j,:}\bW\|_2
 \le C\left\{\sqrt{N\bar p\log(e(n\vee d))}\,\|\bW\|_F
                    +\sqrt K\log(e(n\vee d))\|\bW\|_{2\to\infty}\right\}
\end{equation}
with probability at least $1-(n\vee d)^{-5}$, uniformly over fixed $\bW$.

Let $\widehat\bU^{[j]}$ be defined as $\widehat\bU$ with $\bA$ replaced by
\[
 \bA^{[j]}:=\bA-\be_j\bE_{j,:}-\bE_{:,j}\be_j^\top
                   +\bE_{jj}\be_j\be_j^\top.
\]
Then
\[
 \bW_j:=(\bI-\bV\bV^\top)\widehat\bU^{[j]}
 \ \perp\!\!\!\perp\ \bE_{j,:}\mid T.
\]
By \eqref{fs:row-bernstein} and a union bound, with conditional probability
at least $1-C(n\vee d)^{-3}$, \eqref{fs:op-bernstein} and the following hold:
\begin{align*}
 \|\bE_{j,:}\bW_j\|_2
 &\le C\left\{\sqrt{N\bar p\log(e(n\vee d))}\|\bW_j\|_F
             +\sqrt K\log(e(n\vee d))\|\bW_j\|_{2\to\infty}\right\},\quad j\in[D],\\
 \|\bE\bV\|_{2\to\infty}
 &\le C\sqrt{\frac KD}
       \left(\sqrt{ND\bar p\log(e(n\vee d))}
                   +\sqrt K\log(e(n\vee d))\right),
\end{align*}
where $\|\bV\|_F=\sqrt K$ and $\|\bV\|_{2\to\infty}\le C\sqrt{K/D}$.

Let $\bV_\perp=(\bv_{K+1},\ldots,\bv_D)$ and
$\boldsymbol\Lambda_\perp=\operatorname{diag}
(\lambda_{K+1}(\bA_*),\ldots,\lambda_D(\bA_*))$.
The population eigenvalue bounds give
\[
 \Delta:=|\lambda_K(\bA_*)|-|\lambda_{K+1}(\bA_*)|
 \ge\frac{cND\bar p}{\kappa},\qquad
 \|\boldsymbol\Lambda_\perp\|\le C\epsilon ND\bar p\le c\Delta.
\]
Moreover, since $\bA_*\bV_\perp=\bV_\perp\boldsymbol\Lambda_\perp$,
\begin{align*}
 \|\bV_\perp\boldsymbol\Lambda_\perp\|_{2\to\infty}
 &\le\|\bA_*\|_{2\to\infty}\le\frac{CND\bar p}{\sqrt D},\\
 \|\bV_\perp\boldsymbol\Lambda_\perp\bV_\perp^\top\bE\bV\|_{2\to\infty}
 &\le\frac{CND\bar p}{\sqrt D}\|\bE\|,\qquad
 \|\boldsymbol\Lambda_\perp\bV_\perp^\top\bE\bV\|
 \le C\epsilon ND\bar p\|\bE\|.
\end{align*}
By \eqref{fs:rowwise-smallness} and \eqref{fs:op-bernstein},
\begin{align*}
 &\frac{\|\bE\|+\sqrt K\log(e(n\vee d))}{\Delta}\\
 &\quad\le C\kappa\left[
 \sqrt{\frac{\log(e(n\vee d))}{ND\bar p}}
       +\frac{\sqrt K\log(e(n\vee d))}{ND\bar p}\right]\le c,\\
 &\frac{\sqrt K\log(e(n\vee d))}{|\lambda_K(\bA_*)|}
 +\frac{\left(\sqrt{ND\bar p\log(e(n\vee d))}
                      +\log(e(n\vee d))\right)^2}
       {\Delta|\lambda_K(\bA_*)|}\le c,
\end{align*}
for a sufficiently small constant in \eqref{fs:rowwise-smallness}.
These inequalities give the remaining conditions of the perturbation theorem.
By the Davis--Kahan theorem,
\[
 \|(\bI-\bV\bV^\top)\widehat\bU\|
 \le\frac{C\|\bE\|}{\Delta},\qquad
 \|\bV^\top\bE\bV\|\le\|\bE\|.
\]
Substituting these bounds into Theorem 5 of \citet{tang2025eigenvector},
we obtain
\begin{align*}
 &\min_{\bO\in\mathbb O_K}
       \|\widehat\bU\bO-\bV\|_{2\to\infty}\\
 &\quad\le
 \frac C\Delta\left(\|\bE\bV\|_{2\to\infty}
                         +\|\bV\|_{2\to\infty}\|\bE\|\right)\\
 &\qquad+
 \frac{C\|\bE\|}{\Delta^2}
       \left(\|\bA_*\|_{2\to\infty}
                         +\sqrt{KN\bar p\log(e(n\vee d))}\right).
\end{align*}
Here the terms containing an additional factor
$\|\bE\|/\Delta$, $\|\boldsymbol\Lambda_\perp\|/\Delta$, or
$\sqrt K\log(e(n\vee d))/\Delta$ are absorbed into the displayed bound
by the preceding inequalities. It follows that
\begin{align*}
 &\min_{\bO\in\mathbb O_K}
       \|\widehat\bU\bO-\bV\|_{2\to\infty}\\
 &\quad\le C\sqrt{\frac KD}
       \left(\kappa+\frac{\kappa^2}{\sqrt K}\right)
       \left[\sqrt{\frac{\log(e(n\vee d))}{ND\bar p}}
             +\frac{\sqrt K\log(e(n\vee d))}{ND\bar p}\right]\\
 &\qquad+C\sqrt{\frac KD}\kappa^2
       \left[\sqrt{\frac{\log(e(n\vee d))}{ND\bar p}}
             +\frac{\sqrt K\log(e(n\vee d))}{ND\bar p}\right]^2.
\end{align*}
By \eqref{fs:rowwise-smallness}, the last term is absorbed into the first;
Lemma~\ref{fs:det-bias} then gives \eqref{fs:rowwise-rate}.
The bound is uniform over the permitted sets $T$; averaging gives the
same result when $T$ is random and independent of the aggregated observations.
\end{proof}

\subsubsection{Exact community recovery}

\begin{lemma}\label{fs:kmeans-exact}
Let $\bx_1,\ldots,\bx_D\in\mathbb R^K$ have true assignment
$\bsig:[D]\to[K]$.
Suppose $\boldsymbol{\mu}_1,\ldots,\boldsymbol{\mu}_K\in\mathbb R^K$ satisfy
\[
 \min_{k\ne l}\|\boldsymbol{\mu}_k-\boldsymbol{\mu}_l\|_2
 \ge c_0\sqrt{K/D},\qquad
 \min_{k\in[K]}|\{j:\bsig(j)=k\}|\ge c_1D/K,
\]
where $c_0>0$ and $0<c_1\le1$.
Suppose the $K$-means output has objective at most $C$ times the optimum,
for a fixed $C\ge1$, and assigns each point to a nearest fitted center.
If
\[
 \max_j\|\bx_j-\boldsymbol{\mu}_{\bsig(j)}\|_2\le\delta_D\sqrt{K/D},
 \qquad K\delta_D^2\le\frac{c_1c_0^2}{128C},
\]
then its labels are exactly $\bsig$ up to permutation.
\end{lemma}
\begin{proof}
Write the fitted centers as $\widehat{\boldsymbol{\mu}}_1,\ldots,\widehat{\boldsymbol{\mu}}_K$.
Using the true centers as a feasible choice in the $K$-means objective gives
\[
 \sum_{j=1}^D\min_{l\in[K]}\|\bx_j-\widehat{\boldsymbol{\mu}}_l\|_2^2
 \le C\sum_{j=1}^D\|\bx_j-\boldsymbol{\mu}_{\bsig(j)}\|_2^2
 \le CK\delta_D^2\le c_1c_0^2/128.
\]
Suppose that, for some $k$, no fitted center is within
$c_0\sqrt{K/D}/4$ of $\boldsymbol{\mu}_k$. For every $j$ with $\bsig(j)=k$,
the triangle inequality gives
\[
 \min_{l\in[K]}\|\bx_j-\widehat{\boldsymbol{\mu}}_l\|_2
 \ge \min_{l\in[K]}\|\boldsymbol{\mu}_k-\widehat{\boldsymbol{\mu}}_l\|_2
       -\|\bx_j-\boldsymbol{\mu}_k\|_2
 \ge (c_0/4-\delta_D)\sqrt{K/D}
 \ge \frac{c_0}{8}\sqrt{K/D}
\]
because $K\ge1$, $c_1\le1$, and $C\ge1$ give $\delta_D<c_0/8$.
It follows that
\[
 \sum_{j=1}^D\min_{l\in[K]}\|\bx_j-\widehat{\boldsymbol{\mu}}_l\|_2^2
 \ge\frac{c_1D}{K}\frac{c_0^2K}{64D}
 =\frac{c_1c_0^2}{64},
\]
a contradiction. Each true center therefore has a fitted center within
$c_0\sqrt{K/D}/4$. Since these $K$ balls are disjoint, after relabeling,
\[
 \max_{k\in[K]}\|\widehat{\boldsymbol{\mu}}_k-\boldsymbol{\mu}_k\|_2
 \le\frac{c_0}{4}\sqrt{K/D}.
\]
For a point with $\bsig(j)=k$ and any $l\ne k$, we then have
\begin{align*}
 \|\bx_j-\widehat{\boldsymbol{\mu}}_k\|_2
 &\le\|\bx_j-\boldsymbol{\mu}_k\|_2+\|\boldsymbol{\mu}_k-\widehat{\boldsymbol{\mu}}_k\|_2
 \le(\delta_D+c_0/4)\sqrt{K/D},\\
 \|\bx_j-\widehat{\boldsymbol{\mu}}_l\|_2
 &\ge\|\boldsymbol{\mu}_k-\boldsymbol{\mu}_l\|_2-\|\bx_j-\boldsymbol{\mu}_k\|_2
       -\|\boldsymbol{\mu}_l-\widehat{\boldsymbol{\mu}}_l\|_2
 \ge(3c_0/4-\delta_D)\sqrt{K/D}.
\end{align*}
Since $\delta_D<c_0/8$, these inequalities give
\[
 \|\bx_j-\widehat{\boldsymbol{\mu}}_{\bsig(j)}\|_2
 <\min_{l\ne\bsig(j)}\|\bx_j-\widehat{\boldsymbol{\mu}}_l\|_2,
 \qquad j\in[D].
\]
Thus the nearest-center assignment recovers $\bsig$ after the same permutation.
\end{proof}

{\color{black}For $\bU_0=\bZ_m(\bZ_m^\top \bZ_m)^{-1/2}$, distinct community rows have
squared distance $1/D_k+1/D_l\ge cK/D$.} By
Lemmas~\ref{fs:random-row} and~\ref{fs:kmeans-exact}, the community labels are recovered exactly if
\begin{equation}\label{eq:init-exact-K}
 K\bigg[\kappa\epsilon+
 \left(\kappa+\frac{\kappa^2}{\sqrt K}\right)
 \left(\sqrt{\frac{\log(e(n\vee d))}{ND\bar p}}
       +\frac{\sqrt K\log(e(n\vee d))}{ND\bar p}\right)\bigg]^2
 \le c.
\end{equation}
Here $c>0$ is sufficiently small. This also implies
\eqref{fs:rowwise-smallness}, since
$\sqrt K(\kappa+\kappa^2/\sqrt K)\ge\kappa$.

\subsubsection{Proof of exact community recovery}

Fix $k\in\{0,1\}$ and set $D=|V_{1-k}|$.
{\color{black}The raw labels $\widetilde\bz^{[k]}$ use observations within
$V_k$ and are independent of those within $V_{1-k}$.}
We apply Lemma~\ref{fs:random-row} conditional on these labels.
{\color{black}If their error is at most $\eta_z\le c$, with $c>0$ smaller
than the layer-balance constant, choose
\[
 \pi_k^*\in\operatorname*{argmin}_{\pi\in\mathfrak S_2}
       \sum_{i\in N_k}
           \mathbb I[\pi(\widetilde z_i^{[k]})\ne z_i^*],
\]
where $\pi_k^*$ aligns the raw labels with the truth; $\pi_k$ in the algorithm aligns them with the reference. For each true component $m\in\{1,2\}$,
its selected set for the community fit on $V_{1-k}$ is
\[
 T=\{i\in N_k:\pi_k^*(\widetilde z_i^{[k]})=m\},\qquad
 N=|T|,\qquad N_{-m}=|\{i\in T:z_i^*\ne m\}|.
\]
Since there are at most $|N_k|\eta_z$ incorrectly clustered layers, the balance condition gives
\[
 |T|\ge c|N_k|-|N_k|\eta_z\ge cn,
 \qquad |T|\le |N_k|\le n,
 \qquad
 \frac{N_{-m}}{|T|}\le
 \frac{|N_k|\eta_z}{c|N_k|-|N_k|\eta_z}
 \le C\eta_z.
\]
}

Let $\eta_z$ be the deterministic bound in
\eqref{eq:init-joint-conclusion}. By \eqref{eq:init-raw-layer-rate},
\begin{equation}\label{eq:init-local-error-radius}
 \max_{k\in\{0,1\}}h(\widetilde\bz^{[k]},\bz^{*[k]})
 \le\eta_z=C\max_{k\in\{0,1\}}
 \frac{|V_k|\bar p}
 {\sigma_2^2(\scrM_3(\mathbb E\bcalX_{N_k,V_k}))}.
\end{equation}
The raw layer event has failure $C(n\vee d)^{-3}$, and
\eqref{eq:init-variable-gap} on the two diagonal restrictions gives
$\kappa\sqrt K\eta_z\le Cc$.
By \eqref{fs:rowwise-rate}, it suffices to verify \eqref{eq:init-exact-K}.
Recall that $K\le r\le2K-1$, $D\asymp d$, $N\asymp n$,
and $ND\bar p\asymp\bar a$. The lower bound on $\bar a$ gives
\begin{align*}
 &(\kappa\sqrt K+\kappa^2)
       \sqrt{\frac{\log(e(n\vee d))}{ND\bar p}}\\
 &\quad\le\frac C{\sqrt{C_0}}
 \frac{\sqrt K/(\kappa^3r^2)+1/(\kappa^2r^2)}
      {\log(e\kappa r)\log^{5/2}(e(n\vee d))}
 \le\frac C{\sqrt{C_0}\log^{5/2}(e(n\vee d))},\\
 &(\kappa\sqrt K+\kappa^2)
       \frac{\sqrt K\log(e(n\vee d))}{ND\bar p}\\
 &\quad\le\frac C{C_0}
 \frac{K/(\kappa^7r^4)+\sqrt K/(\kappa^6r^4)}
      {\log^2(e\kappa r)\log^5(e(n\vee d))}
 \le\frac C{C_0\log^5(e(n\vee d))},
\end{align*}
where we used $K\le r$ and $\kappa\ge1$.
For the bias term, $\epsilon\le C\eta_z$ gives
\[
 \sqrt K\kappa\epsilon
 \le C\sqrt K\kappa\eta_z\le Cc.
\]
For sufficiently large $C_0$ and small $c$, \eqref{eq:init-exact-K} holds.
By Lemmas~\ref{fs:random-row} and~\ref{fs:kmeans-exact}, on
$\{h(\widetilde\bz^{[k]},\bz^{*[k]})\le\eta_z\}$,
\[
 \PP\left(
 \max_{a\in[2]}h(\widetilde\bsig_a^{[1-k]},
                    \bsig_{\pi_k^*(a)}^{[1-k]})>0
 \,\middle|\,\widetilde\bz^{[k]}\right)
 \le C(n\vee d)^{-3},\qquad k\in\{0,1\}.
\]

Let $\pi$ be the reference permutation in Theorem~\ref{ra:reference}.
The reference bounds, \eqref{eq:init-raw-smallness}, and
\eqref{eq:init-finite-reference-match} give
\begin{align*}
 |N_k|\eta_z+n h(\check\bz,\bz^*)
 &<\frac14\min_m|\{i\in N_k:z_i^*=m\}|,\\
 d\max_m h(\check\bsig_{\pi(m)},\bsig_m)
 &<\frac14\min_{m,l}|\{u\in V_{1-k}:\bsig_m(u)=l\}|.
\end{align*}
The maximum-overlap argument in the proof of Theorem~\ref{ra:reference}
gives $\pi_k=\pi\circ\pi_k^*$ and aligns the two local community labelings
with the same reference labeling. Averaging the conditional bound and
taking a union bound with the layer and reference events gives
\begin{equation}\label{fs:exact-init}
 \frac1n\sum_{i=1}^n\mathbb I[\widetilde z_i\ne\pi(z_i^*)]
 \le C\max_{k\in\{0,1\}}
 \frac{|V_k|\bar p}
 {\sigma_2^2(\scrM_3(\mathbb E\bcalX_{N_k,V_k}))},
 \qquad \max_m h(\widetilde\bsig_{\pi(m)},\bsig_m)=0
\end{equation}
with probability at least $1-C(n\vee d)^{-3}$, conditional on an admissible
split. This proves \eqref{eq:init-joint-conclusion}.

}

\ZL{
\subsection{Proof of Corollary~\ref{cor:main}}
\label{fs:core-separation}

Condition on a split satisfying Theorem~\ref{thm:init}.
By \eqref{fs:exact-init}, after relabeling the components and communities,
\[
 \frac1n\sum_{i=1}^n\mathbb I(\widetilde z_i\ne z_i^*)\le\eta_z,
 \qquad \widetilde\bsig_m=\bsig_m,\quad m\in[2],
\]
with probability at least $1-C(n\vee d)^{-3}$.
We apply Proposition~\ref{fs:exact-community-refinement} with
$\delta_{\rm init}=C(n\vee d)^{-3}$.

For each true edge cell $A_{mk\ell}$, let
\[
 \mu_{\min}=\min_{m\in[2],\,1\le k\le\ell\le K}
             |A_{mk\ell}|\bB_m(k,\ell).
\]
By community balance and homogeneity,
\begin{equation}\label{fs:cor-cell-scale}
 \mu_{\min}\asymp\frac{d^2\bar p}{K^2}
 =\frac dn\frac{\bar a}{K^2}.
\end{equation}
For Bernoulli probabilities $p,q\le1-c_0$,
\[
 D_{\rm KL}({\rm Bern}(p)\|{\rm Bern}(q))
 \le\frac{(p-q)^2}{q(1-q)}.
\]
Since all probabilities are comparable to $\bar p$, summing over the
independent edges gives
\begin{equation}\label{fs:cor-information-upper}
 I^*\le D_{\rm KL}(\PP_1\|\PP_2)
 \le Cd^2\bar p\le CK^2\mu_{\min}.
\end{equation}
Here $\PP_m$ is the distribution of one layer from component $m$.

The gap condition in Theorem~\ref{thm:init} gives
$\kappa\sqrt K\eta_z\le Cc$. Since $K\le C\kappa^2$ by
\eqref{eq:init-K-condition},
\begin{equation}\label{fs:cor-finite-bias}
 K^2\eta_z^2
 \le C\kappa^2K\eta_z^2\le Cc^2.
\end{equation}

By $K\le r$, $\kappa\ge1$, and the dimension and density conditions, we have
\begin{equation}\label{fs:cor-small-ratios}
 \frac{K^4}{n}\le\frac{C}{K^2\log n}\to0,
 \qquad
 \frac{K^4}{\bar a}\le\frac{C}{\log^6(e(n\vee d))}\to0,
 \qquad \log(eK)\lesssim\log n.
\end{equation}
Taking $A=4$ in Proposition~\ref{fs:exact-community-refinement}, we have
\[
 q=\eta_z\log\frac e{\eta_z}
       +\frac{4\log n+4I^*+\log(2K(K+1))}{n},
 \qquad
 R=\eta_z^2I^*+K(K+1)q\left(1+\frac q{\mu_{\min}}\right).
\]
By \eqref{fs:cor-cell-scale}, \eqref{fs:cor-information-upper}, and $n=O(d)$,
\begin{align}
 \frac{K(K+1)q}{\mu_{\min}}
 &\le C\frac{K^4}{\bar a}
       \left(\eta_z\log\frac e{\eta_z}+\frac{\log n}{n}\right)
       +C\frac{K^4}{n}=o(1).
 \label{fs:cor-stochastic-radius}
\end{align}
In particular, $q/\mu_{\min}=o(1)$. Equations
\eqref{fs:cor-information-upper}--\eqref{fs:cor-stochastic-radius} imply
\[
 \frac R{\mu_{\min}}
 \le CK^2\eta_z^2+o(1)\le Cc^2+o(1).
\]
Thus $R\le c'\mu_{\min}$ for the constant $c'$ in
Proposition~\ref{fs:exact-community-refinement}, provided that $c$ is sufficiently small.

Since $q/\mu_{\min}=o(1)$, \eqref{cond:Jmin-main} and
\eqref{fs:cor-small-ratios} give
\begin{align*}
 \frac R{I^*}
 &\le\eta_z^2+C\frac{K^2}{I^*}
       \left(\eta_z\log\frac e{\eta_z}+\frac{\log n}{n}\right)
       +C\frac{K^2}{n}=o(1),\\
 \frac{\sqrt{RK(K+1)}}{I^*}
 &\le C\sqrt{\frac R{I^*}\frac{K^2}{I^*}}=o(1),\\
 \frac{K(K+1)\sqrt R}{I^*\sqrt{\mu_{\min}}}
 &\le C\frac{K^2}{I^*}\sqrt{\frac R{\mu_{\min}}}=o(1).
\end{align*}
Substituting these bounds into \eqref{fs:exact-refinement-expectation} yields
\begin{align*}
 \EE h(\widehat\bz,\bz^*)
 &\le\exp(-[1-o(1)]I^*)+C(n\vee d)^{-3}
                  +n^{-4}\exp(-4I^*)\\
 &\le\exp(-[1-o(1)]I^*)+C(n\vee d)^{-3}.
\end{align*}
where the last inequality follows from $n^{-4}\exp(-4I^*)=o(\exp(-I^*))$.
Taking $u=2(I^*)^{1-\epsilon}$ in \eqref{fs:exact-refinement-probability}, we obtain
\begin{align*}
 &\PP\left(h(\widehat\bz,\bz^*)>
       \exp(-[1-o(1)]I^*+2(I^*)^{1-\epsilon})\right)\\
 &\qquad\le\exp(-2(I^*)^{1-\epsilon})+C(n\vee d)^{-3}
                  +n^{-4}\exp(-4I^*)\\
 &\qquad\le\exp(-(I^*)^{1-\epsilon})+C(n\vee d)^{-3}
\end{align*}
for all sufficiently large $n$. Since $(I^*)^{1-\epsilon}=o(I^*)$,
the proof is complete.

}

\end{document}